\documentclass[10pt]{amsart}

\usepackage{mlmodern}
\usepackage{microtype}

\usepackage[T1]{fontenc}
\usepackage[all,cmtip]{xy}

\usepackage[cal=boondoxo]{mathalfa}
\DeclareSymbolFontAlphabet{\mathcalorig} {symbols}

\usepackage{amssymb}
\usepackage{amsmath}
\usepackage{amsthm}
\usepackage{amssymb}
\usepackage{amsfonts}
\usepackage{latexsym}
\usepackage{float}

\usepackage{amsmath,amssymb,amsthm,mathtools}

\allowdisplaybreaks
\usepackage{mysects2} 
\usepackage{etoolbox}
\patchcmd{\section}{\scshape}{\large\scshape}{}{} 

\usepackage{longtable}
\usepackage{afterpage}
\usepackage{array}
\usepackage{verbatim}
\usepackage{graphicx} 
\usepackage[svgnames,x11names]{xcolor}
\usepackage{fancyhdr} 
\usepackage[labelformat=simple]{caption}
\usepackage{setspace}
\usepackage{bm}

\usepackage{enumitem}

\usepackage{microtype}
\usepackage{lineno}
\usepackage{xcolor}
 \usepackage[normalem]{ulem}
 \newcommand\rsout{\bgroup\markoverwith{\textcolor{red}{\rule[0.5ex]{2pt}{1pt}}}\ULon}

\usepackage{hhline}

\definecolor{refkey}{rgb}{1,1,1}   
\definecolor{labelkey}{rgb}{0,1,0} 
\definecolor{citekey}{rgb}{0,0,0}

\usepackage[multiple]{footmisc}
\usepackage[hidelinks]{hyperref}
\hypersetup{
    colorlinks=true, 
    linktoc=all,     
    linkcolor=DarkBlue,  
	citecolor=Black,
	}

\usepackage{appendix}

\usepackage[backend=biber,doi=false,isbn=false,maxalphanames=4,maxbibnames=4,maxcitenames=99,url=false,eprint=true,style=alphabetic,sorting=nyt]{biblatex}
\renewbibmacro{in:}{}
\AtBeginBibliography{\small}
\renewcommand{\url}[1]{\nolinkurl{#1}}

\newtheorem{maintheorem}{Theorem}

\newtheorem{innercustomthm}{Theorem}
\newenvironment{thm}[1]
  {\renewcommand\theinnercustomthm{#1}\innercustomthm}
  {\endinnercustomthm}

\newtheorem{theorem}{Theorem}[section]

\newtheorem*{theorem*}{Theorem}
\newtheorem{proposition}[theorem]{Proposition}
\newtheorem*{proposition*}{Proposition}
\newtheorem{lemma}[theorem]{Lemma}
\newtheorem*{lemma*}{Lemma}
\newtheorem{corollary}[theorem]{Corollary}

\newtheorem*{question*}{Question}
\theoremstyle{definition}
\newtheorem{definition}[theorem]{Definition}
\newtheorem*{remark*}{Remark}
\newtheorem{remark}[theorem]{Remark}

\newtheorem*{acknowledgements}{Acknowledgements}

\newtheorem{notation}[theorem]{Notation}

\theoremstyle{definition} 
\newtheorem{example}[theorem]{Example}

\newcommand{\be}{\begin{equation}}
\newcommand{\ee}{\end{equation}}

\newcommand{\0}{{(0)}}

\newcommand{\Aut}{\operatorname{Aut}}

\newcommand{\GL}{\operatorname{GL}}

\newcommand{\SL}{\operatorname{SL}}

\newcommand{\SO}{{\operatorname{SO}}}

\newcommand{\Vol}{\operatorname{Vol}}

\newcommand{\meas}{\operatorname{\meas}}

\newcommand{\half}{\frac{1}{2}}

\newcommand{\nc}{\newcommand}

\nc{\la}{\langle} \nc{\ra}{\rangle}
 \nc{\CA}{\cal A}
 \nc{\CBB}{\cal B}

\nc{\CDD}{\cal D}
\nc{\CE}{\cal E}
\nc{\CF}{\cal F} \nc{\CG}{\cal
G} \nc{\CH}{\cal H} \nc{\CI}{\cal I} \nc{\CJ}{\cal J}
\nc{\CK}{\cal K} \nc{\CL}{\cal L} \nc{\CM}{\cal M} \nc{\CN}{\cal
N} \nc{\CO}{\cal O} \nc{\CP}{\cal P} \nc{\CQ}{\cal Q}
\nc{\CR}{\cal R} \nc{\CS}{\cal S} \nc{\CT}{\cal T} \nc{\CU}{\cal
U} \nc{\CV}{\cal V} \nc{\CW}{\cal W} \nc{\CZ}{\cal Z}

\makeatletter
\g@addto@macro\bfseries{\boldmath}
\makeatother

\nc{\fa}{  \mathfrak a} \nc{\fg}{  \Gen} \nc{\fk}{  \mathfrak k}
\nc{\fh}{  \mathfrak h} \nc{\fm}{  \mathfrak m} \nc{\fn}{  \mathfrak n}
\nc{\fA}{  \mathfrak A} \nc{\fC}{  \Cl} \nc{\fI}{  \mathfrak I}
\nc{\fL}{  \mathfrak L} \nc{\fS}{  \mathfrak S}
\nc{\fz}{  \mathfrak z}

\nc{\nen}{\newenvironment} \nc{\ol}{\overline}
\nc{\ul}{\underline} \nc{\lra}{\longrightarrow}
\nc{\lla}{\longleftarrow} \nc{\Lra}{\Longrightarrow}
\nc{\Lla}{\Longleftarrow} \nc{\Llra}{\Longleftrightarrow}
\nc{\hra}{\hookrightarrow} 

\def\C{\mathbb C}
\def\N{\mathbb N}
\def\Q{\mathbb Q}
\def\R{\mathbb R}
\def\Z{\mathbb Z}

\definecolor{deepmagenta}{rgb}{0.8, 0.0, 0.8}
 	\definecolor{earthyellow}{rgb}{0.88, 0.66, 0.37}
 	\definecolor{deepp}{rgb}{0.5, 0.0, 0.13}
    \definecolor{darkgreen}{rgb}{0.13, 0.8, 0.13}
    \definecolor{rpurple}{rgb}{1.02,0.414,0.421 }
    \definecolor{brown}{rgb}{0.788, 0.294, 0.0}
    \definecolor{deepred}{rgb}{0.7, 0.0, 0.0}

\newcommand*{\rom}[1]{\expandafter\@slowromancap\romannumeral #1@}
\newcommand{\RM}[1]{%
  \textup{\uppercase\expandafter{\romannumeral#1}}%
}

\newcommand{\chiz}{\chi_{\raisebox{-.4ex}{$\scriptstyle 0$}}}

\newcommand{\chiorb}{\chi_{\raisebox{-.4ex}{$\scriptstyle \textnormal{orb}$}}}

\DeclareMathOperator{\sg}{sg}
\DeclareMathOperator{\Mat}{Mat}
\DeclareMathOperator{\Ort}{O}

\newcommand{\Cl}{\textup{Cl}}
\newcommand{\Ge}{\textup{Ge}}
\newcommand{\Gen}{\textup{Gen}}
\newcommand{\SGen}{\textup{SGen}}
\newcommand{\Packet}{\textup{Pkt}}

\DeclareMathOperator{\Sq}{Sq}

\DeclareMathOperator{\Supp}{Supp}
\DeclareMathOperator{\Sym}{Sym}

\makeatletter
\@addtoreset{equation}{section}
\makeatother
\numberwithin{equation}{subsection}
 
 \nc{\bq}{\mathbb Q}
 \nc{\br}{\mathbb R}
 \nc{\bz}{\mathbb Z}
 \nc{\bc}{\mathbb C}
 \nc{\bn}{\mathbb N}
 \nc{\bX}{\mathbb{X}}
 \nc{\cL}{ \mathcal{L}}
 \nc{\G}{\Gamma}
 \nc{\sm}{\setminus}
 \nc{\sub}{\subset}
 \nc{\lm}{\lambda}
  \nc{\Lm}{\Lambda}
 \nc{\al}{\alpha}
 \nc{\bt}{\beta}
 \nc{\om}{\omega}
 \nc{\dl}{\delta}
 \nc{\g}{\gamma}
 \nc{\Dl}{\Delta}
 \nc{\Om}{\Omega}
 \nc{\s}{\sigma}
 \nc{\ro}{\rho}
 \nc{\te}{\theta}
 \nc{\SLZ}{\operatorname{SL}_2(\bz)}
 \nc{\SLZp}{\operatorname{SL}_2(\bz[\frac{1}{p}])}
 \nc{\SLZl}{\operatorname{SL}_2(\bz[\frac{1}{\ell}])}
 \nc{\SLR}{\operatorname{SL}_2(\br)}
 \nc{\GLR}{\operatorname{GL}_2(\br)}
 \nc{\PGLR}{\operatorname{PGL}_2(\br)}
 \nc{\PSLR}{\operatorname{PSL}_2(\br)}
 \nc{\PSLZ}{\operatorname{PSL}_2(\bz)}
 
 \nc{\SLC}{\operatorname{SL}_2(\bc)}
 \nc{\uH}{\mathbb H}
 \nc{\fD}{ \mathcal{D}}
 \nc{\fE}{ \mathcal{E}}
 \nc{\fO}{ \mathcal{O}}
 \nc{\haf}{\frac{1}{2}}
 \nc{\qtr}{\frac{1}{4}}
 \nc{\shaf}{{\scriptstyle\frac{1}{2}}}
 \nc{\hlm}{{\scriptstyle\frac{\lambda}{2}}}

 \nc{\8}{\infty}
 \nc{\7}{{-\infty}}
 \nc{\inv}{^{-1}}
 \nc{\eps}{\varepsilon}
 \nc{\aG}{\mathbf{G}}
 \nc{\spn}{\operatorname{Span}}
 \nc{\Cm}{\operatorname{CM}}

  \nc{\SG}{\textup{SG}}

  \nc{\Reg}{\operatorname{Reg}}

 \nc{\MAR}{\operatorname{MAR}}

  \nc{\std}{\operatorname{std}}

 \nc{\Sub}{\operatorname{Sub}}
 \nc{\CSub}{\operatorname{CSub}}
 \nc{\Frob}{\operatorname{Frob}}

 \nc{\sgn}{\operatorname{sgn}}

 \nc{\Per}{\operatorname{Per}}

\nc{\h}{\operatorname{h}}
\nc{\orb}{\textup{orb}}
\nc{\bin}{\textup{bin}}

 \nc{\Stab}{\operatorname{Stab}}

 \nc{\modulo}{\:\textup{mod}\:}

 \nc{\ac}{\operatorname{ac}}

 \newcommand{\track}[1]{}

 \newcommand{\tmod}[1]{\,(\mathrm{mod}\ #1)}

\def\ve {\varepsilon}

\nc{\new}{{\color{red}{\label{ATTENTION ******}}}}

\nc{\upar}{\upshape\arabic*.\ }
\nc{\upal}{itemsep=0.3cm,label=\upshape\alph*.}

\newcounter{daggerfootnote}

\usepackage{chngcntr}
\usepackage{xurl} 
\newcommand{\hide}[1]{}

\makeatletter
\renewcommand{\nxandlist}[4]{%
	\def\and{\space\ensuremath{\mathbf{\cdot}}\space}%
	\let\@@tmp\and
	\def\@@and{\space\ensuremath{\mathbf{\cdot}}\space}%
}
\makeatother

\begin{document}

\renewcommand{\footnoterule}{{\hrule}\vspace{3.5pt}} 
\renewcommand{\thefootnote}{(\arabic{footnote})}

\title[Moments of genus invariants]{Moments of genus invariants and\\$\ell$-adic variations on Siegel's averages}

\author[A. Gamburd]{\scriptsize Alexander Gamburd}
\address{CUNY Graduate Center, NY, NY 10016, USA }
\email{agamburd@gmail.com}

\author[A. Ghosh]{Amit Ghosh}
\address{Oklahoma State University, Stillwater, OK 74078, USA }
\email{ghosh@okstate.edu}

\author[P. Sarnak]{Peter Sarnak}
\address{Institute for Advanced Study \&  Princeton University, Princeton, NJ 08540, USA }
\email{sarnak@math.ias.edu \vspace{.4em}}

\author[J. P. Whang]{Junho Peter Whang}
\address{\parbox{11cm}{Department of Mathematical Sciences \& Research Institute of Mathematics, Seoul National University,
Seoul  08826, S. Korea }}
\email{jhpwhang@gmail.com}


\begin{abstract}
This work concerns the distribution of class numbers and regulators associated to primitive integral ternary quadratic forms.
For every prime $\ell\geq5$, we prove that the orbifold class number of a genus of positive definite forms is typically divisible by arbitrarily large powers of $\ell$, while the prime-to-$\ell$ part is uniformly distributed among the invertible congruence classes modulo $\ell$. The same holds for the Euler characteristics of the arithmetic hyperbolic surfaces associated to the classes of indefinite forms. These results are accompanied by asymptotic formulas for $\ell$-adic moments of genus invariants. Our approach accommodates two ways of ordering genera: by determinant, and by the orbifold class number or regulator. The latter entails a nonabelian analogue of the prime geodesic theorem in lattice space. We also establish versions of these asymptotic formulas for restricted collections of classes, such as those with prime or squarefree determinant.
\end{abstract}

\maketitle

{\footnotesize \tableofcontents}


\setcounter{footnote}{0}

\pagestyle{plain}


\renewcommand{\thefootnote}{(\arabic{footnote}) } 




\section{Introduction}

\subsection{Main results}\

Let $\mathcal F$ denote the set of nonsingular primitive classically integral ternary quadratic forms $F$ of positive determinant $D(F)>0$. Two forms in $\mathcal F$ are in the same class if they are $\GL_3(\Z)$-equivalent. They are in the same genus if they have the same signature and are $\GL_3(\Z_p)$-equivalent for every prime $p$. Each genus decomposes into a finite union of spinor genera (see Def.~\ref{defcat}.4), each of which is a finite union of classes.
For a genus or spinor genus $G $, its class number is $h(G)$. A theorem of Eichler \cite{eichler} states that a spinor genus of indefinite forms in three or more variables has exactly one class.

Given a form $F\in\mathcal F$, let $\Ort_F$ denote the orthogonal group of $F$. If $G$ is a genus or spinor genus of definite forms, we define  its \emph{orbifold class number} $h_{\orb}(G)$  by
\[h_{\orb}(G):=\sum_{i=1}^{h(G)}\frac{2}{|\Ort_{F_i}(\Z)|},\]
where $F_1,\dots,F_{h(G)}$ is a complete set of class representatives in $G$. The name accounts for the fact that each class of definite forms $F$ determines a special point on the orbifold $\SO_3(\R)\backslash\SL_3(\R)/\SL_3(\Z)$, and the weight $2/|\Ort_F(\Z)|$ corresponds to the effective inertia of that orbifold point. Note that if $G$ is a genus we have $h_{\orb}(G)=2m(G)$ where the latter is the mass of $G$ as defined in \cite{Siegel1935}.

If $C$ is a class of indefinite forms represented by $F\in\mathcal F$, then  $\chi_\orb(C)$ will denote the orbifold Euler characteristic of the arithmetic hyperbolic surface
\[\Sigma_F:=(\Ort_3(\R)\cap\Ort_F(\R))\backslash \Ort_F(\R)/\Ort_F(\Z).\]

Given a finite prime $\ell$ and $a\in\Q$, let $|a|_\ell$ denote the $\ell$-adic absolute value of $a$, and write $a_{\ell'}=a|a|_\ell$ for the prime-to-$\ell$ part of $a$. Let $|-|_\infty$ denote the archimedean absolute value on $\Q$. 

Special cases of our first main result are the following. For any prime $5\leq\ell\leq\infty$ and any $-1<\beta\leq\infty$, there exists a constant $c_{\ell,\beta}>0$ such that, as $X\to\infty$,
\begin{enumerate}[label=\upshape\roman*.\ ]
\item\quad
$\displaystyle \sum_{\substack{D(G)\leq X\\\text{$G$ definite}}}|h_{\orb}(G)|_\infty^\beta \sim \sum_{\substack{D(G)\leq X\\\text{$G$ definite}}}|h(G)|_\infty^\beta \sim c_{\infty,\beta}X^{1+\beta}(\log X)^{2^{1-\beta}-1}$, \\
\hspace*{\fill}$(\beta<\infty)$. \vspace{6pt}
    \item\quad
    $\displaystyle \sum_{\substack{D(G)\leq X\\\text{$G$ definite}}}|h_{\orb}(G)|_\ell^\beta\sim c_{\ell,\beta}X(\log X)^{\left(1-\frac{2}{\ell-1}+\frac{2}{\ell^{1+\beta}-1}\right)},\quad(\ell<\infty)$.\vspace{6pt}
    \item\quad
    $\displaystyle \sum_{\substack{D(C)\leq X\\\text{$C$ indef class}}}|\chiorb(C)|_\infty^\beta\sim c_{\infty,\beta}X^{1+\beta}(\log X)^{2^{1-\beta}-1},\quad(\beta<\infty).$\vspace{6pt}
    \item\quad
    $\displaystyle \sum_{\substack{D(C)\leq X\\\text{$C$ indef class}}}|\chiorb(C)|_\ell^\beta\sim c_{\ell,\beta}X(\log X)^{\left(1-\frac{2}{\ell-1}+\frac{2}{\ell^{1+\beta}-1}\right)},\quad(\ell<\infty).$\vspace{6pt}
\end{enumerate}
In (i) and (ii), the sums can be taken to be over genera of definite forms or spinor genera; their asymptotic laws being the same. The cases (i) and (iii) with $\beta=1$ are due to Minkowski \cite{minkowski} and Siegel \cite{siegel}, and read
\begin{align}\label{asymp}
&\sum_{\substack{D(G)\leq X\\\text{$G$ def}}}h_{\orb}(G)\sim\sum_{\substack{D(G)\leq X\\\text{$G$ def}}}h(G)=\sum_{{\substack{D(C)\leq X\\\text{$C$ def class}}}}1\sim \frac{1}{12}\zeta(3) X^2,\\ \intertext{and} 
\label{asymp'}&\sum_{{\substack{D(C)\leq X\\\text{$C$ indef class}}}}|\chiorb(C)|\sim \frac{1}{12}\zeta(3) X^2.
\end{align}
The cases (i) and (iii) with $\beta=0$ combined give
\begin{align}
\label{ggsw-class}\sum_{\substack{D(G)\leq X\\\text{$G$ def}}}1\sim\sum_{\substack{D(C)\leq X\\\text{$C$ indef class}}}1\sim\frac{19}{20}\frac{\zeta(2)}{\zeta(4)} X\log X,    
\end{align}
which was proved in \cite{GGSW}. 

Our $\ell$-adic versions of (\ref{asymp}) and (\ref{asymp'}) show that the orbifold class number $h_{\orb}(G)$ of a definite genus (or spinor genus) $G$ is divisible by large powers of $\ell$ for most $G$'s, and that the orbifold Euler characteristic $\chiorb(C)$ of an indefinite class $C$ is divisible by large powers of $\ell$ for most $C$'s:
\begin{gather*}
    \sum_{\substack{D(G)\leq X\\\text{$G$ def}}}|h_{\orb}(G)|_\ell=o\Bigl(\,\sum_{\substack{D(G)\leq X \\ \text{$G$ def}}}1\,\Bigr),\\\shortintertext{and} \sum_{\substack{D(C)\leq X\\\text{$C$ indef class}}}|\chiorb(C)|_\ell=o\Bigl(\,\sum_{\substack{D(C)\leq X\\\text{$C$ indef class}}}1\,\Bigr).
\end{gather*}

Complementing the above, we  also show that for $\beta>-1/2$ we have
\[\sum_{\substack{D(G)\leq X\\ \text{$G$ def}}}h_{\orb}(G)_{\ell'}^\beta\sim \sum_{\substack{D(C)\leq X\\ \text{$C$ indef class}}}|\chiorb(C)|_{\ell'}^\beta\sim c\cdot X^{1+\beta}(\log X)^{2^{1-\beta}\left(1-\frac{1}{\ell-1}+\frac{1}{\ell^{1+\beta}-1}\right)-1}\]
for some constant $c>0$ depending on $\ell$ and $\beta$ (see Thm.~\ref{thm1n} below with $\beta=\gamma$).

The distribution of invariants for definite forms and indefinite forms can be given a unified treatment. We introduce the notion of measure for an integral quadratic form, which specializes in the ternary case as follows. Given a form $F\in\mathcal F$, the \emph{measure} of $F$ is
\begin{align*}
    \mu(F)=\left\{\begin{array}{ll} 6\cdot 2/|\Ort_F(\Z)| & \text{if $F$ is definite,}\\ 6\cdot |\chiorb(\Sigma_F)| & \text{if $F$ is indefinite.}\end{array}\right.
\end{align*}
The factor $6$ appearing in the definition of the measure is the reciprocal of the \emph{standard mass} of a definite ternary quadratic form defined in \cite[Eq.(7)]{CSMass}. If $G$ is a genus or spinor genus of forms in $\mathcal F$ with class representatives $F_1,\dots,F_{h(G)}$, we define
\[\mu^*(G)=\sum_{i=1}^{h(G)}\mu(F_i).\]
If $G$ is a definite genus, then $\mu^*(G)$ is precisely twice the ratio between Siegel's mass $m(G)$ and the standard mass. The Siegel mass formula expresses $\mu^*(G)$ for a genus $G$ as a product of local densities.  This, along with the observation\footnote{See  Appendix \ref{appendixb} for a sharper statement, and a proof.} from \cite{GGSW} that most spinor genera belong to genera $G$ with only one spinor genera, is the starting point of our present work. Part of what makes the $\ell$-adic averages (i)-(iv) above possible is the elementary nature of the standard mass, which is in marked contrast with the case of binary quadratic forms. (See also Remark \ref{bintern} below.) A formal definition of the measure and a discussion of its relation to class numbers and regulators are given in Section \ref{sect:mass}.

Our first main result concerns the arithmetic distribution of $\mu^*(G)$ for a genus or spinor genus $G$ when ordered by determinant. 

For any $\ell>1$ and $-1<\beta\leq \infty$, define
    \begin{equation}\label{alpha-def}\alpha=\alpha(\ell,\beta) :=1-\frac{1}{\ell-1}+\frac{1}{\ell^{1+\beta}-1}.\end{equation}
Note that $\alpha$ is a decreasing function of $\beta$ with
    \begin{enumerate}[label=\upshape\arabic*.\ ]
        \item $\lim_{\beta\to -1}\alpha(\ell,\beta)=\infty$,
        \item $\alpha(\ell,0)=1$, and
        \item $\alpha(\ell,\infty)=\lim_{\beta\to\infty}\alpha(\ell,\beta)=1-\frac{1}{\ell-1}$.
    \end{enumerate}
As we shall demonstrate (see Sec.~\ref{primedet}), for primes $\ell\geq5$ the quantity $\alpha(\ell,\beta)$ expresses the average order of the $\beta$th power of the $\ell$-adic absolute value of $\mu^*(G)$ among the genera $G$ of prime determinant in $\mathcal F$:
\begin{equation}
\label{alphameaning}
    \alpha(\ell,\beta)=\lim_{X\to\infty}\frac{\sum_{\text{$D(G)\leq X$, prime}}|\mu^*(G)|_\ell^\beta}{\sum_{\text{$D(G)\leq X$, prime}}1}.
\end{equation}
By an unfolding process involving the Landau--Selberg--Delange method and the method of packets used in \cite{GGSW}, we obtain Theorem \ref{thm1n} below, from the above and related asymptotics (See Section \ref{sect:outline} for an outline). 

In what follows, the archimedean Hasse invariant $c_\infty(F)$ of a form $F\in\mathcal F$ is $-1$ if $F$ is definite, and is $+1$ otherwise.

\begin{maintheorem}\label{thm1n}
Let $\ell\geq5$ be prime. Let $\beta,\gamma\in(-1,\infty]$ with
$\min(\beta,0)+\min(\gamma,0)>-1$.
Let $S$ run through the spinor genera of primitive integral ternary quadratic forms of positive determinant, and let $G$ run through the genera of such forms. There is a constant $C_{\ell,\beta,\gamma}>0$, given as an Euler product in \eqref{exp-clb}, such that 
    \[ \begin{split}\sum_{\substack{D(G)\leq X}}|\mu^*(G)|_\ell^\beta|\mu^*(G)|_\infty^\gamma &\sim \sum_{\substack{D(S)\leq X}}|\mu^*(S)|_\ell^\beta|\mu^*(S)|_\infty^\gamma,\\&\sim C_{\ell,\beta,\gamma}~X^{1+\gamma}(\log X)^{2^{1-\gamma}\alpha(\ell,\beta)-1}.\end{split} \]
Moreover for any $\varpi\in\{\pm1\}$ and $(a,\ell)=1$,
    \[\sum_{\substack{c_\infty(S)=\varpi,\,D(S)\leq X\\\mu^*(S)_{\ell'}\equiv a\tmod\ell}}|\mu^*(S)|_\ell^\beta|\mu^*(S)|_\infty^\gamma\sim\frac{1}{2}\frac{1}{\ell -1}~\sum_{\substack{D(S)\leq X}}|\mu^*(S)|_\ell^\beta|\mu^*(S)|_\infty^\gamma,\]
and the same asymptotic holds if one sums over genera $G$. In particular, when ordered by the determininant, the typical $\mu^*(S)$ are divisible by arbitrarily large powers of $\ell$. Moreover, the prime-to-$\ell$ parts of $\mu^*(S)$ are uniformly distributed in $(\Z/\ell\Z)^*$, and the same holds true if we restrict to those $S$ such that $(\mu^*(S),\ell)=1$. 
\end{maintheorem}

\begin{remark}\
\begin{enumerate}[label=\upshape\arabic*.\ ]
    \item For $\beta=\infty$, the main asymptotic formula in Theorem \ref{thm1n} reads:
\[\sum_{\substack{c_\infty(S)=\varpi,\,D(S)\leq X\\m(S)\equiv a\tmod\ell}}|\mu^*(S)|_\infty^\gamma\sim \frac{1}{2}\frac{1}{\ell-1} C_{\ell,\infty,\gamma}~ X^{1+\gamma}(\log X)^{2^{1-\gamma}(1-\frac{1}{\ell-1})-1}.\]
    \item The analogue of Theorem \ref{thm1n} for genera of prime determinants is given in the first part of Theorem \ref{thm-prime}. For squarefree determinants, see Sec.~\ref{sqfree-det}. 
    \item We only consider primes $\ell\geq5$ since, when $\ell= 2$ or $3$, there are additional factors of $2$ and $3$ in $\mu^*$ which change the nature of the calculations. While the sums for these cases can be evaluated, we do not do so.
\end{enumerate}
\end{remark}

A question of similar interest is 
the distribution of $\mu^*$ when the genera and spinor genera are ordered by $\mu^*$. Our method for proving Theorem \ref{thm1n} allows us to treat the distribution of such sums below.

\begin{maintheorem}\label{thm3n}
Let $\ell\geq5$ be prime and $-1<\beta\leq\infty$. Let $S$ run through the spinor genera of primitive integral ternary quadratic forms of positive determinant, and let $G$ run through the genera. There is an explicit constant $C_{\ell,\beta}'>0$ such that 
    \[\sum_{\substack{\mu^*(G)\leq X}}|\mu^*(G)|_\ell^\beta\sim \sum_{\substack{\mu^*(S)\leq X}}|\mu^*(S)|_\ell^\beta\sim C_{\ell,\beta}'~X(\log X)^{8\alpha(\ell,\beta) -1}.\]
Moreover for any $\varpi\in\{\pm1\}$ and $(a,\ell)=1$,
    \[\sum_{\substack{c_\infty(S)=\varpi,\,\mu^*(S)\leq X\\\mu^*(S)_{\ell'}\equiv a\tmod\ell}}|\mu^*(S)|_\ell^\beta\sim\frac{1}{2}\frac{1}{\ell -1}C_{\ell,\beta}'~X(\log X)^{8\alpha(\ell,\beta) -1}.\]
 In particular, when ordered by $\mu^*$, the typical $\mu^*(S)$ are divisible by arbitrarily large powers of $\ell$. Moreover, the prime-to-$\ell$ parts of $\mu^*(S)$ are uniformly distributed in $(\Z/\ell\Z)^*$, and the same holds true if we restrict to those $S$ such that $(\mu^*(S),\ell)=1$.
\end{maintheorem}

The analogue of Theorem \ref{thm3n} where $G$ runs over genera with prime determinant is given in the second part of Theorem \ref{thm-prime}. For squarefree determinants, see Sec.~\ref{sqfree-mu-remark}. 

The following special case of Theorem \ref{thm3n}, with an explicit constant, recovers an analogue of (\ref{ggsw-class}) for counting by $\mu^*$.

\begin{maintheorem}\label{thm2n}
Let $S$ run through the spinor genera of primitive integral indefinite ternary quadratic forms of positive determinant, and let $G$ run through the genera of such forms. For any $\varpi\in\{\pm1\}$,
\[\sum_{\substack{c_\infty(G)=\varpi\\\mu^*(G)\leq X}}1\sim\sum_{\substack{c_\infty(S)=\varpi\\ \mu^*(S)\leq X}}1\sim  \frac{49}{24\Gamma(8)}\frac{1}{\zeta(2)}\prod_{p\geq 3}\left(1-\frac{1}{p}\right)^4\left(1+\frac{4}{p}-\frac{1}{p^2}\right) X(\log X)^7.\]
\end{maintheorem}

Analogues of (i.)-(iv.) of counting by measure are given below (see Theorem \ref{thm4n} in Sec.~\ref{pf-thm4n}). For any prime $5\leq\ell<\infty$ and any $-1<\beta\leq\infty$, there are constants $K_\infty,K_{\ell,\beta}>0$ such that, as $X\to\infty$,
\begin{enumerate}[label=\upshape\roman*$'$.\ ]
    \item\quad
    $\displaystyle    \sum_{\substack{h_{\orb}(G)\leq X\\\text{$G$ definite}}}|h_{\orb}(G)|_\infty^\beta\sim \frac{6}{(1+\beta)}K_\infty X^{1+\beta}(\log X)^{7},\quad (\beta<\infty)$. \vspace{6pt}
    \item\quad
    $\displaystyle \sum_{\substack{h_{\orb}(G)\leq X\\\text{$G$ definite}}}|h_{\orb}(G)|_\ell^\beta\sim K_{\ell,\beta}X(\log X)^{8\left(1-\frac{1}{\ell-1}+\frac{1}{\ell^{1+\beta}-1}\right)-1},\quad(\ell<\infty)$.\vspace{6pt}
    \item\quad
    $\displaystyle \sum_{\substack{|\chiorb(C)|\leq X\\\text{$C$ indef class}}}|h_{\orb}(C)|_\infty^\beta\sim\frac{6}{(1+\beta)}  K_{\infty} X^{1+\beta}(\log X)^7,\quad(\beta<\infty)$.\vspace{6pt}
    \item\quad
    $\displaystyle \sum_{\substack{|\chiorb(C)|\leq X\\\text{$C$ indef class}}}|\chiorb(C)|_\ell^\beta\sim K_{\ell,\beta}X(\log X)^{8\left(1-\frac{1}{\ell-1}+\frac{1}{\ell^{1+\beta}-1}\right)-1},\quad(\ell<\infty).$\vspace{6pt}
\end{enumerate}
Here, we have explicitly
\[K_\infty=\frac{49}{24\Gamma(8)}\frac{1}{\zeta(2)}\prod_{p\geq 3}\left(1-\frac{1}{p}\right)^4\left(1+\frac{4}{p}-\frac{1}{p^2}\right).\]
Results (i$'$.) and (iii$'$.) are obtained from Theorem \ref{thm3n} by simple partial summation. In (i$'$.) and (ii$'$.), the sums can be taken to be over genera of definite forms or spinor genera; their asymptotic laws being the same. As before we draw the conclusion that, when ordered by measure, the orbifold class number $h_{\orb}(G)$ of a definite genus (or spinor genus) $G$ is divisible by large powers of $\ell$ for most $G$'s, and that the orbifold Euler characteristic $\chiorb(C)$ of an indefinite class $C$ is divisible by large powers of $\ell$ for most $C$'s:
\begin{gather*}
\sum_{\substack{h_{\orb}(G)\leq X\\\text{$G$ def}}}|h_{\orb}(G)|_\ell = o\Bigl(\,\sum_{\substack{h_{\orb}(G)\leq X\\\text{$G$ def}}}1\,\Bigr),  \shortintertext{and}
\sum_{\substack{\;|\chiorb(C)|\leq X\\\text{$C$ indef class}}}|\chiorb(C)|_\ell = o\Bigl(\,\sum_{\substack{|\chiorb(C)|\leq X\\\text{$C$ indef class}}}1\,\Bigr).
\end{gather*}
Complementing the above, for $\beta>-1$ we have
\[\sum_{\substack{h_{\orb}(G)\leq X\\ \text{$G$ def}}}h_{\orb}(G)_{\ell'}^\beta\sim \sum_{\substack{|\chiorb(C)|\leq X\\ \text{$C$ indef class}}}|\chiorb(C)|_{\ell'}^\beta\sim c'\cdot X^{1+\beta}(\log X)^{8\left(1-\frac{1}{\ell-1}+\frac{1}{\ell^{1+\beta}-1}\right)-1}\]
for some constant $c>0$ depending on $\ell$ and $\beta$.

In Section \ref{sect:outline}, we outline our proofs of the main theorems and the contents of this paper. In Section \ref{sect:motivation}, we give the history and motivation behind our work.

\subsection{Outline of the proofs}\label{sect:outline}

We give an outline of the proofs of our main theorems and the contents of the paper.

In Section 2, we collect background on integral quadratic forms and their classes, spinor genera, and genera. Motivated by the work of Siegel on his mass formula \cite{Siegel1935, Siegel1936, Siegel1937} and the explicit  mass computations of Conway--Sloane \cite{CSMass}, we introduce notions of regulator $\Reg(F)$ and standard mass $\std(F)$ associated to an integral quadratic form $F$. In brief, $\Reg(F)$ is the volume of the arithmetic orbifold associated to the unit group $\Ort_F(\Z)$. The standard mass $\std(F)$ generalizes that introduced by Conway--Sloane in \cite{CSMass} for positive definite forms. The measure of $F$ is then introduced as the ratio
\[\mu(F)=\frac{\Reg(F)}{\std(F)}.\]
The measure of a form $F$ is always rational. Given a genus $G$ with class representatives $F_1,\dots,F_{h(G)}$, the Siegel mass formula admits formulation as a product formula
\[\mu^*(G)=\sum_{i=1}^{h(G)}\mu(F_i)=\prod_p\mu_p^*(G).\]
Here each $\mu_p^*(G)$, called $p$-measure, is a rational number that depends only on the $p$-adic integral equivalence class of forms determined by $G$, and is invariant under scaling in $\Z_p^*$. Moreover, each $p$-measure satisfies $\mu_p^*(G)=1$ whenever $p$ is an odd prime that does not divide the determinant $D(G)$ of $G$. These properties inspire us the abstract notion of \emph{mass functions} on the space of genera in $\mathcal F$ of which the determinant $D$, the measure $\mu^*$, and the Hasse invariant $c_\infty$, are examples, along with their compositions with multiplicative functions on $\Q^*$. All mass functions admit product formulas of the type above. Thus the problems of asymptotics considered in this paper belong to the following framework: choose a positive mass function $A$ by which to order the genera, and let $B$ be another mass function. We seek asymptotics to
\begin{align}
\label{samplesum}
    \sum_{A(G)\leq X}B(G)
\end{align}
and linear combinations of such sums such as
\[ \sum_{\substack{c_\infty(G)=\varpi\,\\A(G)\leq X}}B(G) = \half\sum_{\substack{A(G)\leq X}}B(G) + \frac{\varpi}{2}\sum_{\substack{A(G)\leq X}}c_\infty(G)B(G).\]
As is customary, we approach (\ref{samplesum}) by studying the analytic properties of the (generalized) Dirichlet series
\[\sum_{G}\frac{B(G)}{A(G)^s}.\]
which admits a formal Euler-like product (see Prop.~\ref{ds-mass}) by virtue of product formulas for mass functions and the fact that we are working with ternary forms. Using this observation, we seek to apply the Landau--Selberg--Delange method (as refined in \cite{GranK}), which gives summatory asymptotic formulas for multiplicative functions whose values at primes admit a finite (possibly non-integral) average.

In Section \ref{sect:preliminary}, we set up certain Eulerian Dirichlet series, supported on squarefree numbers, that model part of our Dirichlet series in the case $A(G)=D(G)$. We verify that their prime coefficients indeed admit a finite average, using the Siegel--Walfitz theorem. This is enough for us to obtain asymptotic formulas for genera with prime determinant in Section \ref{primedet}. The results for genera of squarefree determinants can also be deduced from this. The complete proof of Theorems \ref{thm1n} and \ref{thm3n}, however, encounters three complications:
\begin{itemize}
    \item Part of the preliminary estimates required for \cite{GranK} involves good control on contributions on high ramification, whose verification in our analysis of (\ref{samplesum}) for genera can be computationally involved.
    \item If one sets up a Dirichlet series of the type (\ref{samplesum}) where one sums over spinor genera, one no longer has an Euler-like product.
    \item In the case of Theorem \ref{thm3n} where we are counting by $A(G)=\mu^*(G)$, the denominators one encounters in (\ref{samplesum}) and its Euler-like product expansion are not necessarily natural numbers and prime numbers.
\end{itemize}
We resolve the third issue by applying the method of \cite{BEPS} on generalized Dirichlet series. To address the first two issues, we modify the method of packets used in \cite{GGSW}. In this prior work, root packets were used to partition the collection of all genera of forms. A simplified analysis of the asymptotics for genera for each fixed root packet, and a dominated convergence theorem, allowed us to derive asymptotic formulas such as (\ref{ggsw-class}). In this paper, we perform a similar decomposition, but the definition of root packets and the associated partition are modified so as to depend on the parameters $\beta$ and $\gamma$. This is carried out in Section \ref{pf-thms}.

The method of packets therefore has two roles: it simplifies the treatment of contributions from high ramifications in genus asymptotics, and also allows us to promote the asymptotic formula for genera to one for spinor genera. Below, we give a toy model of the packet method that addresses its first role.

\subsubsection{A toy model for method of packets:}\label{toy}
\quad Suppose that $f:\N\to\C$ is a multiplicative function and let $S_f(X)=\sum_{d\leq X}f(d)$.  For simplicity of exposition, throughout this discussion  we assume that $f$ is real-valued and nonnegative. We describe an unfolding process that derives an asymptotic formula for $S_f(X)$ in stages. Below, $\alpha>0$ is fixed.
\begin{enumerate}[label=\upshape\arabic*.\ ]
    \item For the sum over primes $p$,
    \[\sum_{p\leq X} f(p)\sim\alpha\frac{X}{\log X}.\]
    \item For the sum over squarefree numbers $s$,
    \[\sum_{s\leq X}f(s)\sim\frac{1}{\Gamma(\alpha)}\prod_{p}\left(1+\frac{f(p)}{p}\right)\left(1-\frac{1}{p}\right)^\alpha X(\log X)^{\alpha-1}.\]
    \item For the sum over $d\in\N$,
    \[\sum_{n\leq X}f(d)\sim\frac{1}{\Gamma(\alpha)}\prod_{p}\left(\sum_{k=0}^\infty \frac{f(p^k)}{p^k}\right)\left(1-\frac{1}{p}\right)^\alpha X(\log X)^{\alpha-1}.\]
\end{enumerate}
Assume that $(1)$ holds. More precisely, we assume
\[\sum_{p\leq X}f(p)\log p=\alpha X+O\left(\frac{X}{(\log X)^A}\right),\quad\text{and}\quad \sum_p\sum_{k=2}^\infty\frac{f(p)^k}{p^{k(1-\delta)}}<\infty.\]
Using the Landau--Selberg--Delange method, as refined for multiplicative functions in \cite{GranK}, we observe that (2) follows from the above hypotheses.

To further obtain (3), one needs good control of high ramification meaning the coefficients $f(p^k)$ for $k\geq2$. One sufficient condition, given in \cite[Theorem 1]{GranK}, is that $|f(d)|\leq \tau_k(d)$ for some $k$, where $\zeta(s)^k=\sum_{d=1}^\infty\tau_k(d)/d^s$. Another set of conditions is provided in \cite[Section 7]{GranK}. The method of packets provides a way to derive (3) from (2) under an additional assumption (\ref{toy-ellone}) below: 
\begin{align}
\label{toy-ellone}
\sum_{r\in R}\frac{f(r)}{r}<\infty,\quad\text{or equivalently}\quad \sum_{p}\sum_{k=2}^\infty\frac{f(p^k)}{p^k}<\infty.
\end{align} 
Here $R$ denotes the set of square-full numbers, i.e., those $r\in\N$ such that $v_p(r)\neq 1$ for every $p$. Then every $d\in\N$ can be written uniquely in the form $d=rs$ with $r\in R$ and $s$ squarefree such that $(r,s)=1$.
We remark that the left hand sum is bounded above by the exponential of the right hand sum. 

Now, we  decompose the summatory function $S_f(X)$ as
\[\sum_{d\leq X}f(d)=\sum_{r\in R}f(r)\sum_{\substack{\text{$s\leq X/r$}\\(s,r)=1}}f(s).\]
Using the LSD method, we can derive asymptotic formulas for the inner sums above. Under (\ref{toy-ellone}), we may apply the dominated convergence theorem to justify the interchanging of the limit and sum in the second equality below:
\begin{align*}
&\lim_{X\to\infty}\frac{1}{X(\log X)^{\alpha-1}}\sum_{d\leq X}f(d)=\lim_{X\to\infty}\frac{1}{X(\log X)^{\alpha-1}}\sum_{r\in R}f(r)\sum_{\substack{\text{$s\leq X/r$}\\(s,r)=1}}f(s),\\
&=\sum_{r\in R}f(r)\lim_{X\to\infty}\frac{1}{X(\log X)^{\alpha-1}}\sum_{\substack{\text{$s\leq X/r$}\\(s,r)=1}}f(s),\\
&=\frac{1}{\Gamma(\alpha)}\sum_{r\in R}\frac{f(r)}{r}\prod_{p\mid r}\left(1-\frac{1}{p}\right)^\alpha\prod_{(p,r)=1}\left(1+\frac{f(p)}{p}\right)\left(1-\frac{1}{p}\right)^\alpha,
\end{align*}
which recovers (3).

\subsection{Motivation and history} \label{sect:motivation}
\subsubsection{Binary quadratic forms}\

Let $\Cl_\bin(\Delta)$ denote the set of $\SL_2(\Z)$-equivalence classes of primitive integral binary quadratic forms $ax^2+bxy+cy^2$ of discriminant $\Delta=b^2-4ac$, and $h_\bin(\Delta)=\#\Cl_\bin(\Delta)$ the class number. For the case $\Delta<0$ of definite forms, Gauss posited the asymptotic formula
\[\sum_{-X\leq \Delta<0}h_\bin^+(\Delta)\sim \frac{\pi}{18\zeta(3)} X^{\frac{3}{2}},\]
where $h_\bin^+(\Delta)=h_\bin(\Delta)/2$ is the number of positive definite classes of discriminant $\Delta$. This was proved by Lipschitz \cite{Lipschitz1865} and Mertens \cite{Mertens1874}; and Minkowski \cite{minkowski} later proved an analogous asymptotic formula for definite forms in $m\geq3$ variables. For indefinite forms, observing the apparent need for a weighting of the sum of class numbers, Gauss conjectured
\[\sum_{\substack{0<\Delta\leq X\\\text{$\Delta$ non-square}}}h_\bin(\Delta)R_\Delta\sim \frac{\pi^2}{18\zeta(3)} X^{\frac{3}{2}}\]
where the weight $R_\Delta$ denotes the regulator. This formula was proved by Siegel \cite{siegel}, who established in the same work an analogous formula for quadratic forms in $m\geq3$ variables that recovers Minkowski's result as a special case. For this, Siegel introduced a global measure of a class that replaces the regulator $R_\Delta$, and his asymptotic formula concerns the sum of these measures over the classes $C$ of bounded determinant.

Following these results, a finer arithmetic study of class numbers and regulators was developed for the binary case. We mention some that motivate the present work.

(I) \textit{Cohen--Lenstra heuristics}: Each set $\Cl_\bin(\Delta)$ forms a finite abelian group under Gauss's composition law. Restricting to fundamental discriminants $\Delta$, the Cohen--Lenstra heuristics \cite{CohenLenstra1984I} propose for each odd prime $\ell$ a limiting probability distribution of the $\ell^\infty$-parts $\Cl_\bin(\Delta)[\ell^\infty]$ of the class groups in the space of finite abelian $\ell$-groups. Several partial results are known, mostly focusing on the primes $2$ and $3$; see \cite{Ellenberg2026} for a recent survey.

(II) \textit{Counting by regulator}: The elements of the set $\Cl_\bin(\Delta)$ with non-square $\Delta>0$ are in bijective correspondence with the primitive closed geodesics on the modular surface, with the lengths of the closed geodesics given by $2R_\Delta$. The prime geodesic theorem of Sarnak \cite{sarnak82} for the modular surface gives the asymptotic formula
\[\sum_{\substack{2R_\Delta\leq X\\\text{$\Delta>0$ non-square}}}h_\bin(\Delta)\sim \frac{e^{X}}{X}.\]
Moreover, by proving and using an analogue of the Chebotarev density theorem for primitive closed geodesics on modular surfaces, Sarnak \cite{sarnak82} established the asymptotic counts of indefinite quadratic form classes ordered by regulator with prescribed divisibility properties of coordinates of fundamental units.

(III) \textit{Frequency of class numbers}: Let $\mathcal{D}(h)$ denote the number of fundamental discriminants $\Delta<0$ with $h_2^+(\Delta)=h$. Soundararajan \cite{sound07} established the following asymptotic formula for the average behavior of $\mathcal{D}(h)$: 
    \[\sum_{h\leq X}\mathcal{D}(h)\sim \frac{3\zeta(2)}{\zeta(3)}X^2.\]

It is natural to seek analogues of the theory laid out in the preceding paragraphs for $m$-ary quadratic forms. For forms in three or more variables, it is more natural to group classes by the genus rather than the determinant. In this work we focus on the arithmetic statistics of class numbers and regulators for ternary quadratic forms grouped by genera or spinor genera.

\begin{remark}\label{bintern}
    The following observation describes the differnce between ternary quadratic forms and binary quadratic forms. If $G$ is a definite genus of odd prime determinant $p$ in $\mathcal F$, then by the Siegel mass formula we have
    \[h_{\orb}(G)=\frac{1}{6}\mu_2^*(G)\mu_p^*(G)=\frac{\mu_2^*(G)}{6}\cdot \frac{p+\varepsilon}{2}\] for some $\varepsilon\in\{\pm1\}$, and $\mu_2^*(G)$ is a well-understood product of bounded powers of $2$ and $3$. This is in contrast with, e.g., for $p\equiv 3\modulo 4$,
    \[h_{\bin}^+(-p)=\frac{\sqrt{p}}{2\pi}L(1,\chi_{-p}).\]
    Thus, in our analysis of ternary quadratic forms, we encounter sums over primes $p\leq X$ of expressions such as $\left|\frac{p+\varepsilon}{2}\right|_\ell^\beta$ which can be handled using the Siegel--Walfitz theorem. On the other hand, $|\frac{\sqrt{p}}{2\pi}L(1,\chi_p)|_\ell$ does not admit such a straightforward analysis.
\end{remark}

\vspace{1.5\baselineskip}
\begin{acknowledgements}\leavevmode 

AGa gratefully acknowledges support from a Simons Foundation grant No. 964948.

AGh thanks the Institute for Advanced Study and Princeton University for  visits during parts of the years 2025-26 when some of this work took place.

JPW was supported by the Samsung Science and Technology Foundation under Project Number SSTF-BA2201-03, and by the National Research Foundation of Korea (NRF) grant funded by the Korea government (MSIT) (RS-2025-02293115).

Some of the computations  in the paper were done using Mathematica. 

AI was not used.

\end{acknowledgements}

\section{Quadratic forms and Siegel mass}\label{sect:mass}

\subsection{Regulator and Siegel mass formula}\label{submass1}\

 For $m,n\in\N$, let $\Mat_{m,n}$ denote the affine space of $m\times n$ matrices. We shall write $\Mat_m=\Mat_{m,m}$. Let $\Sym_m$ denote the affine space of symmetric $m\times m$ matrices. If $F\in\Mat_{m,n}$ or $\Sym_m$,  $F_{ij}$ denotes  the $(i,j)$th coordinate of $F$. Let $D:\Mat_m\to \mathbb A^1$ denote the determinant
\[D(F)=\det(F_{ij}).\]
If $A$ is an integral domain and $F\in\Mat_m(A)$, we say $F$ is \emph{nonsingular} if $D(F)\neq 0$.

In this paper, a quadratic form is synonymous with a symmetric matrix. In particular, an integral quadratic form is an element $F\in\Sym_m(\Z)$. We have a morphism
\[\Sym_m\times\Mat_{m,n}\to\Sym_n\]
given by $(F,X)\mapsto F[X]:=X^tFX\in\Sym_n$. The group $\GL_m$ acts on $\Sym_m$ on the right by $F\cdot g=F[g]$. The subgroup $\SL_m$ leaves the determinant $D:\Sym_m\to\mathbb A^1$ invariant.

\begin{definition}
    Let $F\in\Sym_m$ be a quadratic form.
    \begin{enumerate}[label=\upshape\arabic*.\ ]
        \item The group $\Ort_F=\{g\in\GL_m:F[g]=F\}$ is the \emph{orthogonal group} of $F$. The group $\SO_F=\Ort_F\cap\SL_m$ is the \emph{special orthogonal group} of $F$.
        \item The \emph{unit group} of $F\in\Sym_m(\Z)$ is the group $\Ort_F(\Z)$.
        \item Let $R$ be a principal ideal domain. Then $F\in\Sym_m(R)$ is \emph{primitive} if the greatest common divisor of the coefficients of $F$ is $1$. 
    \end{enumerate}    
\end{definition}

Inspired by Siegel's work on his mass formula \cite{Siegel1935, Siegel1936, Siegel1937}, we introduce a convenient notion of \emph{regulator} for integral quadratic forms in some generality. Let us endow the space $\Mat_m(\R)$ of real matrices with the Euclidean measure given by the standard matrix coordinates, and the space $\Sym_m(\R)$ of real symmetric matrices with the Euclidean measure given by the standard upper triangular matrix coordinates. Let $F\in\Sym_m(\R)$ be nonsingular. Siegel's Haar measure $d\tilde g_F$ on $\Ort_F(\R)$ is obtained by disintegration of the Euclidean measure on $\Mat_m(\R)$ by the Euclidean measure on $\Sym_m(\R)$ under the map
\[\Psi_F:\Mat_m(\R)\to\Sym_m(\R)\]
given by $\Psi_F(X)=F[X]:=X^t FX$. If $F\in\Sym_m(\Z)$ is nonsingular, Siegel's invariant $\rho(F)$ is defined by
\[\rho(F)=\Vol(\Ort_F(\R)/\Ort_F(\Z),d\tilde g_F).\]
Explicitly, $\rho(F)$ is defined as follows. Let $T$ be a small neighborhood of $F$ in $\Sym_m(\R)$ and let $Y=\psi_F^{-1}(T)$. Then $\Ort_F(\Z)$ acts on $Y$ by multiplication on the left and stabilizes $Y$. If $Y_0=Y/\Ort_F(\Z)$, then
\[\rho(F)=\lim_{T\to F}\frac{\Vol(Y_0)}{\Vol(T)}\]
where the volumes are normalized Euclidean volumes in their corresponding spaces.

\begin{definition}
Let $F\in\Sym_m(\R)$ be nonsingular. We define the normalized Haar measure on $\Ort_F(\R)$ by
\begin{equation}
\label{haar}
    dg_{F}=|D(F)|^{\frac{m+1}{2}}d\tilde g_F.
\end{equation}
In particular, if $F\in\Sym_m(\Z)$ then $\Vol(\Ort_F(\R)/\Ort_F(\Z),dg_F)=|D(F)|^{\frac{m+1}{2}}\rho(F)$.
\end{definition}
\begin{remark}
     The normalized Haar measure (not Siegel's) has the advantage of being invariant under scaling and under $\GL_m(\R)$-equivalence of quadratic forms. In the case $F=\mathbb I_m\in\Sym_m(\Z)$ is the standard diagonal form, we shall write $dg_m$ for the normalized Haar measure $dg_F=d\tilde g_F$ on $\Ort_m(\R)$ and we have
\[\Vol(\Ort_m(\R))=\prod_{k=1}^m\frac{\pi^{k/2}}{\Gamma(k/2)},\]
which differs from the usual convention by a factor of $2^{m}$.
\end{remark}

\begin{definition}\
\begin{enumerate}[label=\upshape\arabic*.\ ]
    \item The \emph{regulator} of an integral quadratic form $F\in\Sym_m(\Z)$ of signature $(a,b)$ is
\[\Reg(F)=2\frac{\Vol(\Ort_F(\R)/\Ort_F(\Z))}{\Vol(K_F)},\]
where $K_F=\Ort_m(\R)\cap\Ort_F(\R)\sim\Ort_a(\R)\times\Ort_b(\R)$ is a maximal compact subgroup of $\Ort_F(\R)$. Here, the volumes are taken with respect to the normalized Haar measures (\ref{haar}) on $\Ort_F(\R)$ and $K_F$. 
\begin{enumerate}
    \item If $F$ is a definite form, then $K_F=\Ort_F(\R)$ and we have $\Reg(F)=2/|\Ort_F(\Z)|$.
    \item If $F$ is indefinite, then $\Reg(F)$ is the volume of the locally symmetric space
\[\Sigma_F:=K_F\backslash \Ort_F(\R)/\Ort_F(\Z).\]
\end{enumerate} The regulator is a class invariant. If $C$ is the class of $F$, we write $\Reg(C)=\Reg(F)$.
\item If $G$ is a genus of indefinite forms, it turns out that the regulator is constant on $G$ (Prop.~\ref{equal-mass}). We define the \emph{regulator} of an indefinite genus $G$ to be $\Reg(G)=\Reg(F)$ for any choice of $F\in G$.
    \item If $G$ is a genus or a spinor genus of integral quadratic forms, we define the \emph{regulated class number} of $G$ to be
\[h_{\Reg}(G)=\sum_{C\subseteq G}\Reg(C),\]
where the sum runs over the classes $C$ in $G$. We have
\[h_{\Reg}(G)=\left\{\begin{array}{ll} h_{\orb}(G) & \text{$G$ is definite}\\ h(G)\Reg(G) & \text{$G$ is indefinite.}\end{array}\right.\]
\end{enumerate}
\end{definition}

\begin{example}
When $G$ is a genus of signature $(1,m-1)$, the regulator $\Reg(G)$ is the volume of the $(m-1)$-dimensional hyperbolic orbifold $\mathbf H^{m-1}/\operatorname{PO}_F(\Z)$ for any choice of member $F$ in the genus $G$ (See \cite[Section 7]{Siegel1936}). If $m=2n+1$ is odd, then by the Chern--Gauss--Bonnet theorem we have
\[\Reg(G)=\Vol(\Sigma_F)=\frac{2^{2n}\pi^n n!}{(2n)!}|\chiorb(\Sigma_F)|,\]
where $\chiorb(\Sigma_F)$ is the orbifold Euler characteristic of $\Sigma_F=K_F\backslash \Ort_F(\R)/\Ort_F(\Z)$.
In particular, $\Reg(G)$ is a rational multiple of $\pi^n$.
\end{example}

Let $p$ be prime. We write $\rho_p$ for the standard $p$-adic probability Haar measure on $\Sym_m(\Z_p)\simeq\Z_p^{m(m+1)/2}$.  For each $S\subseteq\Sym_m(\Z_p)$ we thus have
\[\rho_p(S)=\lim_{r\to\infty}\frac{|\{a\modulo p^r:a\in S\}|}{|\Sym_m(\Z/p^r\Z)|},\]
provided that the limit exists. Note that if $F\in\Sym_m(\R_p)$ and $G$ is its $\GL_m(\Z_p)$-orbit, then we have
\begin{align*}
\rho_p(G)&=\lim_{\substack{q = p^r \\ q \to \infty}}\frac{|\GL_m(\Z/q\Z)|}{|\Ort_F(\Z/q\Z)|}\frac{1}{|\Sym_m(\Z/q\Z)|}=\frac{1}{\prod_{n=1}^m\zeta_p(n)}\lim_{\substack{q = p^r \\ q \to \infty}}\frac{q^{m(m-1)/2}}{|\Ort_F(\Z/q\Z)|}.
\end{align*}

\begin{definition}
Let $p$ be prime. Let $F\in\Sym_m(\Z_p)$ and let $G$ be its $\GL_m(\Z_p)$-orbit.
\begin{enumerate}[label=\upshape\arabic*.\ ]
    \item Following \cite{siegel}, the \emph{$p$-density} of $F$ is \[\delta_p(F)=\lim_{q=p^r\to\infty}\frac{q^{m(m-1)/2}}{|\Ort_F(\Z/q\Z)|}=\prod_{n=1}^m\zeta_p(n)\cdot\rho_p(G).\]
    \item Following \cite{CSMass}, the \emph{$p$-mass} of $F$ is
    $m_p(F)=D_p(F)^{(m+1)/2}\delta_p(F)$.
\end{enumerate}
\end{definition}

Let us state the Siegel mass formula using the language of regulators introduced above.

\begin{theorem}[\emph{Siegel mass formula} \cite{siegel,siegelIAS}]\label{massform}
If $G$ is a genus of quadratic forms of signature $(a,b)$, then
\begin{align*}
h_{\Reg}(G)=\frac{4}{\Vol(\Ort_{a,b}(\R))}\prod_{p}2m_p(G).
\end{align*}
\end{theorem}

\subsection{Standard mass and measure}\label{measure}

In these paragraphs, following Conway--Sloan \cite{CSMass} we introduce the notion of standard mass $\std(F)$ of an integral quadratic form (Conway--Sloan only worked with positive definite forms, and our definition is a natural generalization of theirs to forms of arbitrary signature). The regulator $\Reg(F)$ introduced in Section \ref{submass1} turns out to be a rational multiple of $\std(F)$, and their ratio will be called the \emph{measure} of $F$.

\begin{definition}[Standard mass]
Let $F\in\Sym_m(\Z)$ be nonsingular of signature $(a,b)$. The \emph{standard mass} of $F$ is given as follows.
    \begin{enumerate}
        \item If $m=2s-1$ is odd, then
        \[\std(F)=\frac{2}{\Vol(\Ort_{a,b}(\R))}\prod_{k=1}^{s-1}\zeta(2k).\]
        \item If $m=2s$ is even, then
        \[\std(F)=\frac{2}{\Vol(\Ort_{a,b}(\R))}\prod_{k=1}^{s-1}\zeta(2k)\cdot L(s,\chi_{(-1)^{s}D(F)}).\]
        Here, $\chi_{d}=\left(\frac{d}{\cdot}\right)$ is the Legendre symbol.
    \end{enumerate}
Additionally, if $G$ is the class, spinor genus, or genus of $F$, then we write $\std(G)=\std(F)$.
\end{definition}

\begin{remark}
    If $F$ is positive definite, then $\std(F)$ recovers the definition of standard mass given in \cite{CSMass}. We shall write $\std_m(d)$ for the standard mass of a positive definite quadratic form $F\in\Sym_m(\Z)$ of determinant $d$. We observe that $\std_m(d)$ is rational and independent of $d$ if $m$ is odd.
\end{remark}

\begin{definition}
Let $p$ be prime. Let $F\in\Sym_m(\Z_p)$. Following \cite{CSMass}, the \emph{standard $p$-mass} of $F$ is  given as follows.
    \begin{enumerate}
        \item If $m=2s-1$ is odd, then
        \[\std_p(F)=\frac{1}{2}\prod_{n=1}^{s-1}\zeta_p(2n).\]
        \item If $m=2s$ is even, then
        \[\std_p(F)=\frac{1}{2}L_p(s,\chi_{D(F)})\prod_{n=1}^{s-1}\zeta_p(2n) .\]
    \end{enumerate}
Here, $\zeta_p(a)=(1-p^{-a})^{-1}$ and $L_p(a,\chi_d)=(1-(\frac{d}{p})p^{-a})^{-1}$. If $G$ is the $\GL_m(\Z_p)$-equivalence class of $F$, we also write $\std_p(G)=\std_p(F)$.
\end{definition}

\begin{remark}
If $F\in\Sym_m(\Z)$ is nonsingular of signature $(a,b)$, then
    \[\std(F)=\frac{2}{\Vol(\Ort_{a,b}(\R))}\prod_p 2\std_p(F).\]
\end{remark}

\begin{definition}[Measure] We define the following.
\begin{enumerate}
    \item Let $F\in\Sym_m(\Z)$ be nonsingular. The \emph{measure} of $F$ is
    \[\mu(F)=\frac{\Reg(F)}{\std(F)}.\]
    \item If $G$ is a genus or spinor genus of quadratic forms with class representatives $F_1,\dots,F_{h(G)}$, we define the \emph{measure} of $G$ to be
    \[\mu^*(G)=\sum_{i=1}^{h(G)}\mu(F_i)=\frac{h_{\Reg}(G)}{\std(G)}.\]
    \item Let $p$ be prime. Let $F\in\Sym_m(\Z_p)$. The \emph{$p$-measure} of $F$ is
    \[\mu_p^*(F)=\frac{m_p(F)}{e_p\std_p(F)}\]
where $e_p=2$ if $p=2$, and $1$ otherwise. If $G$ is the $\GL_m(\Z_p)$-class of $F$, then we shall write $\mu_p^*(G)=\mu_p^*(F)$.
\end{enumerate}
\end{definition}

\begin{remark}
Given $F\in\Sym_m(\Z)$, we have $\mu_p^*(F)=1$ for all but finitely many primes $p$. The Siegel mass formula for a genus $G$ of quadratic forms containing classes $C_1,\dots,C_{h(G)}$ reads as
\[h_{\Reg}(G)=\std(G)\prod_p\mu_p^*(G),\quad\text{or}\quad\mu^*(G)=\sum_{i=1}^{h(G)}\mu(C)=\prod_p\mu_p^*(G).\]
It follows that $\mu^*(G)$ is rational.
\end{remark}

\begin{proposition}\label{equal-mass}
Let $G$ be a genus of nonsingular integral quadratic forms in $3$ or more variables. Then $\mu^*(S_1)=\mu^*(S_2)$ for any spinor genera $S_1$ and $S_2$ contained in $G$. In particular, for any spinor genus $S\in\SGen(G)$,
\[\mu^*(G) = \sg(G) \mu^*(S).\]
In particular, if $G$ is indefinite, then $\mu^*(G)=h(G)\mu(C)$ for any class $C\subseteq G$. Thus, $\mu(C)$ is rational for every class $C$ of nonsingular integral quadratic forms in $3$ or more variables.
\end{proposition}
\begin{proof}
This is folklore knowledge following work of Kneser \cite{Kneser1961}. We refer the reader to a discussion in \cite{SchulzePillot2014}.
\end{proof}

\begin{proposition}\label{mass-intp}\ 
Let $F\in\Sym_m(\Z)$ be a nonsingular quadratic form.
    \begin{enumerate}[label=\upshape{\arabic*}.]
        \item If $F$ is positive definite, then
        \[\mu(F)=\frac{1}{\std_m(D(F))}\frac{2}{|\Ort_F(\Z)|}.\]
        \item If $m=2n+1$ is odd and $F$ has signature $(1,m-1)$, then
        \[\mu(F)=\frac{1}{\std_m(D(F))}|\chiorb(\Sigma_F)|.\]
    \end{enumerate}
\end{proposition}

\begin{proof}
The first statement follows by definition. The second follows by noting that
\begin{align*}
\std(F) &= \frac{\Vol(\Ort_m(\R))}{\Vol(\Ort_{m-1}(\R))}\std_m(D(F)) = \frac{\pi^{m/2}}{\Gamma(m/2)}\std_m(D(F)), \\
        &= \frac{2^{2n}\pi^n n!}{(2n)!}\std_m(D(F)),
\end{align*}
and combining it with the Chern--Gauss--Bonnet formula.
\end{proof}

\subsection{Invariants of ternary forms}

From this point, we shall focus on ternary quadratic forms unless otherwise stated.

\begin{definition}\label{defcat}
We introduce the following notations.
\begin{enumerate}[label=\upshape\arabic*.\ ]
    \item $\mathcal F$ is the set of nonsingular primitive integral quadratic forms $F\in\Sym_3(\Z)$ of positive determinant $D(F)>0$ and $\Cl=\mathcal F/\GL_3(\Z)$ is the set of equivalence classes. 
    \item For each prime $p<\infty$,  $\mathcal F_{p}$ is the set of nonsingular primitive $\Z_p$-quadratic forms $F\in\Sym_3(\Z_p)$, and the set of equivalence classes is
    \[\Gen_{p}=\mathcal F_{p}/\GL_3(\Z_p).\]
    Let $\mathcal F_{\infty}$ denote the set of nonsingular real ternary quadratic forms of positive determinant and write $\Gen_{\infty}=\mathcal F_{\infty}/\GL_3(\R)$. For $p\leq\infty$, the inclusion $\mathcal F\to\mathcal F_{p}$ induces a map $\Cl\to\Gen_{p}$ denoted by $C\mapsto C_p$.
    \item Two forms $F,F'\in\mathcal F$ are in the same \emph{genus} if they have the same image in $\Gen_{p}$ for every $p\leq\infty$. Let $\Gen$ denote the set of genera of forms in $\mathcal F$.  
    \item  Two forms $F,F'\in\mathcal F$ are in the same \emph{spinor genus} if there exist $A\in\GL_3(\Q)$ and $B_p\in\GL_3(\Z_p)$ for every $p<\infty$ such that:
    \begin{enumerate}[label=\upshape\alph*.\ ]
        \item $F[A]=F'$,
        \item $F[B_p]=F'$ for every prime $p<\infty$,
        \item $\det(AB_p^{-1})=1$ for every $p<\infty$, and
        \item $AB_p^{-1}$ has trivial spinor norm (i.e.,~$\theta_p(AB_p^{-1})=1$) for every $p<\infty$. (See Section \ref{spinor} for the definition of spinor norm $\Ort_F(\Q_p)\to\Sq^*(\Q_p)$.)
    \end{enumerate} Then, $\SGen$ denotes the set of spinor genera of forms in $\mathcal F$. Each genus consists of finitely many spinor genera; see Section \ref{spinor} and references within for details. Each spinor genus contains finitely many classes.
    
\end{enumerate}
\end{definition}

\begin{definition}[Determinant] \ 
\begin{enumerate}[label=\upshape\arabic*.\ ]
    \item The determinant $D:\mathcal F\to\N$ induces functions on $\Cl, \SGen$, and $\Gen$, also denoted by the same letter. For each $d\in\N$, we   write
    \[\mathcal F(d)=\{F\in\mathcal F:D(F)=d\},\]
    and  define the sets $\Cl(d)$, $\SGen(d)$, and $\Gen(d)$ similarly.
    \item  For each prime $p<\infty$, the determinant $D:\mathcal F_p\to\Z_p$ induces a function
    \[D:\Gen_p\to\Z_p/(\Z_p^*)^2,\]
    whose fibers we shall denote $\Gen_p(d)=D^{-1}(d)$. We define the function
    \[D_p:\Gen_p\to\N,\]
    by $D_p(G_p)=p^{v_p(D(G_p))}=|D(G_p)|_p^{-1}$.
\end{enumerate}
\end{definition}

\begin{remark}
    If $F$ is a form, class, spinor genus, or genus, then $D(F)=\prod_{p}D_p(F_p)$.
\end{remark}

\begin{definition}[Hasse invariant]
Let $p\leq\infty$. Following \cite{jones-bk}, given $F\in\mathcal F_p$ we denote the $p$-adic Hasse invariant of $F$ by
\[c_p(F)=(-1,-D)_p(D_1,-D_2)_p(D_2,-D)_p,\]
where $D_1=F_{11}$, $D_2=F_{11}F_{22}-F_{12}^2$, and $D=D(F)$. Here $(-,-)_p$ is the Hilbert symbol. This differs from Cassels' definition \cite{Cassels} by the factor $(-1,-D)_p$. The Hasse invariant descends to a function $c_p:\Gen_p\to\{\pm1\}$.
\end{definition}

\begin{definition}
    For $\varpi\in\{\pm1\}$, we shall write $\mathcal F^{\varpi}=\{F\in\mathcal F:c_\infty(F)=\varpi\}$, and analogously define the sets $\Cl^\varpi$, $\SGen^\varpi$, and $\Gen^\varpi$.
\end{definition}

\begin{remark}
    If $G\in\Gen$, then we have $c_\infty(F)=\prod_p c_p(F)$. The product is finite in the sense that $c_p(G)=1$ for all but finitely many primes $p$.
\end{remark}

\begin{definition}
If $G\in\Gen$ is a genus, we write
\[\Cl(G)=\{C\in\Cl:C\subseteq G\}\quad \text{and}\quad\SGen(G)=\{S\in\SGen:S\subseteq G\}.\]
We shall write \[h(G)=\#\Cl(G)\quad  \text{and}\quad  \sg(G)=\#\SGen(G).\]
By Eichler's theorem, $\Cl^1=\SGen^1$. Thus, if $G\in\Gen^1$ is the genus of an indefinite form then $\Cl(G)=\SGen(G)$ and therefore $h(G)=\sg(G)$.
\end{definition}

\subsubsection{Indefinite classes and periodic orbits}\label{subperiodic}\ 

In these paragraphs, we describe the bijective correspondence between the set $\Cl^1$ of primitive indefinite integral ternary quadratic form classes and the set of periodic $\SO_{1,2}(\R)$-orbits in the lattice space $\SL_3(\R)/\SL_3(\Z)$. Let us write
\[\mathbf I_{1,2}=\begin{bmatrix}1 & 0 & 0\\0 & -1 & 0 \\0 & 0 & -1\end{bmatrix}\in\Sym(\Z),\]
so that $\SO_{1,2}=\SO_{\mathbf I_{1,2}}$. Let $F\in\mathcal F^1$ be an indefinite primitive integral ternary quadratic form of positive determinant $D(F)>0$. Note that the signature of $F$ is $(1,2)$. The real quadratic form
\[\tilde F=\frac{F}{D(F)^{1/3}},\]
has determinant $1$. Choose $g_F\in\SL_3(\R)$ such that $\mathbf I_{1,2}[g_F]=\tilde F$. Then
\[\{g\in \SL_3(\R):\mathbf I_{1,2}[g]=\tilde F\}=\SO_{1,2}(\R) g_F,\]
is a torsor over $\SO_{1,2}(\R)$. Consider the natural map
\[\Phi_{g_{F}}:\SO_{1,2}(\R)\to \SL_3(\R)/\SL_3(\Z),\]
given by $\Phi_{g_F}(A)=Ag_{F}\SL_3(\Z)$. Its image is independent of the choice of $g_F$ and depends only on the $\GL_3(\Z)$-equivalence class $C$ of $F$ (note that $-1\in\GL_3(\Z)$ acts trivially on $\mathcal F$). We shall therefore denote by image by
\[X_C=\Phi_{g_F}(\SO_{1,2}(\R))\subseteq\SL_3(\R)/\SL_3(\Z).\]
We note that $X_C$ is a periodic (i.e.~closed and of finite volume) $\SO_{1,2}(\R)$-orbit for the natural left action of $\SO_{1,2}(\R)$ on $\SL_3(\R)/\SL_3(\Z)$ by left multiplication. The map $\SO_F(\R)\to\SL_3(\R)$ given by $h\mapsto g_F h$ induces an embedding
\begin{align*}
    \Ort_F(\R)/\Ort_F(\Z)=\SO_F(\R)/\SO_F(\Z)&\simeq X_C\subseteq \SL_3(\R)/\SL_3(\Z).
\end{align*}
Indeed, this map is well-defined. Writing $K_F=\SO(n)\cap \SO_F(\R)$ for the maximal compact subgroup of $\SO_F(\R)$, we have a totally geodesic immersion
\[\Sigma_F=K_F\backslash\SO_F(\R)/\SO_F(\Z)\to\SO(3)\backslash\SL_3(\R)/\SL_3(\Z)\]
of the hyperbolic orbifold $\Sigma_F$, of which $X_C$ is the unit tangent bundle. The orbifold Euler characteristic $\chiorb(C):=\chiorb(S_F)$  is a rational invariant of the class $C$. See Prop.~\ref{mass-intp} for the relation between this invariant and the normalized Siegel mass, introduced in Section \ref{submass1}.
\begin{remark}
    One can show that the above assignment $C\to X_C$ sets up a bijection
    \[\Cl^1\leftrightarrow\{\textup{periodic $\SO_{1,2}(\R)$-orbits in $\SL_3(\R)/\SL_3(\Z)$.}\}\]
    One may compare the above with the bijection between the equivalence classes of primitive integral indefinite binary quadratic forms and the periodic $\SO_{1,1}(\R)$-orbits in $\SL_2(\R)/\SL_2(\Z)$. The latter are precisely the periodic orbits of the geodesic flow on the unit tangent bundle of the modular surface.
\end{remark}

\subsubsection{Computation of local measures}\ 

We will need various instances of Prop.~11.2 of \cite{GGSW} in our analysis and so restate it here.

\begin{proposition}\label{prelim-6}
Let $p\geq 3$ be a prime and $d\in\Z\backslash\{0\} $ and write $d=p^k d_0$ where $p\nmid d_0$  $($so $k=v_p(d)\geq 0$$)$.  Then, the classes in $\text{Gen}_p(p^k)$ satisfy the following:
\begin{enumerate}[label=\upshape{\arabic*}.]
\item If $k=0$, there is one class  with $C_p=1$, and $\mu_p^* =1$.
\item If $k=1$, there are two classes  with $C_p=\pm1$, and
\[
    \mu_p^*=\frac12\left(p+C_p\right).
\]
 \item If $k=2$, there are a total of four classes. The first two with $C_p=\pm1$, and
\[
    \mu_p^*=\frac12\left(p+C_p\right);
\]
and the remaining two with both having $C_p=1$, and
\[
    \mu_p^*=\frac{p}{2}\left(p \pm 1\right).
\]
\item If $k\geq 3$, we have the following types:\medskip \\
\noindent\underline{Type I.}\quad
If $k$ is odd, there are two classes with $C_p=\pm1$, and
\[
    \mu_p^*=\frac12\left(1+C_pp^{-1}\right)p^k.
\]
If $k$ is even, there are two forms with $C_p=1$, and
\[
    \mu_p^*=\frac12\left(1\pm p^{-1}\right)p^k.
\]
\noindent\underline{Type II.}\quad
 If $k\equiv0\pmod4$, there are two forms with
$C_p=1$, and
\[
    \mu_p^*=\frac12\left(1\pm p^{-1}\right)p^{k/2}.
\]
If $k\equiv2\pmod4$, there are two forms with $C_p=\pm1$, and
\[
    \mu_p^*=\frac12\left(1+C_pp^{-1}\right)p^{k/2}.
\]
\noindent\underline{Type III.}\quad  If $ k=t_2+t_3$ with  $1\leq t_2<t_3$, then in all cases
\[
    \mu_p^*=\frac14p^{t_3}\left(1-p^{-2}\right).
\]
If $t_2,t_3$ are both even, there are four forms, all with
$C_p=1$. Otherwise, there are two forms for each of
$C_p=\pm1$.

\end{enumerate}
\end{proposition}

\begin{corollary}\label{prelim-7a}
    Let $\ell\geq5$ be a finite prime. Let $\beta,\gamma\in\R$. Write $\beta^*=\max(-\beta,0)$ and $\gamma^*=\max(-\gamma,0)$.
    \begin{enumerate}[label=\upshape{\arabic*}.]
        \item For every $G\in\Gen_2$, we have $|\mu_2(G)|_\ell^\beta=1$.
        \item For every $G\in\Gen_\ell$ with $v_\ell(D(G))\geq1$, we have
        \[|\mu_\ell(G)|_\ell^\beta\leq \left(\frac{D_\ell(G)}{\ell}\right)^{\beta^*}.\]
        \item For any prime $p\neq 2,\ell$ and $G\in\Gen_p$ with $v_p(D(G))\geq2$ we have
        \[|\mu_p^*(G)|_\ell^\beta\leq p^{\beta^*}.\]
        \item There is an absolute constant $C>0$ such that for every prime $p$ and $G\in\Gen_p$ with $v_p(G)\geq2$ we have
        \[\left|\frac{\mu_p^*(G)}{D_p(G)}\right|_\infty^\gamma\leq C\cdot D_p(G)^{\gamma^*/2}\]
    \end{enumerate}
    In particular, if $\beta^*+\gamma^*/2<1$, then for every prime $p<\infty$ we have
    \[\sum_{G\in\Gen_p}\frac{|\mu_p^*(G)|_\ell^\beta}{D_p(G)}\left|\frac{\mu_p^*(G)}{D_p(G)}\right|_\infty^\gamma<\infty.\]
\end{corollary}
\begin{proof}
Part (1) follows from the fact, observed from \cite{CSMass}, that $\mu_2(G)$ is a product of integral powers of $2$ and $3$ for each $G\in\Gen_2$. Parts (2) and (3) follow by examination of Prop.~\ref{prelim-6}, noting that $\ell$ divides at most one of $p-1$ and $p+1$ for $p\neq 2,\ell$. Part (4) similarly follows from Prop.~\ref{prelim-6}. Now, assume that $\beta^*+\gamma^*/2<1$. To conclude the last statement, we see by the analysis in \cite[Proof of Theorem 3.1]{GGSW} that $\#\Gen_p(p^k)\leq 8 k$, say, for $k\geq1$. It therefore follows that
\[\sum_{\substack{G\in\Gen_2\\v_2(D(G))\geq 1}}\frac{|\mu_2^*(G)|_\ell^\beta}{D_2(G)}\left|\frac{\mu_2^*(G)}{D_2(G)}\right|_\infty^\gamma\leq C\sum_{k=1}^\infty\frac{8k}{(2^{1-\gamma^*/2})^k}<\infty,\]
for some finite $A>0$. Next,
\begin{align*}
\sum_{\substack{G\in\Gen_\ell\\v_\ell(D(G))\geq1}}\frac{|\mu_\ell^*(G)|_\ell^\beta}{D_\ell(G)}&\leq C\sum_{k=1}^\infty\frac{8k}{\ell^k}\left(\frac{\ell^k}{\ell}\right)^{\beta^*}(\ell^k)^{\gamma^*/2},\\
&\leq 8\ell^{-\beta^*}C\sum_{k=1}^\infty\frac{k}{(\ell^{1-\beta^*-\gamma^*/2})^k}<\infty.
\end{align*}
Lastly, if $p\neq 2,\ell$, then we have
\[\sum_{\substack{G\in\Gen_p\\v_p(D(G))\geq 2}}\frac{|\mu_p^*(G)|_\ell^\beta}{D_p(G)}\left|\frac{\mu_p^*(G)}{D_p(G)}\right|_\infty^\gamma\leq 8p^{\beta^*}C\sum_{k=2}^\infty\frac{k }{(p^{1-\gamma^*/2})^k}<\infty.\]
This proves the   result.
\end{proof}

\begin{corollary}\label{prelim-7b}
Let $\ell\geq5$ be a finite prime. Let $-1<\beta<\infty$. Write $\beta^*=\max(-\beta,0)$. There is an absolute constant $C>0$ such that, for all primes $p$, we have
\[\sum_{G\in\Gen_p\setminus\Gen_p^\tau}\frac{|\mu_p^*(G)|_\ell^\beta}{\mu_p^*(G)}\leq C\cdot p^{\delta}\frac{2p-1}{p(p-1)^2}.\]
\end{corollary}
\begin{proof}
Consider the first the case where $p\neq 2,\ell$. Note that for $p\geq3$,
\[p\geq\frac{p+1}{2}\geq \frac{p-1}{2}\geq \frac{p}{3}.\]
 Then examination of Prop.~\ref{prelim-6} shows that
\[\sum_{G\in\Gen_p\setminus\Gen_p^\tau}\frac{|\mu_p^*(G)|_\ell^\beta}{\mu_p^*(G)}\ll p^{\beta^*}\sum_{k=2}^\infty\frac{k}{p^k}=p^{\beta^*}\frac{2p-1}{p(p-1)^2},\]
where the implied constant is independent of $p$. To complete the proof, it suffices to show that the sum is finite for $p\in\{2,\ell\}$. For $p=\ell$ we have the estimate, again using Prop.~\ref{prelim-6},
\[\sum_{G\in\Gen_\ell\setminus\Gen_\ell^\tau}\frac{|\mu_\ell^*(G)|_\ell^\beta}{\mu_\ell^*(G)}\ll\sum_{k=2}^\infty\frac{k^2}{\ell^k}\left(\frac{\ell^k}{\ell}\right)^{\delta^*}<\infty,\]
since we have $1-\delta^*>0$.
The argument for $p=2$ can be deduced from \cite{CSMass}.\end{proof}

\subsection{Mass functions}\label{submass2}\

\begin{definition}
    A \emph{mass function} is a function $f:\Gen\to\C$ together with a function $f_p:\Gen_p\to\C$ for each prime $p$ such that:
    \begin{enumerate}[label=\upshape{\roman*}.]
        \item when $D(G)=1$, for all but finitely many $p$ we have $f_p(G)=1$ ,
        \item for every $p$, we have $f_p(uG)=f_p(G)$ for every $u\in\Z_p^*$ and $G\in\Gen_p$, and
        \item $f(G)=\prod_pf_p(G)$ for all $G\in\Gen$.
    \end{enumerate}
\end{definition}

Mass functions are closed under multiplication, and under composition by totally multiplicative functions.

\begin{definition}\label{sel-mult}A mass function $f:\Gen\to\C$ is said to be {\it Selberg multiplicative} if for each $d\geq 1$, the function 
\[F(d):= \sum_{G\in \text{Gen}(d)} f(G),\] is a Selberg multiplicative function (see Appendix \ref{prelim}). Then, we write $F(d)=\prod_pF_p(\nu_p(n))$, where for a prime $p$, and $r\geq0$
\[F_p(r):=\sum_{G\in\Gen_p(p^r)}f_p(G).\]\end{definition}

\begin{example}\label{examp-1}
    Some examples of mass functions are the following:
\begin{enumerate}[label=\upshape{\arabic*}.]
    \item The determinant map $D:\Gen\to\N$ is a mass function, with $D_p=|D(-)|_p^{-1}$.
    \item The archimedean Hasse invariant $c_\infty:\Gen\to\{\pm1\}$ is a mass function. We have $c_\infty=\prod_p c_p$ where $c_p$ denotes the $p$-adic Hasse invariant.
    \item By the Siegel mass formula, the measure $\mu^*$ is a mass function, and $\mu^*=\prod_p\mu_p^*$.
\end{enumerate}
They are all  Selberg multiplicative.
\end{example}

\begin{proposition}\label{ds-mass}
Let $M$ and $f$ be mass functions. Then formally
    \[\sum_{G\in\Gen}\frac{f(G)}{M(G)^s}=\prod_p\left(\sum_{k=0}^\infty\sum_{G\in\Gen_p(p^k)}\frac{f_p(G)}{M_p(G)^s}\right).\]
\end{proposition}

\begin{proof}
The natural map
\[\Gen\to \sideset{}{'}\prod_p\left(\coprod_{k\geq0}\Gen_p(p^k),\Gen_p(1)\right)\]
to the restricted product given by $G\mapsto \left(\frac{D_p(G)}{D(G)}G_p\right),$ is a bijection.
\end{proof}

\begin{proposition}\
\begin{enumerate}[label=\upshape{\arabic*}.]
    \item We have
    \[\sum_{G\in\Gen}\frac{1}{D(G)^s}=A(2^{-s})\cdot\frac{\zeta(s)^2}{\zeta(2s)}\frac{\zeta(2s)^2}{\zeta(4s)}\]
    for some explicit rational function $A(t)\in\Q(t)$.
    \item We have
    \[\sum_{G\in\Gen}\frac{\mu^*(G)}{D(G)^s}=\frac{\zeta(s-1)\zeta(2s-1)}{\zeta(3s)}.\]
\end{enumerate}
\end{proposition}
    
\begin{remark}
The above result shows that our normalization of Siegel mass enjoys a pleasant multiplicative property when genera of definite forms and indefinite forms are considered together. For each $d\in\N$, we write $H(d)=\sum_{D(G)=d}\mu^*(G)$. Then $H(d)$ is a multiplicative function of $d$.
\end{remark}

\subsection{Spinor genera, genera and packets}\label{spinor-pkts}\

\subsubsection{Spinor genera}\label{spinor}\ 

Our work in \cite{GGSW} relies on the structure theory of spinor genera (See \cite{Conway1993} and references within), reviewed below. The presentation here largely follows that given in \cite{GGSW}, with some modifications necessitated by consideration of spinor genera for definite forms and a finer analysis of ramification (needed later for the analytic study of spinor genera ordered by $\mu^*$).

Given a commutative ring $\mathcal{A}$, we shall write $\mathcal{A}^*$ for the group of units in $\mathcal{A}$, and write $\Sq^*(\mathcal{A})=\mathcal{A}^*/(\mathcal{A}^*)^2$ for the group of unit square classes in $\mathcal{A}$. We note that $\Sq^*(\mathcal{A}_1\times \mathcal{A}_2)=\Sq^*(\mathcal{A}_1)\times\Sq^*(\mathcal{A}_2)$ for all commutative rings $\mathcal{A}_1$ and $\mathcal{A}_2$. For our application $\mathcal{A}$ is  $\mathbb{Z}$, $\mathbb{Z}_p$, $\Q_p$, or $\Z/n\Z$. We have
\begin{align*}
\Sq^*(\Z_2)&\simeq\Sq^*(\Z/8\Z)=(\Z/8\Z)^*,\quad\text{and}\\
\Sq^*(\Z_p)&\simeq\Sq^*(\Z/p\Z)\simeq\{\pm1\}\quad\text{for $p\geq3$ prime.}
\end{align*}
If $p\geq2$ prime, we have a group isomorphism 
$\Sq^*(\Q_p)\simeq\Z/2\Z\times\Sq^*(\Z_p)$
given by $p^\alpha a\mapsto ([\alpha],a)$ for all $\alpha\in\Z$ and $a\in\Z_p^*$.

In these paragraphs, following \cite{Conway1993} we treat $-1$ as a ``prime'' so that each nonzero $d\in\Z$ has a unique prime factorization $d=(-1)^a\prod_{p}p^{v_p(d)}$ with $a=v_{-1}(d)\in\{0,1\}$. Thus $-1$ is a prime divisor of a nonzero integer $d$ if and only if $d< 0$. We write $\Q_{-1}=\R$ so that $\Sq^*(\Q_{-1})\simeq\{\pm1\}$. We formally define the symbol $\Sq^*(\Z_{-1})$ to mean the trivial group \[\Sq^*(\Z_{-1})=\{1\}\leq\Sq^*(\Q_{-1}).\]Below, when we say $p$ is a prime, we shall include $p=-1$ unless otherwise specified.

\begin{definition}[\emph{Spinor operators}]
Let $n\geq1$. Let $p$ be prime. Let $F\in\Sym_m(\Z)$ be nonsingular.
\begin{enumerate}[label=\upshape\arabic*.\ ]
	\item The \emph{spinor norm} is the homomorphism $sn:\Ort_F(\Q_p)\to\Sq^*(\Q_p)$ given by
	\[sn(B)=F(v_1)\cdots F(v_k)\in\Sq^*(\Q_p),\] where $v_1,\dots,v_k$ is any sequence of vectors in $\Q_p^n$ such that $B$ is a product of reflections $B=r_{v_1}\circ\cdots\circ r_{v_k}$ in said vectors.
	\item The group $\mathcal A_p(F)\leq\Sq^*(\Q_p)$ of \emph{$p$-adically automorphous numbers} for $F$ is
	\[\mathcal A_p(F)= sn(\SO_F(\Z_p))\leq\Sq^*(\Q_p).\]
    Note that $\mathcal A_p(F)$ only depends on the $\GL_m(\Z_p)$-equivalence class of $F$.
	\item Let $\Pi$ be a finite set of primes. We define the group homomorphism
\[\Delta_p^{(\Pi)}:\Sq^*(\Q_p)\to\prod_{q\in\Pi}\Sq^*(\Z_q),\quad\Delta_p^{(\Pi)}(A)_q=\left\{\begin{array}{l l} A/p^{v_p(A)} & q=p,\\p^{v_p(A)} & q\neq p.\end{array}\right.\]
If $A\in\mathcal A_p(F)$, we call $\Delta_p^{(\Pi)}(A)$ the \emph{spinor operator} associated to $(p,F,A,\Pi)$.
\end{enumerate}
\end{definition}

The content of the following lemma is found in \cite[Ch.~9]{Conway1993} Sections 9.4--9.6. See also \cite[Ch.3]{Montesinos}.

\begin{lemma}\label{p-autonum}
Let $n\geq3$ and $d\geq1$ be integers. Let $F\in\Sym_n(\Z)$ be a nonsingular integral $n$-ary quadratic form of determinant $d$. Let $p \! \mid \!  -2d$ be prime satisfying $p\neq2$. Assume that $F$ is $\GL_n(\Z_p)$-equivalent to a diagonal form
\[a_1p^{t_1}x_1^2+\dots+a_np^{t_n}x_n^2,\] with $a_i\in\Z_p^*$ for all $i=1,\dots,n$. Then the group of $p$-adically automorphous numbers for $F$ is given by $\mathcal A_p(F)=\langle S_1\cup S_2\rangle$ where:
\begin{enumerate}[label=\upshape\arabic*.\ ]
	\item $S_1=\{a_ip^{t_i}/a_jp^{t_j}:1\leq i,j\leq n\}$, and
	\item $S_2=\Sq^*(\Z_p)$ if $t_i=t_j$ for some $i\neq j$, and $S_2=\emptyset$ otherwise.
\end{enumerate}
\end{lemma}

\begin{definition}
    Given $d\in\N$ and $k\in\N$, write $W_k(d)=\{\text{$p\geq2$ prime}:p^k\mid d\}$. We shall define $\omega_k(d)=|W_k(d)|$. As is customary, we also write $\omega(d)=\omega_1(d)$.
\end{definition}
\begin{definition}\label{gen-tau}
    Let $p\geq3$ be prime. We write $\Gen_p^\tau$ for the set of $\GL_3(\Z_p)$-classes of the primitive diagonal forms
    \[ax_1^2+bx_2^2+cx_3^2\quad\text{or}\quad ax_1^2+bx_2^2+cpx_3^2\quad\text{or}\quad ax_1^2+bpx_2^2+cpx_3^2,\]
    for some $a,b,c\in\Z_p^*$.
\end{definition}
\begin{definition}
    Let $G\in\Gen$ be a genus. We write $T(G)$ for the set of primes $p\geq3$ such that $p\mid D(G)$ and $G_p\in\Gen_p^\tau$. Note that $T(G)\subseteq W_1(D(G))\setminus W_3(D(G))$.
\end{definition}

\begin{proposition}\label{propE1}
Let $d\in\N$. Let $P=8\prod_{p\in W_3(d)}p$, and let $\Pi=\{2\}\cup W_3(d)$ be the set of prime divisors of $P$. Then 
$\prod_{q\in\Pi}\Sq^*(\Z_q)=\Sq^*(\Z/P\Z)$, and for each $G\in\Gen(d)$ there is a bijection
\[\SGen(G)\simeq\frac{\Sq^*(\Z/P\Z)}{\langle\Delta_p^{(\Pi)}(\mathcal A_p(G_p)):\textup{$p\in \{-1,2\}\cup W_1(d)\setminus T(G)$}\rangle\cdot \langle p:p\in T(G)\rangle}.\]
In particular, we have 
\begin{enumerate}[label=\upshape\arabic*.\ ]
	\item $\sg(G)$ is a power of $2$ satisfying $\sg(G)\leq 2^{2+\omega_3(d)}$, and
	\item If $\sg(G)\geq2$ then there is a proper subgroup $\mathcal C\leq \Sq^*(\Z/P\Z)$ of index $2$ such that $p\in\mathcal C$ for every prime $p\in T(G)$, which includes those primes $p\geq3$ with $v_p(d)=1$.
\end{enumerate}
\end{proposition}

\begin{proof}
Let $\Pi'=\{-1,2\}\cup W_1(d)$. It follows by \cite[Ch.~9, Thm.~17 and Sec.~9.4]{Conway1993} (See also \cite[Sec.~3.5]{Montesinos}) that the set of spinor genera in $G\in\Gen(d)$ is given by
\[\SGen(G)\simeq\frac{\prod_{q\in\Pi'}\Sq^*(\Z_q)}{\langle\Delta_p^{(\Pi')}(\mathcal A_p(F)):p\in\Pi'\rangle}.\]
Now, suppose that $p=-1$ or $p\in W_1(d)\setminus W_3(d)$. Then $G_p$ is the $\GL_3(\Z_p)$-class of a primitive diagonal form
\[ax_1^2+bp^{t_2}x_2^2+cp^{t_3}x_3^2\]
with $a,b,c\in\Z_p^*$ such that:
\begin{enumerate}[label=\upshape\arabic*.\ ]
    \item $t_1=t_2=0$ and $t_3\in\{1,2\}$, or
    \item $t_1=0$ and $t_2=t_3=1$.
\end{enumerate}
By Lemma~\ref{p-autonum}, we have $\Delta_p^{(\Pi')}(\mathcal A_p(F))\supseteq\Delta_p(\Sq^*(\Z_p))=\Sq^*(\Z_p)$; in the terminology of \cite{Conway1993}, such primes are called \emph{tractable}. Moreover, if $p$ is an odd prime such that $v_p(d)\leq 2$, then we have $\mathcal A_p(F)=\langle p,\Sq^*(\Z_p)\rangle$ again by Lemma~\ref{p-autonum}. We thus have
\begin{align*}
\Cl(\Ge(F))&\simeq\frac{\prod_{p\in\Pi'}\Sq^*(\Z_p)}{\langle\Delta_p^{(\Pi')}(\mathcal A_p(F)):p\in\Pi'\rangle}\ ,\\
&=\frac{\prod_{p\in\Pi'}\Sq^*(\Z_p)}{\langle\Delta_p^{(\Pi')}(\mathcal A_p(F)):p\in\Pi'\setminus T(G)\rangle\cdot \langle\Delta_p^{(\Pi')}(p),\Sq^*(\Z_p):p\in T(G)\rangle}\ ,\\
&\simeq\frac{\prod_{q\in\Pi}\Sq^*(\Z_q)}{\langle\Delta_p^{(\Pi)}(\mathcal A_p(F)):p\in\Pi'\setminus T(G)\rangle\cdot \langle p:p\in T(G)\rangle}\ ,
\end{align*}
from which the  result follows.
\end{proof}

\subsubsection{Packets}\label{packets2}\

\begin{definition}
    Let $A$ be a finite set of primes containing $2$. A \emph{packet} $H$ with \emph{support $\Supp(H)=A$} is a sequence of local equivalence classes
    \[H=(H_p)_{p\in S}\in\prod_{p\in A}\Gen_p.\]
    Let $\Packet$ denote the set of all packets.
\end{definition}

\begin{definition}
    Given a mass function $f:\Gen\to\C$ and $H\in\Packet$, we denote
    \[f(H)=\prod_{p\in\Supp(H)}f_p(H_p).\]
In particular, we define the following.
\begin{enumerate}[label=\upshape\arabic*.\ ]
    \item We define $D(H)=\prod_{p\in\Supp(H)} D_p(H_p)$. Note that $D(H)$ need not have the same class as $D(H_p)\modulo(\Z_p^*)^2$.
    \item We define $\mu^*(H)=\prod_{p\in\Supp(H)}\mu_p^*(H_p)$.
    \item We define $c_\infty(H)=\prod_{p\in\Supp(H)}c_p(H_p)$.
\end{enumerate}
\end{definition}

\begin{definition}[\mbox{\mathversion{normal}$A$}-root packets]
    Let $A$ be a finite set of primes containing $2$. An $A$-root packet is a packet $H\in\Packet$ such that:
    \begin{enumerate}[label=\upshape\arabic*.\ ]
        \item $A\subseteq\Supp(H)$, and
        \item for any prime $p\in\Supp(H)\setminus A$, we have $v_p(D_p(H_p))\geq2$.
    \end{enumerate}
    We write $\mathfrak R_A$ for the set of all $A$-root packets. If $A=\{p\text{ prime}:p\leq B\}$ for some real number $B\geq2$, we shall also write $\mathfrak R_A=\mathfrak R_B$ and say $B$-root packets to mean $A$-root packets.
\end{definition}

\begin{remark}
    The above definition generalizes that of a root packet introduced in \cite{GGSW}, which coincides with a $2$-root packet.
\end{remark}

\begin{definition}
    Let $H\in\Packet$ be a packet. We write $\Gen_H$ for the set of genera $G\in\Gen$ such that:
    \begin{enumerate}[label=\upshape\arabic*.\ ]
        \item $G_p=H_p$ for every $p\in\Supp(H)$, and
        \item $v_p(D_p(G_p))\leq1$ for every $p\notin\Supp(H)$.
    \end{enumerate}
    We denote by $\SGen_H$ the set of those spinor genera whose genera lie in $\Gen_H$. Note that for any finite set $A$ of primes we have
    \[\Gen=\coprod_{H\in\mathfrak R_{A}}\Gen_H,\quad \SGen=\coprod_{H\in\mathfrak R_{A}}\SGen_H.\]
    We shall write $\Gen_H^\varpi=\Gen^\varpi\cap \Gen_H$ and $\SGen_H^\varpi=\SGen^\varpi\cap \SGen_H$.
\end{definition}

\section{Preliminary asymptotics}\label{sect:preliminary}
Throughout this section, we fix the following data:
\begin{enumerate}[label=\upshape\arabic*.\ ]
    \item a prime $\ell\geq5$;
    \item $\beta\in(-1,\infty]$ and $\beta^*=\max(-\beta,0) <1$;
    \item $\gamma\in\R$; and
    \item an odd squarefree positive integer $P$.
\end{enumerate}

\begin{notation}\label{notation} We fix the following notation.
     \begin{enumerate}[label=\upshape{\roman*}.]
        \item Let $q$ and $a$ be integers such that 
	$1 \le a \le q$ and $(q,a) = 1$. We write for $t\geq 2$
\[ \pi(t;q,a) = \sum_{\substack{p\leq t\\p\equiv a\tmod{q}}}1 = \frac{\text{Li}(t)}{\phi(q)} + E(t;q,a), \] and \[ \psi(t;q,a) = \sum_{\substack{p\leq t\\p\equiv a\tmod{q}}}\log{p} =  \frac{t}{\phi(q)} + E^*(t;q,a). \]
        \item For a Dirichlet character $\chi$ modulo $q$, $\delta_{\chi}=1$ if $\chi$ is principal, and is zero otherwise.
        \item $\delta_j = \frac{1}{2}(1+(-1)^j)$.
    \end{enumerate}
\end{notation}

\begin{definition}\label{def3.1}\
\begin{enumerate}[label=\upshape\arabic*.\ ]
    \item We define
\[\alpha=\alpha(\ell,\beta)=1-\frac{1}{\ell-1}+\frac{1}{\ell^{1+\beta}-1}.\]
    \item Given an Dirichlet character $\chi$ modulo $\ell$, given a quadratic Dirichlet character $\rho$ modulo $8P$, given $j\in\{0,1\}$, and $\ve\in\{\pm1\}$, define
    \[f_{\chi,\rho,j,\ve}(p)=\ve^j\rho(p)\mathbf\chi\left(\left(\frac{p+\ve}{2}\right)_{\ell'}\right)\left|\frac{p+\ve}{2}\right|_\ell^\beta\left|\frac{p+\ve}{2p}\right|_\infty^\gamma\]
    for each prime $p$. Define $f_{\chi,\rho,j}(p)=\sum_{\ve=\pm1}f_{\chi,\rho,j,\ve}(p)$ and define the Eulerian Dirichlet series
    \begin{align}\label{building}
    A_{\chi,\rho,j}(s)&=\sum_{d=1}^\infty\frac{f_{\chi,\rho,j}(d)}{d^s}=\prod_{p}\left(1+\frac{1}{p^s}\sum_{\ve=\pm1}f_{\chi,\rho,j,\ve}(p)\right).
    \end{align}
\end{enumerate}    
\end{definition}

\begin{proposition}\label{pre1}
For each $(\chi,\rho,j)$, there is a constant $\alpha_{\chi,\rho,j}\in\C$ such that for all $X\geq2$ and all $A>0$ sufficiently large we have
\[\sum_{p\leq X}f_{\chi,\rho,j}(p)\log p=\alpha_{\chi,\rho,j} X+O\left(\frac{X}{(\log X)^A}\right).\]
If $\chi=\chiz$ and $\rho=\rho^{}_0$ are principal characters and $j=0$, we have $\alpha_{\chiz,\rho^{}_0,0}=2^{1-\gamma}\alpha(\ell,\beta)$. If $(\chi,\rho,j)\neq(\chiz,\rho^{}_0,0)$, then $|\alpha_{\chi,\rho,j}|<2^{1-\gamma}\alpha(\ell,\beta)$.
\end{proposition}

\begin{proof} We shall write $f(d)=f_{\chi,\rho,j}(d)$ for simplicity. Setting $L=[\log X/\log \ell]$ we have, using $\sum_{p\leq X}(\log p)/p\ll\log X$,
\begin{gather}
 \label{3a}\\
    \sum_{p\leq X}f(p)\log{p}=\frac{\overline{\chi(2)}}{2^\gamma}\sum_{k=0}^L\ell^{-k\beta}\sum_{\ve=\pm1}\ve^j\sum_{\substack{p\leq X\\\ell^k\mid p+\ve}}\chi\left(\frac{p+\ve}{\ell^k}\right)\rho(p)\log p +O( X^{\beta^*}\log X).\nonumber
\end{gather}
Consider first the case $\rho=\rho^{}_0$. The main term in \eqref{3a} is $\frac{1}{2^\gamma}$ times
    \[\overline{\chi(2)}
    \sum_{\ve=\pm1}\ve^j\left(\sum_{a\in(\Z/\ell\Z)^*}\chi(a)\sum_{k=1}^L\ell^{-k\beta} \psi(X;\ell^{k+1},a\ell^k-\ve) + \sum_{\substack{p\leq X}}\chi\left(p+\ve\right)\log p\right), \]
    which is $\frac{1}{2^\gamma}$ times
    
 \begin{align*}
     2\delta_j\delta_{\chi}(\ell-1)X&\left(\sum_{k=1}^L\ell^{-k\beta}\frac{1}{\phi(\ell^{k+1})}\right) + \overline{\chi(2)}\sum_{\ve\in\{\pm1\}}\ve^j\sum_{\substack{p\leq X}}\chi\left(p+\ve\right)\log p \nonumber\\ &+\overline{\chi(2)}\sum_{\ve\in\{\pm1\}}\ve^j\sum_{a\in(\Z/\ell\Z)^*}\chi(a)\left(\sum_{k=1}^L\ell^{-k\beta} E^*(X;\ell^{k+1},a\ell^k-\ve)\right), 
 \end{align*} 
 
 \begin{align}\label{3c}
 =\  2\delta_j\delta_{\chi}&\left(\ell^{(1+\beta)}-1\right)^{-1}X  + O(X^{-\beta})+ \overline{\chi(2)}\sum_{\ve\in\{\pm1\}}\ve^j\sum_{\substack{p\leq X}}\chi\left(p+\ve\right)\log p\nonumber \\ &+ \overline{\chi(2)}\sum_{\ve=\pm1}\ve^j\sum_{a\in(\Z/\ell\Z)^*}\chi(a)\left(\sum_{k=1}^L\ell^{-k\beta} E^*(X;\ell^{k+1},a\ell^k-\ve)\right). 
\end{align} 
Applying Siegel-Walfitz (as in Lemma \ref{prelim-1}) to the second sum in \eqref{3c}, and summing over the $\ve$, we see that we get 
\begin{equation}\label{3d}
\sum_{\ve=\pm1}\ve^j\sum_{\substack{p\leq X}}\chi\left(p+\ve\right)\log p =\begin{cases}
 2\delta_j\frac{\ell-2}{\ell-1}X + \mathcal{E}(X),&\quad \text{if}\ \chi=\chiz,\\ \\ -\left[1+(-1)^j\chi(-1)\right]\frac{X}{\ell-1} +\mathcal{E}(X),&\quad  \text{if}\ \chi\neq\chiz,
\end{cases} 
\end{equation}
     with 
     $\mathcal{E}(X) = O(\frac{X}{(\log{X})^{A}})$. To estimate the sum over $k\geq 1$ in \eqref{3c}, we split the sum onto two ranges: $1\leq k \leq K$ and $K<k\leq L$, with 
\[K=C\frac{\log\log X}{\log \ell}.\]
For $k\geq K$, we use the trivial bound that \[E^*(X;\ell^{k+1},a\ell^k-\ve) \ll \frac{X}{\phi(\ell^{k+1})}+ \frac{X\log{X}}{\ell^{k+1}},\]
so that the contribution to the sum is $O(X(\log{X})^{1-C(1+\beta)})$. For the range $1\leq k\leq K$, we use Siegel-Walfitz  to get a bound $O(X\log{X})^{-A}(1+ (\log{X})^{-C\beta})$. Choosing $C$ appropriately, say, $C= \frac{A+1}{\beta +1}$ if $\beta\geq 0$ and $C=A+1$ if $\beta <0$ gives an error term of the form $X(\log{X})^{-A'}$ for $A'$ suitably large. Collecting the estimates for the three terms gives the result for $\rho=\rho^{}_0$.

It remains to treat the case where $\rho\neq\rho^{}_0$ is nonprincipal. As the above analysis yields the same type of upper bound on the error terms, we need only concern ourselves with the main term of the asymptotics. The main term contribution to (\ref{3a}) for fixed $k$ and fixed $\ve$ is $\frac{1}{2^\gamma}$ times
    \[\overline{\chi(2)}\ell^{-k\beta}\sum_{a\in(\Z/\ell\Z)^*}\chi(a)\sum_{\substack{p\leq X\\ p\equiv a\ell^k-\ve\modulo \ell^{k+1}}}\rho(p)\log p.\]
    If the character $\rho$ is nonconstant along the arithmetic progression $a\ell^k-\ve+\Z\ell^{k+1}$ for every $a\in(\Z/\ell\Z)^*$, then by nontriviality of $\rho$ modulo $N/\gcd(\ell^{k+1},N)$ the above sum is zero. So we suppose the character $\rho$ is constant along $a\ell^k-\ve+\Z\ell^{k+1}$ for some $a\in(\Z/\ell\Z)^*$. This implies that $\rho$ is induced from a character modulo $\gcd(\ell^{k+1},N)$. Since $\rho$ is quadratic and nontrivial, this implies that $\rho$ is induced from the Legendre symbol $(-/\ell)$. Then the leading coefficient of $X$ in the main term contribution to (\ref{3a}) becomes $\frac{1}{2^\gamma}$ times
    \[\left(1+(-1)^j\left(\frac{-1}{\ell}\right)\right)\delta_\chi (\ell^{(1+\beta)}-1)^{-1}+\frac{1}{\ell-1}\overline{\chi(2)}\sum_{\ve=\pm1}\sum_{a\in(\Z/\ell\Z)^*}\chi(a+\epsilon)\left(\frac{a}{\ell}\right).\]
    We need to show that the above is strictly smaller in absolute value than
    \[\alpha=\frac{2}{2^\gamma}\left(\frac{\ell-2}{\ell-1}+\frac{1}{\ell^{1+\beta}-1}\right).\]
    This is clear if $\chi\neq\chiz$. If $\chi=\chiz$, note for example $|\sum_{\ve\neq a\in(\Z/\ell\Z)^*}(a/\ell)|<\ell-2$ for each $\epsilon=\pm1$, since $\ell\geq5$. The desired result follows.
\end{proof}

\begin{proposition}\label{pre3}
Let $(\chi,\rho,j)$ be given, and let $\ve\in\{\pm1\}$.
    \begin{enumerate}[label=\upshape\arabic*.\ ]
        \item There exists $k\in\N$ with
        \[\sum_{p\leq X}\frac{|f_{\chi,\rho,j,\ve}(p)|}{p}\leq k\log\log x+O(1)\]
        \item
        \[\sum_{p\leq x}\frac{|f_{\chi,\rho,j,\ve}(p)|\log p}{p}\ll\log x.\]
        \item Assume that $\gamma=0$. Then there exists some small $\delta>0$ such that
        \[\sum_{p}\sum_{k=2}^\infty\left|\frac{f_{\chi,\rho,j,\ve}(p)}{p^{1-\delta}}\right|^k<\infty.\]
    \end{enumerate}
\end{proposition}

\begin{proof}
For simplicity, we shall write $f(p)=f_{\chi,\rho,j,\ve}(p)$ for each prime $p$. Since $(1+\ve/p)^\gamma=1+\gamma\ve/p+O_\gamma(p^{-2})$, it is enough to establish the result for $\gamma=0$, which we shall now assume. Let $L=[\log X/\log\ell]$. Let $\theta(x) = 1 $ or $\log{x}$, and consider for $\ve = \pm 1$ \[S(X;\theta):=\sum_{p\leq X}\frac{|p+\ve|_\ell^\beta\theta(p)}{p} \ll \sum_{t=0}^{L}\ell^{t|\beta|}\sum_{\substack{p\leq X\\p\equiv \ve \tmod{\ell^t}}}\frac{\theta(p)}{p}.\] 
When $t=0$, the desired bounds in (1) and (2) are satisfied. For $t\geq 1$, we apply Prop.~\ref{prelim-3} and  Cor.~\ref{prelim-5}. We have \[S(X,1)= \log\log{X}\left(\sum_{t=1}^L \frac{\ell^{t|\beta|}}{\phi(\ell^t)}\right) + O\left(\sum_{t=1}^L \frac{\ell^{t|\beta|}}{\phi(\ell^t)}\log{\ell^t}\right) + \sum_{t=1}^L \ell^{t|\beta|}F(X,\ell^t,\ve).\] The first and last term gives a bound $O(\log\log{X})$. To estimate the last term, we consider bounds for sums of the type $T(u):=\sum_{t=1}^L \ell^{t|\beta|}E(u,\ell^t,\ve)$ with $3\leq u \leq X$. If $t\geq \frac{\log{u}}{\log{\ell}}$, $E(u,\ell^t,\ve)= -\text{Li}(u)\phi(\ell^t)^{-1}$ so that the contribution to $T(u)$ is $O(\text{Li}(u)u^{-(1-|\beta|)})$. For $A'\frac{\log\log{u}}{\log{\ell}}\leq t\leq \frac{\log{u}}{\log{\ell}}$ we use Brun-Titchmarsh (Prop.~\ref{prelim-1}) to get a bound of $u(\log{u})^{-A'(1+\beta)}$, while for the remaining range, we use Siegel-Walfitz to get the bound $u(\log{u})^{-A+A'|\beta|}$.  Choosing $A'$ appropriately and combining these estimates gives $T(u)\ll \frac{u^{|\beta|}}{\log{u}} + u(\log{u})^{-C}$ for $C$ a sufficiently large constant. Substituting these estimates into the third sum in $S(X,1)$ above and using the definition of $F(X,\ell^t,\ve)$ gives a bound of $o(1)$ as $X$ grows. This proves (1). For (2), we use the same analysis as the case above, using Cor.~\ref{prelim-5} to obtain that $S(X,\log)$ is bounded by 
\[ \sum_{t=0}^L \frac{\log{X}}{\ell^{t(1-|\beta|)}} + \left(\frac{1}{X^{1-|\beta|}} + \frac{1}{(\log{X})^{C}}\right)+ \int_3^X\log{\frac{u}{e}}  \left(\frac{u^{|\beta|}}{u^{2}\log{u}} + \frac{1}{(\log{u})^{C}}\right)du,\] which is $O(\log{X})$ as required. It thus remains to prove (3). The claim is clear if $\beta\geq0$, so we shall assume $-1<\beta<0$. We claim that it suffices to take $1-\delta=\sigma>\max(|\beta|,\frac{1}{2})$. We have
\[\sum_{p}\sum_{k=2}^\infty\frac{|f(p)|^k}{p^{k\sigma}}=\sum_{t=0}^\infty\sum_{k=2}^\infty\sum_{\substack{p,\\\ell^t||p+\ve}}\frac{\ell^{-t\beta k}}{p^{k\sigma}}.\]
Extracting $t=0$ separately, that sum is
\[\ll\sum_{k\geq2}\sum_p\frac{1}{p^{k\sigma}}\ll\sum_p\frac{1}{p^{2\sigma}}<\infty.\]
For $t\geq1$, we write $p+\ve=\ell^t m$ (and necessarily $m\geq2$). Also, $p\gg\ell^tm$, so the sum is
\[\ll\sum_{t\geq1}\sum_{k\geq2}\ell^{-t(\sigma-\beta)k}\sum_m\frac{1}{m^{k\sigma}}\ll\zeta(2\sigma)\sum_{k\geq2}\ell^{-k(\sigma-\beta)}\ll1\]
since $\sigma>\max(\frac{1}{2},|\beta|)$. This completes the proof.
\end{proof}

\begin{corollary}\label{pre4}
Fix a prime $\ell\geq5$. Let $-1<\beta,\gamma\leq\infty$. Let $c\geq1$. Then there exist $\delta>0$ and $B>0$ such that
\[\sum_{p> B}\sum_{k=2}^\infty\frac{(c\sum_{\ve=\pm1}|p+\ve|_\ell^\beta|\frac{p+\ve}{2p}|_\infty^\gamma)^k}{p^{k(1-\delta)}}<\infty.\]
\end{corollary}

\begin{proof}
We note that $c\sum_{\ve\pm1}|p+\ve|_\ell^\beta|\frac{p+\ve}{2p}|_\infty^\gamma\leq 2\cdot 4^{\max(-\gamma,0)}c\max_{\ve}|p+\ve|_\ell^\beta$. By Prop.\ref{pre3}, we can find $\delta'>0$ such that, for every $\ve=\pm1$,
\[\sum_p\sum_{k\ge2}\frac{|p+\ve|_\ell^{\beta k}}{p^{k(1-\delta')}}<\infty.\]
We choose $B>0$ such that
\[\frac{p^{1-\delta'/2}}{2\cdot 4^{\max(-\gamma,0)}c}\geq p^{1-\delta'},\]
for all $p\geq B$. Then we have
\[\sum_{p> B}\sum_{k\geq2}\frac{(c\sum_{\ve=\pm1}|p+\ve|_\ell^\beta|\frac{p+\ve}{2p}|_\infty^\gamma)^k}{p^{k(1-\delta'/2)}}\leq\sum_{\ve=\pm1}\sum_{p>B}\sum_{k\geq2}\frac{|p+\ve|_\ell^{\beta k}}{p^{k(1-\delta')}}<\infty,\]
which gives the result, setting $\delta=\delta'/2$.
\end{proof}

\begin{corollary}[\cite{GranK}, Sec.~2]\label{pre5}
Let $c\geq1$, and choose $B>0$ so that the conclusion of Corollary \ref{pre4} holds. Then the Dirichlet series
    \[D(s)=\prod_{p>B}\left(1+\frac{c\cdot f_{\chi,\rho,j}(p)}{p^s}\right),\]
has a holomorphic extension to $\Re(s)\geq1$ except for a singularity at $s=1$, where
\[ (s-1)^{c\cdot \alpha_{\rho,\chi,j}}D(s) = A +O(|s-1|^\theta),\] for some $0<\theta \leq 1$ and some constant $A\neq 0$, valid in a neighborhood of $s=1$.
\end{corollary}

\section{Asymptotics for prime determinants}\label{primedet}

Throughout this section, we fix the following data:
\begin{enumerate}[label=\upshape\arabic*.\ ]
    \item a prime $\ell\geq5$;
    \item $\beta\in(-1,\infty]$ and $\beta^*=\max(-\beta,0) <1$;
    \item $\gamma\in(-1,\infty)$.
\end{enumerate}
As in the introduction, we define
\[\alpha(\ell,\beta)=1-\frac{1}{\ell-1}+\frac{1}{\ell^{1+\beta}-1}.\]

\begin{theorem}\label{thm-prime}
Let $G$ run through the genera of primitive integral ternary quadratic forms of positive determinant. 
\begin{enumerate}[label=\upshape\arabic*.\ ]
\item For each $\varpi\in\{\pm1\}$,
\[\sum_{\substack{c_\infty(G)=\varpi\\D(G)\leq X\text{ prime}}}|\mu^*(G)|_\ell^\beta|\mu^*(G)|_\infty^\gamma\sim\frac{1}{1+\gamma}8^{-\gamma}(3^\gamma +1)\alpha(\ell,\beta)\frac{X^{1+\gamma}}{\log X}.\]
Moreover, for each $\varpi\in\{\pm1\}$ and $(a,\ell)=1$, we have
\[\sum_{\substack{c_\infty(G)=\varpi\\\mu^*(G)_{\ell'}\equiv a\tmod{\ell}\\D(G)\leq X\text{ prime}}}|\mu^*(G)|_\ell^\beta|\mu^*(G)|_\infty^\gamma\sim\frac{1}{1+\gamma}\frac{8^{-\gamma}}{(\ell -1)}\left[\alpha(\ell,\beta,\gamma) - \sigma(a,\ell,\gamma,\varpi)\right]\frac{X^{1+\gamma}}{\log X},\]
where
\[\alpha(\ell,\beta,\gamma)= (3^\gamma +1)\left(1+ \frac{1}{\ell^{1+\beta}-1}\right),\]and
\[\sigma(a,\ell,\gamma,\varpi):=   \begin{cases} 3^\gamma, & \text{if }\; 8a\equiv 3\varpi \tmod{\ell},\\ 1, & \text{if }\; 8a\equiv -\varpi \tmod{\ell},\\ 0, & \text{otherwise}.
    \end{cases}\]
\item For each $\varpi\in\{\pm1\}$
\[\sum_{\substack{c_\infty(G)=\varpi\\\mu^*(G)\leq X, D(G)\text{ prime}}}|\mu^*(G)|_\ell^\beta\sim 9\,\alpha(\ell,\beta)\frac{X}{\log X}.\]
Moreover, for each $\varpi\in\{\pm1\}$ and $(a,\ell)=1$, we have
\[\sum_{\substack{c_\infty(G)=\varpi\\\mu^*(G)_{\ell'}\equiv a\tmod{\ell}\\\mu^*(G)\leq X, D(G)\text{ prime}}}|\mu^*(G)|_\ell^\beta\sim \frac{1}{(\ell -1)}\left[9\left(1+\frac{1}{\ell^{1+\beta}-1}\right) - \eta(a,\ell,\varpi)\right]\frac{X}{\log X},\]
where
\[\eta(a,\ell,\varpi):=   \begin{cases} 8, & \text{if }\; 8a\equiv 3\varpi \tmod{\ell},\\ 1, & \text{if }\; 8a\equiv -\varpi \tmod{\ell},\\ 0, & \text{otherwise}.
    \end{cases}\]    
    \end{enumerate}
\end{theorem}
\begin{remark}
    Setting $\gamma=0$ and comparing the cases $\beta=0$ and $\beta=1$ in the above theorem, we confirm the interpretation (\ref{alphameaning}) of $\alpha(\ell,\beta)$ as the average order of the $\beta$th power of the $\ell$-adic absolute value of $\mu^*(G)$ among the genera of prime determinant in $\mathcal F$.
\end{remark}
\begin{proof}\ 

If $D(G)=p>2 $ is a prime, then $G$ is determined by a pair of local factors we denote by $G_2^\pm$ and $G_p^\pm$, where the $\pm$ denotes the value of the local Hasse invariant. From our prior data, we have $\mu^*(G_2^+)= \frac{3}{4}$, $\mu^*(G_2^-)= \frac{1}{4}$, $\mu^*(G_p^+)= \half(p+1)$ and $\mu^*(G_p^-)= \half(p-1)$, giving four choices of $G$ for each $p$ (which we will assume exceeds $\ell$).

\noindent \textit{Case 1:}\quad 

We write
 \[\sum_{\substack{c_\infty(G)=\varpi\\\mu^*(G)_{\ell'}\equiv a\tmod\ell \\D(G)= \leq X,\;\text{prime} }}|\mu^*(G)|_\ell^\beta\left|\frac{\mu^*(G)}{D(G)}\right|_\infty^\gamma = \frac{1}{2(\ell -1)}\sum_{\chi \tmod{\ell}}\overline{\chi(a)}\sum_{j=0}^1\varpi^jS_j(X,\chi),\]
 where 
\[\begin{split} S_j(X,\chi)&= \sum_{\substack{D(G)\leq X \\D(G)=\text{prime} }}c_\infty(G)^j\chi(\mu^*(G)_{\ell'})|\mu^*(G)|_\ell^\beta\left|\frac{\mu^*(G)}{D(G)}\right|_\infty^\gamma,\\
&= \sum_{\substack{p\leq X \\p=\text{prime} }}\sum_{G\in \Gen(p)}c_\infty(G)^j\chi(\mu^*(G)_{\ell'})|\mu^*(G)|_\ell^\beta\left|\frac{\mu^*(G)}{D(G)}\right|_\infty^\gamma. \end{split}\]
By multiplicativity, the sum over genera is a product coming from the primes two and $p$; the contribution from the former independent of $p$.

With $f_{\chi,\rho,j}(p)$ be as in Def.~\ref{def3.1}, where we take  $\rho =1$, we have 
\[ S_j(X,\chi) = \chi(\tfrac{1}{4})(\tfrac{1}{4})^\gamma\left[3^\gamma\chi(3)+(-1)^j\right]\sum_{p\leq X}f_{\chi,1,j}(p) +O(1).  \] The proof of Prop.~\ref{pre1} gives \[\sum_{p\leq X}f_{\chi,1,j}(p)=\alpha_{\chi,j} \frac{X}{\log X}+O\left(\frac{X}{(\log X)^2}\right),\]
where 
\[ \alpha_{\chiz,j}= 2^{1-\gamma}\delta_j\alpha(\ell,\beta),\]
and for non-principal $\chi$,
\[\alpha_{\chi,j} =- 2^{-\gamma}\chi(\tfrac{1}{2})\left[1+(-1)^j\chi(-1)\right]\frac{1}{\ell -1}. \]
    Then, \[ \sum_{j=0}^1\varpi^j\left[3^\gamma\chi(3)+(-1)^j\right]\alpha_{\chi,j} = \begin{cases} 2^{1-\gamma}\left[3^\gamma+1\right]\alpha(\ell,\beta), &  \text{if } \chi=\chiz,\\[6pt] -\frac{2^{1-\gamma}}{\ell -1}\left[3^\gamma\chi(\tfrac{3\varpi}{2})+ \chi(\tfrac{-\varpi}{2})\right] , &  \text{if } \chi\neq\chiz. \end{cases}\] 
    Summing over the characters gives us
    \[2^{\gamma -1}\sum_\chi \overline{\chi(4a)}\sum_{j=0}^1\varpi^j\left[3^\gamma\chi(3)+(-1)^j\right]\alpha_{\chi,j}= \alpha(\ell,\beta,\gamma) - \sigma(a,\ell,\beta,\gamma) \]
    with \[\alpha(\ell,\beta,\gamma)= (3^\gamma +1)\left(1+ \frac{1}{\ell^{1+\beta}-1}\right),\]
    and
    \[ \sigma(a,\ell,\gamma,\varpi):=   \begin{cases} 3^\gamma, & \text{if }\; 8a\equiv 3\varpi \tmod{\ell},\\ 1, & \text{if }\; 8a\equiv -\varpi \tmod{\ell},\\ 0, & \text{otherwise}.
    \end{cases}\]
    We then get
    \[ \lim_{X\to \infty} \frac{\log X}{X}\smashoperator[r]{\sum_{\substack{c_\infty(G)=\varpi\\\mu^*(G)_{\ell'}\equiv a\tmod\ell \\D(G)= \leq X,\;\text{prime} }}}|\mu^*(G)|_\ell^\beta\left|\frac{\mu^*(G)}{D(G)}\right|_\infty^\gamma = \frac{8^{-\gamma}}{(\ell -1)}\left(\alpha(\ell,\beta,\gamma) - \sigma(a,\ell,\gamma,\varpi)\right). \]
    If we sum over the reduced classes $a \pmod{\ell}$, we get
    \[ \lim_{X\to \infty} \frac{\log X}{X}\sum_{\substack{c_\infty(G)=\varpi\\D(G)= \leq X,\;\text{prime} }}|\mu^*(G)|_\ell^\beta\left|\frac{\mu^*(G)}{D(G)}\right|_\infty^\gamma = 8^{-\gamma}(3^\gamma +1)\alpha(\ell,\beta). \] The theorem then follows.

    \noindent \textit{Case 2:}\quad The proof follows that given above. Let 
\[ \tilde{S}_j(X,\chi)= \sum_{\substack{\mu^*(G)\leq X \\D(G)=\text{prime} }}c_\infty(G)^j\chi(\mu^*(G)_{\ell'})|\mu^*(G)|_\ell^\beta.\]
Then,
 \[\sum_{\substack{c_\infty(G)=\varpi,\,\mu^*(G)\leq X\\\mu^*(G)_{\ell'}\equiv a\tmod\ell \\D(G)=\text{prime} }}|\mu^*(G)|_\ell^\beta = \frac{1}{2(\ell -1)}\sum_{\chi}\overline{\chi(a)}\sum_{j=0}^1\varpi^j\tilde{S}_j(X,\chi).\]
We have $\mu^*(G)$ being of the form $\tfrac{3}{8}(p\pm 1)$ or $\tfrac{1}{8}(p\pm 1)$, with the restriction that they not exceed $X$. The $\pm 1$ occurring in the interval boundary for the prime restrictions give an error of $O(X^{\beta^*})$, with $\beta^*=\max(-\beta,0) <1$, so that we can disregard them.

    Then,
    \[  \tilde{S}_j(X,\chi) =\chi(\tfrac{3}{4})T_j(\tfrac{8}{3}X,\chi)+ \chi(\tfrac{1}{4})(-1)^jT_j(8X,\chi) +O(X^{\beta^*}),\]
    where $ T_j(Y,\chi)= \sum_{p\leq Y}f_{\chi,1,j}(p)$ as above. Following the same analysis as above, after substituting for the asymptotics of $T_j(Y,\chi)$ with $\gamma=0$, we get
    \[ \tilde{S}_j(X,\chi) = \tfrac{1}{3}\left[8\chi(\tfrac{3}{4})+(-1)^j\chi(\tfrac{1}{4})\right]\alpha_{\chi,\ell}\frac{X}{\log X} + O(X(\log X)^{-2}).  \] Following the analysis above,  summing over $j$ and  $\chi$ gives
    \[\lim_{X\to \infty}\frac{\log X}{X}\smashoperator[r]{\sum_{\substack{c_\infty(G)=\varpi,\,\mu^*(G)\leq X\\\mu^*(G)_{\ell'}\equiv a\tmod\ell \\D(G)=\text{prime} }}}|\mu^*(G)|_\ell^\beta = \frac{1}{2(\ell -1)}\left[18\left(1+\frac{1}{\ell^{1+\beta}-1}\right) - \eta(a,\ell,\varpi)\right], \]where
    \[\eta(a,\ell,\varpi):=   \begin{cases} 16, & \text{if }\; 8a\equiv 3\varpi \tmod{\ell},\\ 2, & \text{if }\; 8a\equiv -\varpi \tmod{\ell},\\ 0, & \text{otherwise}.
    \end{cases}\] the result then follows.
\end{proof}

\section{Proofs of Theorems \ref{thm1n} - \ref{thm4n}}\label{pf-thms}
\subsection{Proof of Theorem \ref{thm1n}}\label{pf-thm1n}\
Throughout this subsection, we fix the following data:
\begin{enumerate}[label=\upshape\arabic*.\ ]
    \item a prime $\ell\geq5$;
    \item $\beta,\gamma\in(-1,\infty]$ such that $\min(\beta,0)+\min(\gamma,0)>-1$;
    \item a constant $B>\ell$ so that the conclusion of Corollary \ref{pre4} holds with $c=1$.
\end{enumerate}
Let
\[\alpha(\ell,\beta)=1-\frac{1}{\ell-1}+\frac{1}{\ell^{1+\beta}-1}\]
as before.

\begin{proposition}\label{pointwise}
Fix a packet $H\in\Packet$ such that $\{p:p\leq B\}\subseteq\Supp(H)$. Then for each $\varpi\in\{\pm1\}$ and $a\in(\Z/\ell\Z)^*$ we have
\[\sum_{\substack{G\in\Gen_H^\varpi,D(G)\leq X\\\mu^*(G)_{\ell'}\equiv a\modulo \ell}}|\mu^*(G)|_\ell^\beta\left|\frac{\mu^*(G)}{D(G)}\right|_\infty^\gamma\sim\frac{1}{2}\frac{1}{\ell-1}C_{\ell,\beta,\gamma,H} X(\log X)^{\alpha-1},\]
where
\[\begin{split}C_{\ell,\beta,\gamma,H}=\frac{1}{\Gamma(\alpha)}&\prod_{p\in\Supp(H)}\frac{1}{|\Z_p/(\Z_p^*)^2|}\frac{|\mu_p^*(H_p)|_\ell^\beta}{D_p(H)}\left(\frac{\mu_p^*(H_p)}{D_p(H)}\right)^\gamma\left(1-\frac{1}{p}\right)^\alpha  \\ &\times \prod_{p\notin\Supp(H)}\left(1+\frac{1}{p}\sum_{\ve=\pm1}|p+\ve|_\ell^\beta\left(\frac{p+\ve}{2p}\right)^\gamma\right)\left(1-\frac{1}{p}\right)^{\alpha}.\end{split}\]
\end{proposition}

\begin{proof}
Let $X(H)$ denote the set of quadratic characters of $\prod_{p\in\Supp(H)}\Z_p/(\Z_p^*)^2$. For each Dirichlet character $\chi$ mod $\ell$, $\rho\in X(H)$, and $j\in\{0,1\}$, we write for simplicity
\[A_{\ell,\beta,\gamma,H,\chi,\rho,j}(s)=\prod_{p\notin\Supp(H)}\left(1+\frac{1}{p^s}\sum_{\ve=\pm1}\ve^j\chi\left(\frac{p+\ve}{2}\right)\rho(p)|p+\ve|_\ell^\beta\left(\frac{p+\ve}{2p}\right)^\gamma\right).\]
We write $[H]$ for the class $(D(H_p)/D(H))\in \prod_{p\in\Supp(H)}\Z_p/(\Z_p^*)^2$. Then
\[\begin{split}&\sum_{\substack{G\in\Gen_H^{\varpi}\\\mu^*(C)_{\ell'}\equiv a\modulo \ell}}\frac{|\mu^*(G)|_\ell^\beta}{D(G)^s}=\frac{1}{2}\frac{1}{|X(H)|}\sum_{j=0}^1\sum_{\chi\modulo\ell}\overline{\chi(a)}\\ & \hspace{4em}\times \sum_{\rho\in X(H)}\varpi^j\frac{c_\infty(H)|\mu^*(H)|_{\ell}^\beta\chi(D(H))\rho([H])}{D(H)^s}\left|\frac{\mu^*(H)}{D(H)}\right|_\infty^\gamma A_{\ell,\beta,\gamma,H,\rho,j}(s).\end{split}\]
The result then follows by applying Proposition \ref{pre1}, Proposition \ref{pre3}(1),(2), Corollary \ref{pre4} (which is applicable by our choice of $B$ and $H$), and Prop.~\ref{simple} to each of the Eulerian Dirichlet series $A_{\ell,\beta,\gamma,H,\rho,j}(s)$ and adding up their contributions, with the leading term given by $\rho=\rho^{}_0$ and $j=0$.
\end{proof}

\begin{proposition}\label{dominator}
    There exists a constant $C=c_{\ell,\beta,\gamma,B}>0$ such that, for any packet $H\in\Packet$ such that $\{p:p\leq B\}\subseteq\Supp(H)$, and for any $X\geq2$,
    \[\frac{1}{X(\log X)^{\alpha-1}}\smashoperator[r]{\sum_{\substack{S\in\SGen_H,\\D(S)\leq X}}}|\mu^*(S)|_\ell^\beta\Bigl|\frac{\mu^*(S)}{D(S)}\Bigr|_\infty^\gamma\leq C(2^{2+\omega_3(D(H))})^{1-\gamma}\frac{|\mu^*(H)|_\ell^\beta}{D(H)}\Bigl|\frac{\mu^*(H)}{D(H)}\Bigr|_\infty^\gamma.\]
\end{proposition}

\begin{proof}
By Proposition \ref{equal-mass} and the fact that $\sg(G)$ is a power of $2$, we have
    \[\sum_{S\in\SGen(G)}|\mu^*(S)|_\ell^\beta\left|\frac{\mu^*(S)}{D(S)}\right|_\infty^\gamma= \sg(G)^{1-\gamma}|\mu^*(G)|_\ell^\beta\left|\frac{\mu^*(G)}{D(G)}\right|_\infty^\gamma\]
    for every $G\in\Gen$. Since $\sg(G)\leq 2^{2+\omega_3(D(H))}$ for $G\in\Gen_H$, to prove the proposition it suffices to show the existence of $c_{\ell,\beta,\gamma,B}'>0$ such that, for any packet $H\in\Packet$ such that $\{p:p\leq B\}\subseteq \Supp(H)$, and for any $X\geq2$,
    \[\frac{1}{X(\log X)^{\alpha-1}}\sum_{\substack{G\in\Gen_H\\D(G)\leq X}}|\mu^*(G)|_\ell^\beta\left|\frac{\mu^*(G)}{D(G)}\right|_\infty^\gamma\leq c_{\ell,\beta,\gamma,B}'\frac{|\mu^*(H)|_\ell^\beta}{D(H)}\left|\frac{\mu^*(H)}{D(H)}\right|_\infty^\gamma.\]
    Consider now the Dirichlet series
    \[\sum_{d=1}^\infty\frac{g(d)}{d^s}=\prod_{p>B}\left(1+\frac{1}{p^s}\sum_{\ve=\pm1}|p+\ve|_\ell^\beta\left|\frac{p+\ve}{2p}\right|_\infty^\gamma\right).\]
    Then we see that
    \[\sum_{\substack{G\in\Gen_H\\D(G)\leq X}}|\mu^*(G)|_\ell^\beta\left|\frac{\mu^*(G)}{D(G)}\right|_\infty^\gamma\leq |\mu^*(H)|_\ell^\beta\left|\frac{\mu^*(H)}{D(H)}\right|_\infty^\gamma\sum_{d\leq X/D(H)}g(d).\]
    Applying Proposition \ref{pre1}, Proposition \ref{pre3}(1),(2), Corollary \ref{pre4} (which is applicable by our choice of $B$ and $H$), and Prop.~\ref{simple}, we deduce that
    \[\sum_{d\leq X}g(d)\sim c_{\ell,\beta,\gamma,B} X(\log X)^{\alpha-1}\]
    for some constant $c_{\ell,\beta,\gamma,B}>0$. The desired result then follows.
\end{proof}

\begin{proposition}\label{splitting}
Fix a packet $H\in\Packet$ with $\{p:p\leq B\}\subseteq\Supp(H)$. Then
\[\sum_{\substack{G\in\Gen_H,\\D(G)\leq X\\\sg(G)\geq 2}}|\mu^*(G)|_\ell^\beta\left|\frac{\mu^*(G)}{D(G)}\right|_\infty^\gamma\ll_H X(\log X)^{\frac{1}{2}\alpha-1}.
\]
\end{proposition}

\begin{proof}
We first assume that $-1<\beta<0$. In what follows, we write
\[f(p)=\sum_{\ve=\pm1}|p+\ve|_\ell^\beta\left|\frac{p+\ve}{2p}\right|_\infty^\gamma\]
for simplicity. Below, given $G\in\Gen_H$ let us write $s(G)=\prod_{p\mid D(G), p\notin\Supp(H)}p$. If $G\in\Gen_H$ satisfies $\sg(G)\geq2$, this implies by Prop.~\ref{propE1} that there exists a proper subgroup $\mathcal C\leq(\Z/8r(H)\Z)^*$ of index $2$ such that $p\in\mathcal C$ for any $p\mid s(G)$. Since there are at most finitely many subgroups $\mathcal C$, it suffices to fix one such $\mathcal C$ and show that
\[\sum_{\substack{G\in\Gen_H,\\D(G)\leq X,\\p\mid s(G)\Rightarrow p\in\mathcal C}}|\mu^*(G)|_\ell^\beta\left|\frac{\mu^*(G)}{D(G)}\right|_\infty^\gamma\ll_{\mathcal C} X(\log X)^{\frac{1}{2}\alpha-1}.\]
As in the proof of Prop.~\ref{pointwise}, the result follows by applying Proposition \ref{pre1}, Proposition \ref{pre3}(1),(2), Corollary \ref{pre4} (which is applicable by our choice of $B$ and $H$), and Prop.~\ref{simple}.

\end{proof}

\begin{proposition}\label{ellone}
Given any packet $H$, define
\[N(H)=(2^{2+\omega_3(D(H))})^{1-\gamma}\frac{|\mu^*(H)|_\ell^\beta}{D(H)}\left|\frac{\mu^*(H)}{D(H)}\right|_\infty^\gamma.\]
Then $\sum_{H\in\mathfrak R_B} N(H)<\infty$.
\end{proposition}
\begin{proof}
Let $Q$ denote the set of root packets $H\in\Packet$ such that $\ell\notin\Supp(H)$. Then
\[\sum_{H\in\mathfrak R_B}N(H)\leq 2^{4+B}\prod_{p\leq B}\left(\sum_{G\in\Gen_p}\frac{|\mu_p^*(G)|_\ell^\beta}{D_p(G)}\left|\frac{\mu_p^*(G)}{D_p(G)}\right|_\infty^\gamma\right)\cdot\sum_{H\in Q}N(H).\]
In light of Corollary \ref{prelim-7a}, it suffices to show that $\sum_{H\in Q}N(H)<\infty$.
But we have
\begin{align*}
\sum_{H\in Q} N(H)&\leq 4^{1-\gamma}\prod_{p\neq \ell}\left(1+2^{1-\gamma}\sum_{\substack{G\in\Gen_p\\v_p(D(G))\geq2}}\frac{|\mu_p^*(G)|_\ell^\beta}{D_p(G)}\left|\frac{\mu_p^*(G)}{D_p(G)}\right|_\infty^\gamma\right)\\
&\leq4^{1-\gamma}\exp\left(2^{1-\gamma}\sum_{p\neq\ell}\sum_{\substack{G\in\Gen_p\\v_p(D(G))\geq2}}\frac{|\mu_p^*(G)|_\ell^\beta}{D_p(G)}\left|\frac{\mu_p^*(G)}{D_p(G)}\right|_\infty^\gamma\right).
\end{align*}
Again by Cor. \ref{prelim-7a}, writing $\beta^*=\max(-\beta,0)$ and $\gamma^*=\max(-\gamma,0)$, we have
\begin{align*}
\sum_{p\neq\ell}\sum_{\substack{G\in\Gen_p\\v_p(D(G))\geq2}}\frac{|\mu_p^*(G)|_\ell^\beta}{D_p(G)}\left|\frac{\mu_p^*(G)}{D_p(G)}\right|_\infty^\gamma&\leq\sum_{p\neq\ell}p^{\delta}\sum_{k=2}^\infty\frac{8k}{(p^{1-\gamma^*/2})^k}\\
&=8\sum_{p\neq\ell}p^{\delta}\frac{2p^{1-\gamma^*/2}-1}{p^{1-\gamma^*/2}(p^{1-\gamma^*/2}-1)^2}<\infty
\end{align*}
showing that $\sum_{H\in Q}N(H)<\infty$ as desired.
\end{proof}

We now give a more precise statement of Theorem \ref{thm1n} and prove it.

\begin{thm}{A$'$}
Let $\ell\geq5$ be prime. Let $\beta,\gamma\in(-1,\infty]$ with
$\min(\beta,0)+\min(\gamma,0)>-1$.
Let $S$ run through the spinor genera of primitive integral ternary quadratic forms of positive determinant, and let $G$ run through the genera of such forms. Let
\begin{equation}\label{exp-clb}
C_{\ell,\beta,\gamma}=\frac{1}{\Gamma(2^{1-\gamma}\alpha(\ell,\beta))}\prod_{p}\Bigl(\mathop{\sum\sum}_{\substack{k\geq 0\\ G\in\Gen_p(p^k)}}\frac{|\mu_p^*(G_p)|_\ell^\beta}{D_p(G_p)}\left|\frac{\mu_p^*(G_p)}{D_p(G_p)}\right|_\infty^\gamma\Bigr)\left(1-\frac{1}{p}\right)^{2^{1-\gamma}\alpha(\ell,\beta)}.
\end{equation}
Then $C_{\ell,\beta,\gamma}\in\R_{>0}$ and
    \[\begin{split} \sum_{\substack{D(G)\leq X}}|\mu^*(G)|_\ell^\beta|\mu^*(G)|_\infty^\gamma &\sim \sum_{\substack{D(S)\leq X}}|\mu^*(S)|_\ell^\beta|\mu^*(S)|_\infty^\gamma \,, \\ &\sim\frac{C_{\ell,\beta,\gamma}}{1+\gamma}~X^{1+\gamma}(\log X)^{2^{1-\gamma}\alpha(\ell,\beta)-1}.\end{split}\]
Moreover for any $\varpi\in\{\pm1\}$ and $(a,\ell)=1$,
    \[\sum_{\substack{c_\infty(S)=\varpi,\,D(S)\leq X\\\mu^*(S)_{\ell'}\equiv a\tmod\ell}}|\mu^*(S)|_\ell^\beta|\mu^*(S)|_\infty^\gamma\sim\frac{1}{2}\frac{1}{\ell -1}~\sum_{\substack{D(S)\leq X}}|\mu^*(S)|_\ell^\beta|\mu^*(S)|_\infty^\gamma,\]
and the same asymptotic holds if one sums over genera $G$. In particular, when ordered by the determininant, the typical $\mu^*(S)$ are divisible by arbitrarily large powers of $\ell$. Moreover, the prime-to-$\ell$ parts of $\mu^*(S)$ are uniformly distributed in $(\Z/\ell\Z)^*$, and the same holds true if we restrict to those $S$ such that $(\mu^*(S),\ell)=1$. 
\end{thm}

\begin{proof}
Recall that we have chosen $B>0$ so that the conclusion of Corollary \ref{pre4} holds with $c=1$. For each $X\geq2$, we consider $f_X:\mathfrak R_B\to\R_{\geq0}$ given by
    \[f_X(H)=\frac{1}{X(\log X)^{\alpha-1}}\sum_{\substack{S\in\operatorname{SGen}_H^\varpi,D(S)\leq X\\\mu^*(S)_{\ell'}\equiv a(\modulo \ell)}}|\mu^*(S)|_{\ell}^\beta\left|\frac{\mu^*(S)}{D(S)}\right|_\infty^\gamma\]
    and
    \[f_X^g(H)=\frac{1}{X(\log X)^{\alpha-1}}\sum_{\substack{G\in\operatorname{Gen}_H^\varpi,D(G)\leq X\\\mu^*(G)_{\ell'}\equiv a(\modulo \ell)}}|\mu^*(G)|_{\ell}^\beta\left|\frac{\mu^*(G)}{D(G)}\right|_\infty^\gamma\]
    so that what we need to show is
    \[\lim_{X\to\infty}\|f_X^g\|_1=\lim_{X\to\infty}\|f_X\|_1=\frac{1}{2}\frac{1}{\ell-1}C_{\ell,\beta}\in\R_{>0}.\]
    Let $N:\mathfrak R_B\to\R_{\geq0}$ be given by
    \[N(H)=2^{2+\omega_3(r(H))}\frac{|\mu^*(H)|_\ell^\beta}{D(H)}\left|\frac{\mu^*(H)}{D(H)}\right|_\infty^\gamma.\]
    Then we have $N\in\ell^1(\mathfrak R_B)$ by Prop.~\ref{ellone}. By Proposition \ref{dominator}, we have bounds $f_X^g(H)\leq f_X(H)\leq c_{\ell,\beta} N(H)$
    for every $H\in\mathfrak R_B$, where $c_{\ell,\beta}>0$ is is the constant given in Prop.~\ref{dominator}. By Propositions \ref{pointwise} and \ref{splitting} we have
    \[\lim_{X\to\infty}f_X^g(H)=\lim_{X\to\infty}f_X(H)=\frac{1}{2}\frac{1}{\ell-1}C_{\ell,\beta,\gamma,H}\]
    where the constant $C_{\ell,\beta,\gamma,H}$ is defined as in Proposition \ref{pointwise}. By the dominated convergence theorem, the pointwise limit $f:\mathfrak R_B\to\R_{\geq0}$ given by $f(H)=\lim_{X\to\infty} f_X(H)$ is integrable, and we have
    \begin{align*}
    \lim_{X\to\infty}\|f_X\|_1=\|f\|_1=\sum_{H\in\mathfrak R_B}f(H)=\frac{1}{2}\frac{1}{\ell-1}\sum_{H\in\mathfrak R_B} C_{\ell,\beta,H}=\frac{1}{2}\frac{1}{\ell-1}C_{\ell,\beta}
    \end{align*}
    which is the desired result.
\end{proof}

\subsubsection{The squarefree case:}\label{sqfree-det}\ 

We also include a proof of the analogous asymptotic formula when the genera $G$ vary over squarefree determinants.

\begin{theorem}\label{squarefree-det}
Let $\ell\geq5$ be prime. Let $\beta,\gamma\in(-1,\infty]$ with
$\min(\beta,0)+\min(\gamma,0)>-1$.
Let $G$ run through the genera of primitive integral ternary quadratic forms of positive, squarefree, determinant. For any $\varpi\in\{\pm1\}$ and $(a,\ell)=1$,
    \begin{align*}
    \sum_{\substack{c_\infty(G)=\varpi,\,D(G)\leq X\\\mu^*(G)_{\ell'}\equiv a\tmod\ell\\\text{$D(G)$ squarefree}}}|\mu^*(G)|_\ell^\beta|\mu^*(G)|_\infty^\gamma\sim\frac{1}{2}\frac{1}{\ell -1}\frac{c_{\ell,\beta,\gamma}}{1+\gamma}X^{1+\gamma}(\log X)^{2^{1-\gamma}\alpha(\ell,\beta)-1}
    \end{align*}
    where
    \[c_{\ell,\beta,\gamma}=\frac{1}{\Gamma(2^{1-\gamma}\alpha(\ell,\beta))}\prod_p\Bigl(\mathop{\sum\sum}_{\substack{k\in\{0,1\}\\G\in\Gen_p(p^k)}}\frac{|\mu_p^*(G)|_\ell^\beta}{|D_p(G)|}\left|\frac{\mu_p^*(G)}{D_p(G)}\right|_{\infty}^\gamma\Bigr)\left(1-\frac{1}{p}\right)^{2^{1-\gamma}\alpha(\ell,\beta)}.\]
In particular, when ordered by the determininant, the typical $\mu^*(S)$ are divisible by arbitrarily large powers of $\ell$. Moreover, the prime-to-$\ell$ parts of $\mu^*(S)$ are uniformly distributed in $(\Z/\ell\Z)^*$, and the same holds true if we restrict to those $S$ such that $(\mu^*(S),\ell)=1$. 
\end{theorem}
\begin{proof}
    The proof proceeds in as in the proof of Theorem \ref{thm1n}, except that one is only summing over the $B$-root packets $H\in\mathfrak R_B$ of squarefree determinant. We note that since there are only finitely many such $B$-root packets, no dominated convergence theorem is needed in this case.
\end{proof}

\subsection{Proof of Theorem \ref{thm3n}}\label{pfthm2n}\ 

The proof of Theorem \ref{thm3n} has the same underlying structure as the proof of Theorem \ref{thm1n}. But there are some changes to be made. Throughout this subsection, we fix the following data:
\begin{enumerate}[label=\upshape\arabic*.\ ]
    \item a prime $\ell\geq5$;
    \item $\beta,\gamma\in(-1,\infty]$ such that $\min(\beta,0)+\min(\gamma,0)>-1$;
    \item a constant $B>\ell$ so that the conclusion of Corollary \ref{pre4} holds with $c=2$.
\end{enumerate}
Let
\[\alpha(\ell,\beta)=1-\frac{1}{\ell-1}+\frac{1}{\ell^{1+\beta}-1},\]
as before.

\begin{definition}[\mbox{\mathversion{normal}$\mu^*$}-weighted  \mbox{\mathversion{normal}$A$}-root packets]
    Let $A$ be a finite set of primes containing $2$. A $\mu^*$-weighted $A$-root packet is a packet $H\in\Packet$ such that:
    \begin{enumerate}[label=\upshape\arabic*.\ ]
        \item $A\subseteq\Supp(H)$.
        \item For any prime $p\in\Supp(H)\setminus A$, we have $\mu_p^*(H_p)\notin\{1,\frac{p\pm1}{2}\}$.
    \end{enumerate}
    We write $\mathfrak R_A'$ for the set of all $\mu$-weighted $A$-root packets. If $A=\{p\text{ prime}:p\leq B\}$ for some real number $B\geq2$, we shall also write $\mathfrak R_A'=\mathfrak R_B'$ and say $\mu$-weighted $B$-root packets to mean $\mu$-weighted $A$-root packets.
\end{definition}

\begin{definition}
    Let $H\in\Packet$ be a packet. We write $\Gen_H'$ for the set of genera $G\in\Gen$ such that:
    \begin{enumerate}[label=\upshape\arabic*.\ ]
        \item $G_p=H_p$ for every $p\in\Supp(H)$, and
        \item $\mu_p^*(G_p)\in\{1,\frac{p\pm1}{2}\}$ for every $p\notin\Supp(H)$.
    \end{enumerate}
    We denote by $\SGen_H'$ the set of those spinor genera whose genera lie in $\Gen_H'$. Note that for any finite set $A$ we have
    \[\Gen=\coprod_{H\in\mathfrak R_{A}'}\Gen_H',\quad \SGen=\coprod_{H\in\mathfrak R_{A}'}\SGen_H'.\]
    We shall write $\Gen_H^{\varpi,v}=\Gen^\varpi\cap \Gen_H'$ and $\SGen_H^{\varpi,v}=\SGen^\varpi\cap \SGen_H'$.
\end{definition}

\begin{proposition}\label{pointwise-v}
Fix a packet $H\in\Packet$ such that $\{p:p\leq B\}\subseteq\Supp(H)$. Then there is a constant $C_{\ell,\beta,H}'$ such that
\[\sum_{\substack{G\in\Gen_H^\varpi,\mu^*(G)\leq X\\\mu^*(G)_{\ell'}\equiv a\modulo \ell}}|\mu^*(G)|_\ell^\beta\sim\frac{1}{2}\frac{1}{\ell-1}C_{\ell,\beta,H}' X(\log X)^{8\alpha-1},\]
for all each $\varpi\in\{\pm1\}$ and $a\in(\Z/\ell\Z)^*$.
\end{proposition}

\begin{proof}
Let $X(H)$ denote the set of quadratic characters of $\prod_{p\in\Supp(H)}\Z_p/(\Z_p^*)^2$. For each Dirichlet character $\chi$ mod $\ell$, $\rho\in X(H)$, and $j\in\{0,1\}$, we write for simplicity the generalized Eulerian Dirichlet series
\[A_{\ell,\beta,H,\chi,\rho,j}'(s)=\prod_{p\notin\Supp(H)}\left(1+2\sum_{\ve=\pm1}\ve^j\chi\left(\frac{p+\ve}{2}\right)\rho(p)|p+\ve|_\ell^\beta\left(\frac{p+\ve}{2}\right)^{-s}\right).\]
Writing
\[B_{\ell,\beta,H,\chi,\rho,j}'(s)=\prod_{p\notin\Supp(H)}\left(1+\frac{2}{p^s}\sum_{\ve=\pm1}\ve^j\chi\left(\frac{p+\ve}{2}\right)\rho(p)|p+\ve|_\ell^\beta\right),\]
we see that
\[A_{\ell,\beta,H,\chi,\rho,j}'(s)=B_{\ell,\beta,H,\chi,\rho,j}'(s)^{2^s}\cdot C_{\ell,\beta,H,\chi,\rho,j}'(s),\]
where $C_{\ell,\beta,H,\chi,\rho,j}'(s)$ is regular for $\Re(s)>1-\delta$ for some suitable $\delta>0$. Indeed, note that by our choice of $B$ and $H$ we may take, omitting subscripts,
\begin{align*}
    \log C'(s)&=\log A'(s)-2^s\log B'(s)\\
    &=\sum_{p\notin\Supp(H)}\sum_{n=1}^\infty\frac{(-1)^{n+1}}{n}\left[2\sum_{\ve=\pm1}\ve^j\chi\left(\frac{p+\ve}{2}\right)\rho(p)|p+\ve|_\ell^\beta\left(\frac{p+\ve}{2}\right)^{-s}\right]^n\\
    &\quad -2^s\sum_{p\notin\Supp(H)}\sum_{n=1}^\infty\frac{(-1)^{n+1}}{n}\left[\frac{2}{p^s}\sum_{\ve=\pm1}\ve^j\chi\left(\frac{p+\ve}{2}\right)\rho(p)|p+\ve|_\ell^\beta\right]^n.
\end{align*}
Since the sums $n\geq2$ are finite it suffices to consider contributions from $n=1$, which give
\[2\cdot 2^s\sum_{p\notin\Supp(H)}\sum_{\ve=\pm1}\ve^j\chi\left(\frac{p+\ve}{2}\right)\rho(p)|p+\ve|_\ell^\beta\left[\frac{1}{(p+\ve)^s}-\frac{1}{p^s}\right].\]
Since we have
\[|p+\ve|_\ell^\beta\left[\frac{1}{(p+\ve)^s}-\frac{1}{p^s}\right]=O\left(\frac{p^{\max(-\beta,0)}p^{\sigma-1}}{p^{2\sigma}}\right)\]
we get the desired claim that $C'(s)$ is regular for $\Re(s)= \sigma >1-\delta$.

Now, by application of Proposition \ref{pre1}, Proposition \ref{pre3}(1),(2), Corollary \ref{pre4} (which is applicable by our choice of $B$ and $H$), and Corollary \ref{pre5} to $B_{\ell,\beta,H,\chi,\rho,j}'(s)$, it follows that $A_{\ell,\beta,H,\chi,\rho,j}'(w)$ admits a holomorphic extension to $\Re(w)\geq1$ except for a singularity at $w=1$, where
\[(w-1)^{8\alpha_{\rho,\chi,j}}A_{\ell,\beta,H,\chi,\rho,j}'(w)=A+O(|w-1|^\theta),\]
for some $0<\theta\leq 1$ and some constant $A\neq0$, valid in a neighborhood of $w=1$. We write $[H]$ for the class $(D(H_p)/D(H))\in \prod_{p\in\Supp(H)}\Z_p/(\Z_p^*)^2$. Then the same analytic conclusion for the Dirichlet series
\[\begin{split}\sum_{\substack{G\in\Gen_H^{\varpi,v}\\\mu^*(G)_{\ell'}\equiv a\modulo \ell}}&\frac{|\mu^*(G)|_\ell^\beta}{\mu^*(G)^s}=\frac{1}{2}\frac{1}{|X(H)|}\sum_{j=0}^1\sum_{\chi\modulo\ell}\overline{\chi(a)}\\ &\times \sum_{\rho\in X(H)}\varpi^j\frac{c_\infty(H)|\mu^*(H)|_{\ell}^\beta\chi(D(H))\rho([H])}{\mu^*(H)^s}A_{\ell,\beta,H,\rho,j}'(s),\end{split}\]
which has nonnegative coefficients. The desired result then follows by application of Proposition \ref{prelim-7}.
\end{proof}

\begin{proposition}\label{dominator-v}
    There exists a constant $c_{\ell,\beta,B}'>0$ such that, for any packet $H\in\Packet$ such that $\{p:p\leq B\}\subseteq\Supp(H)$, and for any $X\geq2$,
    \[\frac{1}{X(\log X)^{\alpha-1}}\sum_{\substack{S\in\SGen_H',\mu^*(S)\leq X}}|\mu^*(S)|_\ell^\beta\leq c_{\ell,\beta,B}'2^{2+\omega_3(D(H))}\frac{|\mu^*(H)|_\ell^\beta}{\mu^*(H)}.\]
\end{proposition}

\begin{proof}
    Since $\sg(G)\leq 2^{2+\omega_3(D(H))}$ for $G\in\Gen_H'$, to prove the proposition it suffices to show the existence of $c_{\ell,\beta,B}'>0$ such that, for any packet $H\in\Packet$ such that $\{p:p\leq B\}\subseteq \Supp(H)$, and for any $X\geq2$,
    \[\frac{1}{X(\log X)^{\alpha-1}}\sum_{\substack{G\in\Gen_H',\mu^*(G)\leq X,}}|\mu^*(G)|_\ell^\beta\leq c_{\ell,\beta,B}'\frac{|\mu^*(H)|_\ell^\beta}{\mu^*(H)}.\]
    Consider now the generalized Dirichlet series
    \[\sum_{d}^\infty\frac{g(d)}{d^s}=\prod_{p>B}\left(1+2\sum_{\ve=\pm1}|p+\ve|_\ell^\beta\left(\frac{p+\ve}{2}\right)^{-s}\right).\]
    Then we see that
    \[\sum_{\substack{G\in\Gen_H',\mu^*(G)\leq X,}}|\mu^*(G)|_\ell^\beta\leq |\mu^*(H)|_\ell^\beta\sum_{d\leq X/\mu^*(H)}g(d).\]
    Arguing as in the proof of Proposition \ref{pointwise-v}, we deduce that \[\sum_{d\leq X}g(d)\sim c_{\ell,\beta,B} X(\log X)^{8\alpha-1},\]
    for some constant $c_{\ell,\beta,B}>0$. The result then follows.
\end{proof}

\begin{proposition}\label{pre2}
    Fix a prime $\ell\geq5$, and let $-1<\beta,\gamma\leq\infty$. Let $\mathcal C\subseteq(\Z/N\Z)^*$ be a set of invertible congruence classes modulo an even integer $N\in\N$. Consider the Dirichlet series
    \[D(s)=\sum_{d=1}^\infty\frac{f(d)}{d^s}=\prod_{p}\left(1+\frac{1}{p^s}\sum_{\ve=\pm1}\left|\frac{p+\ve}{2}\right|_\ell^\beta\left|\frac{p+\ve}{2p}\right|_\infty^\gamma\right).\]
    Then there exist constants $C,C'\in\R$ such that for any $X\geq2$
\begin{flalign*} \quad \rm{a}.&& \sum_{\substack{p\leq X\\ p\in\mathcal C}}\frac{f(p)}{p}=\frac{|\mathcal C|}{\varphi(N)}\alpha\log\log X+C+o(1),&&\end{flalign*}
and
\begin{flalign*} \quad \rm{b}.&&\sum_{p\in \mathcal{C}} \frac{f(p)}{p^{1+1/\log X}}=\frac{|\mathcal C|}{\varphi(N)}\alpha \log\log X+C'+o(1).&& \end{flalign*}
\end{proposition}

\begin{proof}
Part (b) follows from (a) by partial integration. The proof of (a) is as follows: if $p\equiv \alpha \pmod{N}$ represents a member of the class $\mathcal{C}$,  we require an asymptotic formula for \[\sum_{\substack{p\leq X\\ p\equiv \alpha\tmod{N}}}\frac{1}{p}\sum_{\ve =\pm 1}|p-\ve|_\ell^\beta.\] This can be written as
\[ \sum_{\ve =\pm 1}\left(\sum_{\substack{p\leq X\\ p\equiv \alpha\tmod{N}}}\frac{1}{p}- \sum_{\substack{p\leq X\\ p\equiv \alpha\tmod{N}\\ p\equiv \ve\tmod{\ell}}}\frac{1}{p}+(1-\ell^\beta)\sum_{t=2}^L\ell^{-t\beta}\sum_{\substack{p\leq X\\ p\equiv \alpha\tmod{N}\\ p\equiv \ve\tmod{\ell^t}}}\frac{1}{p}\right),\] with $L=[\frac{\log{X}}{\log{\ell}}]$.  The analysis now follows the same strategy as the proof of Proposition \ref{pre1} with Prop.~\ref{prelim-3} replacing Siegel-Walfitz. However, we make the following comment: a problem may arise if $q=[N,\ell^t]$ in Prop.~\ref{prelim-3} with $g^*(q,a)=\frac{1}{a}$ contributing to the third sum. If $t\leq \nu_\ell(N)$, then there are no issues. If $t> \nu_\ell(N)$, then one sees that $a\equiv \ve\tmod{\ell^t}$, so that if $a < \ell^t$ then $a=\pm 1$, so that $g^*(q,a)=0$. If $a>\ell^t$, then $g^*(q,a)=\frac{1}{a}$ is small and   there are again no issues.
\end{proof}

\begin{proposition}\label{splitting-v}
Fix a packet $H\in\Packet$ with $\{p:p\leq B\}\subseteq\Supp(H)$. Then for all small $\gamma>0$ we have
\[\sum_{\substack{G\in\Gen_H,\mu^*(G)\leq X\\\sg(G)\geq2}}|\mu^*(G)|_\ell^\beta\ll_{H,\gamma}\left\{\begin{array}{ll}X(\log X)^{2\alpha(1+\gamma)} & -1<\beta<0,\\ X(\log X)^{4(1+\gamma)} &\quad\quad\beta\geq 0.\end{array}\right.\]
\end{proposition}

\begin{proof}
Fix $\gamma>0$. It suffices to establish the result in the case $-1<\beta\leq0$. If $G\in\Gen_H'$ satisfies $\sg(G)\geq2$, this implies that there exists a proper subgroup $\mathcal C\leq(\Z/8r(H)\Z)^*$ of index $2$ such that $p\in\mathcal C$ for any $p\in T(G)$. Since there are at most finitely many subgroups $\mathcal C$, it suffices fix one such $\mathcal C$ and show that
\[\sum_{\substack{G\in\Gen_H,\mu^*(G)\leq X,\\p\in T(G)\Rightarrow p\in\mathcal C}}|\mu^*(G)|_\ell^\beta\ll_{\mathcal C,\gamma} X(\log X)^{2\alpha(1+\gamma)}.\]
We shall use Rankin's method. Let $A(\epsilon)>0$ be such that
   \[\left(\frac{p}{p-1}\right)^{2}\leq 1+\gamma,\]
for all $p\geq A(\epsilon)$. For $\sigma=1/\log X>0$ and $X\geq 2$ we have
\begin{align*}
&\sum_{\substack{G\in\Gen_H,\\D(G)\leq X,\\p \in T(G)\Rightarrow p\in\mathcal C}}|\mu^*(G)|_\ell^\beta
\leq\sum_{\substack{G\in\Gen_H,\\p\in T(G)\Rightarrow p\in\mathcal C}}|\mu^*(G)|_\ell^\beta\left(\frac{X}{D(G)}\right)^{1+\sigma},\\
&\leq|\mu^*(H)|_\ell^\beta X^{1+\sigma}\prod_{p\in\mathcal C}\left(1+2\sum_{\ve=\pm1}|p+\ve|_\ell^\beta\left(\frac{p+\ve}{2}\right)^{-(1+\sigma)}\right),\\
&\leq |\mu^*(H)|_\ell^\beta X^{1+\sigma}\exp\left(2\cdot 2^{1+\sigma}\sum_{p\in\mathcal C}\sum_{\ve=\pm1}|p+\ve|_\ell^\beta(p+\ve)^{-(1+\sigma)}\right),\\
&\leq |\mu^*(H)|_\ell^\beta X^{1+\sigma}\exp\left(2\cdot  2^{1+\sigma}(1+\gamma)\sum_{p\in\mathcal C}\sum_{\ve=\pm1}|p+\ve|_\ell^\beta p^{-(1+\sigma)}\right),\\
&\ll|\mu^*(H)|_{\ell}^\beta X^{1+\sigma}\exp\left(2\cdot 2^{1+\sigma}(1+\gamma)\alpha\log\log X\right)\ll X(\log X)^{2\alpha(1+\gamma)},
\end{align*}
by Proposition \ref{pre2}.
\end{proof}

\begin{proposition}\label{ellone-v}
Given any packet $H$, define
\[N'(H)=2^{2+\omega_3(r(H))}\frac{|\mu^*(H)|_\ell^\beta}{\mu^*(H)}.\]
Then $\sum_{H\in\mathfrak R_B'} N'(H)<\infty$.
\end{proposition}
\begin{proof}
Let $Q'$ be the set of $\mu^*$-weighted root packets $H\in\mathfrak R'$ such that $\ell\notin\Supp(H)$. Then
\[\sum_{H\in\mathfrak R_B'}N(H)\leq 2^{4+B}\prod_{p\leq B}\left(\sum_{G\in\Gen_p}\frac{|\mu_p^*(G)|_\ell^\beta}{\mu_p^*(G)}\right)\cdot\sum_{H\in Q}N(H).\]
In light of Corollary \ref{prelim-7b}, it suffices to show that $\sum_{H\in Q'}N'(H)<\infty$.
But we have
\begin{align*}
\sum_{H\in Q'} N'(H)&\leq 4\prod_{p\neq \ell}\left(1+2\sum_{\substack{G\in\Gen_p\setminus\Gen_p^\tau}}\frac{|\mu_p^*(G)|_\ell^\beta}{\mu_p^*(G)}\right),\\
&\leq4\exp\left(2\sum_{p\neq\ell}\sum_{\substack{G\in\Gen_p\setminus\Gen_p^\tau}}\frac{|\mu_p^*(G)|_\ell^\beta}{\mu_p^*(G)}\right).
\end{align*}
Again by Cor. \ref{prelim-7b}, writing $\delta=\max(-\beta,0)\in(0,1)$, we have
\[\sum_{p\neq\ell}\sum_{\substack{G\in\Gen_p\setminus\Gen_p^\tau}}\frac{|\mu_p^*(G)|_\ell^\beta}{\mu_p^*(G)}\ll p^\delta\frac{2p-1}{p(p-1)^2}<\infty,\]
showing that $\sum_{H\in Q'}N'(H)<\infty$.
\end{proof}

We now restate and prove Theorem \ref{thm3n}.

\begin{thm}{\ref{thm3n}}
Let $\ell\geq5$ be prime and $-1<\beta\leq\infty$. Let $S$ run through the spinor genera of primitive integral ternary quadratic forms of positive determinant, and let $G$ run through the genera. There is an explicit constant $C_{\ell,\beta}'>0$ such that 
    \[\sum_{\substack{\mu^*(G)\leq X}}|\mu^*(G)|_\ell^\beta\sim \sum_{\substack{\mu^*(S)\leq X}}|\mu^*(S)|_\ell^\beta\sim C_{\ell,\beta}'~X(\log X)^{8\alpha(\ell,\beta) -1}.\]
Moreover for any $\varpi\in\{\pm1\}$ and $(a,\ell)=1$,
    \[\sum_{\substack{c_\infty(S)=\varpi,\,\mu^*(S)\leq X\\\mu^*(S)_{\ell'}\equiv a\tmod\ell}}|\mu^*(S)|_\ell^\beta\sim\frac{1}{2}\frac{1}{\ell -1}C_{\ell,\beta}'~X(\log X)^{8\alpha(\ell,\beta) -1}.\]
\end{thm}

\begin{proof}
Let $B>0$ be chosen so that the conclusion of Corollary \ref{pre4} holds with $c=2$.
For each $X\geq2$, we consider $f_X':\mathfrak R_B'\to\R_{\geq0}$ given by
    \[f_X'(H)=\frac{1}{X(\log X)^{8\alpha-1}}\sum_{\substack{S\in\operatorname{SGen}_H^\varpi,\; \mu^*(S)\leq X\\\mu^*(S)_{\ell'}\equiv a(\modulo \ell)}}|\mu^*(S)|_{\ell}^\beta,\]
    and
    \[f_X''(H)=\frac{1}{X(\log X)^{8\alpha-1}}\sum_{\substack{G\in\operatorname{Gen}_H^\varpi,\; \mu^*(G)\leq X\\\mu^*(G)_{\ell'}\equiv a(\modulo \ell)}}|\mu^*(G)|_{\ell}^\beta.\]
    so that what we need to show is
    \[\lim_{X\to\infty}\|f_X''\|=\lim_{X\to\infty}\|f_X'\|_1=\frac{1}{2}\frac{1}{\ell-1}C_{\ell,\beta}'\in\R_{>0},\]
    for some constant $C_{\ell,\beta}>0$. For each $X\geq2$, we define $N':\mathfrak R_B'\to\R_{\geq0}$ by
    \[N'(H)=2^{2+\omega_3(r(H))}\frac{|\mu^*(H)|_\ell^\beta}{\mu^*(H)}.\]
    Then we have $N'\in\ell^1(\mathfrak R_B')$ by Prop.~\ref{ellone-v}. By Proposition \ref{dominator-v}, we have bounds $f_X''(H)\leq f_X'(H)\leq c_{\ell,\beta}'\cdot N'(H)$
    for every $H\in\mathfrak R_B'$, where $c_{\ell,\beta}'>0$ is is the constant given in Prop.~\ref{dominator-v}. By Propositions \ref{pointwise-v} and \ref{splitting-v} we have
    \[\lim_{X\to\infty}f_X''(H)=\lim_{X\to\infty}f_X'(H)=\frac{1}{2}\frac{1}{\ell-1}C_{\ell,\beta,H},\]
    where the constant $C_{\ell,\beta,H}'$ is defined as in Proposition \ref{pointwise-v}. By the dominated convergence theorem, the pointwise limit $f':\mathfrak R_B'\to\R_{\geq0}$ given by $f'(H)=\lim_{X\to\infty} f_X'(H)$ is integrable, and we have
    \begin{align*}
    \lim_{X\to\infty}\|f_X'\|_1=\|f'\|_1=\sum_{H\in\mathfrak R_B}f'(H)=\frac{1}{2}\frac{1}{\ell-1}\sum_{H\in\mathfrak R_B} C_{\ell,\beta,H}'.
    \end{align*}
    Setting $C_{\ell,\beta}=\sum_{H\in\mathfrak R_B'}C_{\ell,\beta,H}'$, we obtain the desired result.
\end{proof}

\subsubsection{The squarefree case:}\label{sqfree-mu-remark}\ 

One can obtain a similar asymptotic formula for genera of squarefree determinant. The proof proceeds as in the proof of Theorem \ref{thm3n}, except that: (i) one is only summing over the finitely many $\mu^*$-weighted $B$-root packets $H\in\mathfrak R_B'$ of squarefree determinant (hence, no dominated convergence required), and (ii) instead of Eulerian Dirichlet series such as $A_{\ell,\beta,H,\chi,\rho,j}'(s)$ in Prop.~\ref{pointwise-v} one considers
\[A_{\ell,\beta,H,\chi,\rho,j}''(s)=\prod_{p\notin\Supp(H)}\left(1+\sum_{\ve=\pm1}\ve^j\chi\left(\frac{p+\ve}{2}\right)\rho(p)|p+\ve|_\ell^\beta\left(\frac{p+\ve}{2}\right)^{-s}\right)\]
et cetera (note the absence of the coefficient $2$ in the factor). This leads to the following result, whose proof we omit.
\begin{theorem}
Let $\ell\geq5$ be prime and $-1<\beta\leq\infty$. Let $G$ run through the genera of primitive integral ternary quadratic forms of positive, squarefree, determinant. There is an explicit constant $c_{\ell,\beta}>0$ such that, for any $\varpi\in\{\pm1\}$ and $(a,\ell)=1$,
    \[\sum_{\substack{\mu^*(G)\leq X,c_\infty(G)=\varpi,\\\mu^*(G)_{\ell'}\equiv a\modulo \ell\\\text{$D(G)$ squarefree}}}|\mu^*(G)|_\ell^\beta\sim \frac{1}{2}\frac{1}{\ell-1} c_{\ell,\beta}X(\log X)^{4\alpha(\ell,\beta)-1}.\]
\end{theorem}

\subsection{Proof of Theorem \ref{thm2n}}\label{pf-thm2n}\leavevmode \ 

Our purpose in this section is to compute the constant $C^*$ for the asymptotic formula in Thm.~\ref{thm2n}. The asymptotic follows from Thm.~\ref{thm3n} with $\beta=0$. For $j\in\{0,1\}$, we consider the sums 
\begin{equation}\label{Th12-a}
    T_j(X):= \sum_{\mu^*(G)\leq X} c^j_\infty(G),
\end{equation} so that 
\begin{equation}\label{Th12-b}
    S_{\mu^*}^\varpi(1,X):=\sideset{}{^\varpi}\sum_{\mu^*(G)\leq X} 1 = \half\left(T_0(X) +\varpi T_1(X)\right).
\end{equation} To study $T_j(X)$, we consider the generalized Dirichlet series 
\begin{equation}\label{Th12-c}U^{(j)}(w):= \sum_{G\in \text{Gen}} c^j_\infty(G)\mu^*(G)^{-w} .\end{equation} 

From the discussion in Sec.~\ref{sect:mass} and the 
 tables of \cite{GGSW}, if $ \nu_2(D(G))=2u+v$,
then
$
    \mu_2^*(G)
    =3m_2(G)
    =2^{u+v-6}\lambda,
$
where $\lambda\in\mathbb N$ and $3\leq \lambda\leq 2^4$, so that $ 2^6\mu^*(G)\in\mathbb Z$. We then  see that $U^{(j)}(w)$ converges absolutely for $\Re(w) = \sigma >1$ as follows:
\begin{align*}
    U^{(j)}(w)
    &=
    \sum_{G\in\Gen}
    \prod_{p\mid 2D(G)}c_p^j(G)\mu_p^*(G)^{-w}, \\
    &=
    \sum_d\sum_{G\in\Gen_d}
    c_2^j(G)\mu_2^*(G)^{-w}
    \prod_{\substack{p\mid d\\p\geq 3}}
    c_p^j(G)\mu_p^*(G)^{-w}.
\end{align*}
One can bound this by summing over all packets
\[
   U^{(j)}(w) \ll
    \sum_{r,s}
    \sum_{\substack{H\\D(H)=r}}\;\;
     \sideset{}{^\#}\sum_{\substack{L\\D(L)=s}}
    \mu_2^*(H)^{-\sigma}
    \prod_{\substack{p\mid r\\p\geq 3}}
    \mu_p^*(H)^{-\sigma}
    \prod_{p\mid s}\mu_p^*(L)^{-\sigma}.
\]
From the Prop.~\ref{prelim-6} and the tables in \cite{GGSW}, one has
\[
    \mu_p^*(L)\geq Cp,
    \qquad
    \mu_p^*(H)\geq C\sqrt p \; ,
\]
where $C$ is an absolute positive constant, and the number of $H$'s and
$L$'s for fixed $r$ and $s$ are $\ll (rs)^\eps$. Hence
\[
  U^{(j)}(w)  \ll
    \sum_r\bigl(r^{1/2-\eps}\bigr)^{-\sigma}
    \sum_s\bigl(s^{1-\eps}\bigr)^{-\sigma}
    <\infty,
\]
if $\sigma>1$.

\medskip
\begin{remark}\label{rem-1.2}
If $d=rs$ in the usual notation and $G\in\Packet_d$, then
\[
    (\sqrt r\,s)^{1-\eps}
    \ll \mu^*(G)
    \ll (rs)^{1+\eps}.
\]
Hence there are only finitely many $\mu^*(G)$'s for
$G\in\Packet_d$.\end{remark}

For any prime $p$ and $\alpha \geq 0$, let 
\[ A_p^{(j)}(w,\alpha):=\sum_{G\in \text{Gen}_p(p^{\alpha})}c_p^j(G)\mu_p^*(G)^{-w}.\]
Then, from  Prop.~\ref{prelim-6}, $A_p^{(j)}(w,0)=1$ if $p\nmid 2d$.
Let \[A^{(j)}(d,w)= \prod_{p}A_p^{(j)}(w,\nu_p(d)) = \prod_{p\mid 2d}A_p^{(j)}(w,\nu_p(d)).\] Then, $A^{(j)}(d,w)$ is Selberg multiplicative, so that there is a $b$ such that \[A^{(j)}(d,w)=A^{(j)}(b,w)B^{(j)}(\frac{d}{b},w),\] with $b|d$, and $B^{(j)}$ multiplicative (see Prop.~\ref{prelim-4} and the notation there). To compute $b$, we see that $A_p^{(j)}(w,0)=1$ except when $p=2$ so that $k_p=0$ if $p\geq 3$. From the tables, we have $A_2^{(j)}(w,0)= 4^w(3^{-w} + (-1)^j)\neq 0$ if $\sigma >0$. Hence we may take $b=1$ if $\sigma >0$, which we assume. By Prop.~\ref{ds-mass} we then have with $M=\mu^*$ and $f=c_\infty^j$, that
\begin{equation}\label{Thm12-d}
\begin{split}
U^{(j)}(w) = &:= U_{(2)}^{(j)}(w)U_{(3+)}^{(j)}(w) \\ & = \left(\sum_{k=0}^\infty A_2^{(j)}(w,k)\right)\prod_{p\geq 3}\left[1 + \sum_{k=1}^\infty A_p^{(j)}(w,k)\right] ,
\end{split}
\end{equation}say.
We analyse the case $p\geq 3$ first. From  Prop.~\ref{prelim-6}, by considering the various cases, one sees that if

\[A_p^{(j)}(w) := 2^w \left[ \left(1 + \frac{1}{p}\right)^{-w} + (-1)^j \left(1 - \frac{1}{p}\right)^{-w} \right],
\]
then 
\[A_p^{(j)}(w,1) = p^{-w}A_p^{(j)}(w),
\]
and
\[A_p^{(j)}(w,2)= A_p^{(j)}(w,1) + p^{-2w}A_p^{(0)}(w).\]
We write
\[1 + \sum_{k=1}^\infty A_p^{(j)}(w,k) = 1 + 2A_p^{(j)}(w,1) + B_p^{(j)}(w).\] Then
\[
B_p^{(0)}(w)=p^{-2w}\frac{ \left[ \left(A_p^{(0)}(w)\right)^2 - \left(A_p^{(1)}(w)\right)^2 + 2A_p^{(0)}(w)\left(1 - p^{-w}\right) \right]}{\left(1 - p^{-w}\right)^2}
,
\]
and
\[
  B_p^{(1)}(w) = 2p^{-2w}\frac{\left[ A_p^{(0)}(w) + A_p^{(1)}(w)p^{-w} \right]}{1 - p^{-2w}} 
  + p^{-4w}\frac{\left[ \left(A_p^{(0)}(w)\right)^2 - \left(A_p^{(1)}(w)\right)^2 \right]}{\left(1 - p^{-2w}\right)^2}\ .\]
so that $B_p^{(j)}(w) \leq c_1^\sigma  p^{-2\sigma}$ for an absolute constant $c_1 >0$ if say $\sigma >\half$. One also has
 \[
     B_p^{(0)}(1) = \frac{8p(3p - 1)}{(p-1)^3 (p+1)}
, \quad\text{and} \quad B_p^{(1)}(1) =  \frac{8(p^4 + 1)}{(p^2 - 1)^3}.
\]
 We have in \eqref{Thm12-d}
\begin{equation}\label{Thm12-e} U_{(3+)}^{(j)}(w) = \prod_{p\geq 3}\left[1 - A_p^{(j)}(w,1)\right]^{-2} \times V_1^{(j)}(w), \end{equation} where $V_1^{(j)}(w)$ 
converges absolutely for $\sigma > \frac{1}{2}$.

One has  \[ V_1^{(0)}(1) = \prod_{p\geq 3} \frac{(p^2 - 4p - 1)^2 (p^2 + 4p - 1)}{(p-1)^5 (p+1)}
 \neq 0,\]
and
\[V_1^{(1)}(1)= \prod_{p\geq 3}  = \frac{(p^2 + 3)^2 (p^8 - 12p^6 + 22p^4 - 28p^2 + 1)}{(p^2 - 1)^6}
 \neq 0.
\]

We next consider the Euler-product $V_0^{(j)}(w):= \prod_{p\geq 3}\left[1 - A_p^{(j)}(w,1)\right]^{-2}$ in \eqref{Thm12-e}. It is easy to check that there is a $\delta>0$ such that $A_p^{(j)}(w,1)\neq 1$ if $\sigma >1-\delta$, which we now assume.  We use the analysis from \cite{BEPS}, in particular Thm.5.2 to get

\begin{itemize}[leftmargin=2em]
    \item If $j=1$, then $V_0^{(1)}(w)$ is regular for
    $\sigma>1-\delta$.
\medskip
    \item If $j=0$,
    \[
        V_0^{(0)}(w)
        =\prod_{p\geq 3}
        \left[1-\frac{2^{w+1}}{p^w}\right]^{-2}Z_1(w)^2,
    \]
    where $Z_1(w)$ is regular for $\sigma>1-\delta$, with
    \[
        Z_1(1)
        =\prod_{p\geq 3}
        \left[
           1+  \frac{4}{p\left(p^2 - 4p - 1\right)}
        \right],
    \]
\end{itemize}
and
\[
    \prod_{p\geq 3}
    \left[1-\frac{2^{w+1}}{p^w}\right]^{-2}
    =\zeta(w)^{2^{w+2}}Z_2(w)^2,
\]
where $Z_2$ is regular for $\sigma>1-\delta$, with
\[
    Z_2(1)
    =-\prod_{p\geq 2}
    \left(1-\frac1p\right)^4
    \left(1-\frac4p\right)^{-1}.
\]
Thus $V_0^{(0)}(w) = \zeta(w)^{2^{w+2}}Z_1(w)^2Z_2(w)^2 $ is regular for $\sigma \geq 1$ except for the singularity at $w=1$. This gives us in \eqref{Thm12-e}
\[ U_{(3+)}^{(j)}(w) \sim (w-1)^{-8} \times \frac{1}{192\zeta(2)}\prod_{p\geq 3}\left(1-\frac{1}{p}\right)^4\left(1+\frac{4}{p}-\frac{1}{p^2}\right).\]

For the prime $p=2$ in \eqref{Thm12-d}, it is seen (from the tables in \cite{GGSW}), that the series $U_{(2)}^{(j)}(w)$ converges absolutely for $\sigma >0$, and moreover that $U_{(2)}^{(0)}(1)= \frac{2^37^2}{3}$. 

Combining the analysis above for all $p$'s, we therefore have that $U^{(1)}(w)$ is regular for $\sigma>1-\delta$, so that by Kronecker's Lemma\footnote{See the proof of Thm.~2 in the book by Hardy-Riesz ``The General Theory of Dirichlet's Series''.} \hide{(see \cite[Thm.~K]{PTBZ})} one has $T_1(X) = o(X)$ in \eqref{Th12-a}. 

The implications for $T_0(X)$ from $U^{(0)}(w)$ is obtained from the singularity of order 8 at $w=1$. An appeal to Lemma~5.1 and Thm.~5.1 of \cite{BEPS} (see Prop.~\ref{prelim-7} with $r=7$) shows that  $T_0(X) \sim \alpha X(\log X)^7$, with \[\alpha = \frac{49}{72\Gamma(8)}\frac{1}{\zeta(2)}\prod_{p\geq 3}\left(1-\frac{1}{p}\right)^4\left(1+\frac{4}{p}-\frac{1}{p^2}\right).\]

\subsubsection{The squarefree case:} \ 

If we were counting genera with squarefree determinants, the Euler product is $\prod_{p\geq 3}(1 + A_p^{(j)}(w,1))$, together with the $p=2$ factor. There is a  a pole of order 4 at $w=1$ when $j=0$, in which case the analog of the asymptotic for $T_0(X)$ is of the form $\widetilde{\alpha} X(\log X)^3$ where
 \[ \widetilde{\alpha} = \frac{15}{2\Gamma(4)}\zeta(2)\prod_{p\geq 3}\left(1-\frac{1}{p}\right)^4\left(1+\frac{4}{p}-\frac{1}{p^2}\right).\]

 \subsection{Consequences and refinements of Theorems \ref{thm1n} and \ref{thm3n}}\label{pf-thm4n} \ 

Here, we collect some consequences and refinements of our main results. We first consider  indefinite quadratic forms.

\begin{maintheorem}\label{thm4n}\
Let $\ell\geq5$ be  prime. Let $\beta,\gamma\in(-1,\infty]$ with $\min(\beta,0)+\min(\gamma,0)>-1$.
Let  $\alpha$ be as in \eqref{alpha-def}. There is a constant $K_{\ell,\beta,\gamma}>0$ such that, for all $(a,\ell)=1$,
\[\sum_{\substack{\text{$C$ indef class}\\ \chiorb(C)_{\ell'}\equiv a\modulo\ell\\D(C)\leq X\\}}|\chiorb(C)|_\ell^\beta|\chiorb(C)|_\infty^\gamma\sim K_{\ell,\beta,\gamma} X^{1+\gamma}(\log X)^{2^{1-\gamma}\alpha(\ell,\beta)-1}.\]
Similarly, there is a constant $K_{\ell,\beta}'>0$ such that, for all $(a,\ell)=1$,
\[\sum_{\substack{\text{$C$ indef class}\\\chiorb(C)_{\ell'}\equiv a\modulo\ell\\|\chiorb(C)|\leq X\\}}|\chiorb(C)|_\ell^\beta\sim K_{\ell,\beta}' X(\log X)^{8\alpha(\ell,\beta)-1}.\]
In particular, we have the following.
\begin{enumerate}[label=\upshape\arabic*.\ ]
\vspace{6pt}
\item \track{a} \quad\quad $\displaystyle \sum_{\substack{\text{$C$ indef class}\\ D(C)\leq X}}|\chiorb(C)|_\ell^\beta \sim  K_{\ell,\beta,0}\,X(\log X)^{2\alpha -1},\quad (\beta>-1). $
\vspace{6pt}
\item \track{a} \quad\quad $\displaystyle \sum_{\substack{\text{$C$ indef class}\\ D(C)\leq X}}|\chiorb(C)|_\infty^\gamma \sim  K_{\ell,0,\gamma}\,X^{1+\gamma}(\log X)^{2^{1-\gamma} -1},\quad (\gamma>-1). $
\vspace{6pt}
\item \track{b} \quad\quad $\displaystyle \sum_{\substack{\text{$C$ indef class}\\D(C)\leq X\\ \chiorb(C)\equiv a \tmod{\ell} }}|\chi_{\orb}(C)|_\infty^\gamma \sim K_{\ell,\infty,\gamma}\,X^{1+\gamma}(\log X)^{2^{1-\gamma}\left(1-\frac{1}{\ell-1}\right)-1}$,\: \\ \hspace*{\fill}$(\gamma>-1)$ . 
\vspace{6pt}
\item \track{d} \quad\quad $\displaystyle  \sum_{\substack{\text{$C$ indef class}\\ |\chiorb(C)|\leq X}}1 \sim  \frac{49}{4\Gamma(8)}\frac{1}{\zeta(2)}\prod_{p\geq 3}\left(1-\frac{1}{p}\right)^4\left(1+\frac{4}{p}-\frac{1}{p^2}\right)~X(\log X)^7 .  $
\vspace{6pt}
\end{enumerate}
\end{maintheorem}

We now discuss the case of definite ternary forms.

\begin{maintheorem}\label{thm5n}\
Let $\ell\geq5$ be  prime. Let $\beta,\gamma\in(-1,\infty]$ with $\min(\beta,0)+\min(\gamma,0)>-1$.
Let  $\alpha$ be as in \eqref{alpha-def}. There is a constant $K_{\ell,\beta,\gamma}>0$ such that, for all $(a,\ell)=1$,
\[\sum_{\substack{\text{$G$ def}\\h_{\orb}(G)_{\ell'}\equiv a\modulo\ell\\D(G)\leq X\\}}|h_\orb(G)|_\ell^\beta|h_\orb(G)|_\infty^\gamma\sim K_{\ell,\beta,\gamma} X^{1+\gamma}(\log X)^{2^{1-\gamma}\alpha(\ell,\beta)-1}.\]
Similarly, there is a constant $K_{\ell,\beta}'>0$ such that, for all $(a,\ell)=1$,
\[\sum_{\substack{\text{$G$ def}\\h_{\orb}(G)_{\ell'}\equiv a\modulo\ell\\|h_{\orb}(G)|\leq X\\}}|h_\orb(G)|_\ell^\beta\sim K_{\ell,\beta}' X(\log X)^{8\alpha(\ell,\beta)-1}.\]
In particular, we have the following.
\begin{enumerate}[label=\upshape\arabic*.\ ]
\vspace{6pt}
\item \track{a} \quad\quad $\displaystyle \sum_{\substack{\text{$G$ def}\\ D(G)\leq X}}|h_\orb(G)|_\ell^\beta \sim  K_{\ell,\beta,0}\,X(\log X)^{2\alpha -1},\quad (\beta>-1). $
\vspace{6pt}
\item \track{a} \quad\quad $\displaystyle \sum_{\substack{\text{$G$ def}\\ D(G)\leq X}}|h_\orb(G)|_\infty^\gamma \sim  K_{\ell,0,\gamma}\,X^{1+\gamma}(\log X)^{2^{1-\gamma} -1},\quad (\gamma>-1). $
\vspace{6pt}
\item \track{b} \quad\quad $\displaystyle \sum_{\substack{\text{$G$ def, }D(G)\leq X\\ h_\orb(G)\equiv a \tmod{\ell} }}|h_{\orb}(G)|_\infty^\gamma \sim K_{\ell,\infty,\gamma}\,X^{1+\gamma}(\log X)^{2^{1-\gamma}\left(1-\frac{1}{\ell-1}\right)-1}$,\quad \\ \hspace*{\fill}$(\gamma>-1)$ . 
\vspace{6pt}
\item \track{d} \quad\quad $\displaystyle  \sum_{\substack{\text{$G$ def}\\ |h_\orb(G)|\leq X}}1 \sim \frac{49}{4\Gamma(8)}\frac{1}{\zeta(2)}\prod_{p\geq 3}\left(1-\frac{1}{p}\right)^4\left(1+\frac{4}{p}-\frac{1}{p^2}\right)~X(\log X)^7 .  $
\vspace{6pt}
\end{enumerate}
Furthermore, we have
\[\sum_{\substack{\text{$G$ def}\\D(G)\leq X\\}}|h_\orb(G)|_\ell^\beta|h_\orb(G)|_\infty^\gamma\sim\sum_{\substack{\text{$G$ def}\\D(G)\leq X\\}}|h_\orb(G)|_\ell^\beta|h(G)|_\infty^\gamma.\]
\end{maintheorem}

\proof

Only the last statement needs proof. We first consider the case of  $d$  squarefree as the analysis is simpler. Let $h(d)$ be the number of classes of determinant $d$, while $h_0(d)$ the same but weighted. We have  (1).  $h(d) \asymp \phi(d)$ ; (2). the number of classes with non-trivial  automorphism group is about $d^{\half +\ve}$  (see Thm.~8.10 of \cite{ChanOh}), and so (3). $0\leq h(d)-h_0(d) \ll d^{\half +\ve}$, (note that $h(d) \asymp h_0(d)$). Also, since $d$ is squarefree,  $d^{1-\ve}\ll \mu^*(G) \ll d^{1+\ve}$. All sums below are restricted to squarefree $d$'s.

Suppose we know the asymptotics for $\sum_{D(G)\leq X}|\h_\orb(G)|^\beta_\ell \h_\orb(G)^\gamma$, and suppose it has size $\approx X^{\gamma+1}(\log{X})^\theta$, for some $\theta$. Then $\sum_{D(G)\leq X}|\h_\orb(G)|^\beta_\ell\h(G)^\gamma$ has the same asymptotics. 

To see this first note that $\h_\orb(G)\asymp \h(G)$ and \[|\h_\orb(G)^\gamma - \h(G)^\gamma|\ll \h_\orb(G)^{\gamma -1}(\h(G)-\h_\orb(G)).\] As $D(G)$ is squarefree, we have  $D(G)^{1 -\ve} \ll \h_\orb(G)\ll D(G)^{1 +\ve}$. Using  
\[\h(G)-\h_\orb(G) \ll \sum_{C\in G,|\Aut(C)|>2}1 \ll D(G)^{\half +\ve},\] 
we see that if $\gamma > -1$ 
\[ \begin{split}&\sum_{D(G)\leq X}|\h_\orb(G)|^\beta_\ell \h_\orb(G)^\gamma - \sum_{D(G)\leq X}|\h_\orb(G)|^\beta_\ell\h(G)^\gamma \\  &\quad\quad\ll \sum_{D(G)\leq X}|\h_\orb(G)|^\beta_\ell \h_\orb(G)^\gamma D(G)^{-\half + 2\ve}= o\Bigl( \sum_{D(G)\leq X}|\h_\orb(G)|^\beta_\ell \h_\orb(G)^{\gamma}\Bigr), \end{split}\] if $\ve>0$ is chosen small enough.

In a similar manner for squarefree determinants, one has
\[ \sum_{\h(G)\leq X}|\h_\orb(G)|^\beta_\ell  \sim \sum_{\h_\orb(G)\leq X}|\h_\orb(G)|^\beta_\ell. \] This follows from 
\[ \sum_{\h_\orb(G)\leq X - O(X^{\half +\ve})}  |\h_\orb(G)|^\beta_\ell \leq \sum_{\h(G)\leq X}|\h_\orb(G)|^\beta_\ell \leq \sum_{\h_\orb(G)\leq X}|\h_\orb(G)|^\beta_\ell. \]

 The same conclusion also holds without the restriction to squarefree $d$'s as follows. By reduction theory and the Siegel mass formula, one has  $\frac{d}{\log\log d} \ll \h(d) \ll d\log\log d$. Also,  the number of classes with non-trivial  automorphism group is between $d^{\half -\ve}$ and $d^{\half +\ve}$ (see \cite{Kit73}). Lastly, from Prop.\ref{prelim-6}, one has $(\sqrt{r_G}s_G)^{1-\ve}\ll \mu^*(G) \ll d^{1+\ve}$, where $d=D(G)=r_Gs_G$, with $s_G$ the odd squarefree-part of $d$. Then we have
\[ \h(G)^\gamma- \h_\orb(G)^\gamma \ll \h_\orb(G)^\gamma s_G^{-\half}D(G)^\ve,\] so that 
\[ \begin{split}\sum_{D(G)\leq X}|\h_\orb(G)|^\beta_\ell \h_\orb(G)^\gamma -& \sum_{D(G)\leq X}|\h_\orb(G)|^\beta_\ell\h(G)^\gamma \\  
&\ll X^\ve \sum_{D(G)\leq X}|\h_\orb(G)|^\beta_\ell \h_\orb(G)^\gamma s_G^{-\half}.\end{split} \]
To estimate the latter sum, we split the genera into the packet decomposition: \[ \begin{split}&\sum_{D(G)\leq X}|\h_\orb(G)|^\beta_\ell \h_\orb(G)^\gamma s_G^{-\half} \\ &\qquad\ll \sum_{D(H)\leq X}|\h_\orb(H)|^\beta_\ell \h_\orb(H)^\gamma \sum_{\substack{G\in \Gen_H \\ D(G)\leq X}}|\h_\orb(G)|^\beta_\ell \h_\orb(G)^\gamma s_G^{-\half},\\ &\qquad \ll X^{\half +\gamma+\ve}\sum_{\substack{ D(H)\leq X}}|\h_\orb(H)|^\beta_\ell \left(\frac{\h_\orb(H)}{D(H)}\right)^\gamma D(H)^{-\half}.\end{split}\]
To estimate the last sum, let $\delta = \frac{1}{10}(1-\beta^* -\alpha^*) >0$. Then
\[\begin{split} &\sum_{D(G)\leq X}|\h_\orb(G)|^\beta_\ell \h_\orb(G)^\gamma s_G^{-\half} \\&\qquad\ll  X^{1+\gamma-\delta +\ve}\sum_{ D(H)\leq X}|\h_\orb(H)|^\beta_\ell \left(\frac{\h_\orb(H)}{D(H)}\right)^\gamma D(H)^{-1+\delta}\ll X^{1+\gamma - \half \delta},\end{split}\]
the latter following as a conclusion from the analysis given in the proof of Corollary \ref{prelim-7a}. Thus, subject to the conditions of Theorem \ref{thm1n} on $\beta$ and $\gamma$, one may replace $|\h_\orb(G)|^\gamma_\infty$ with $|\h(G)|^\gamma_\infty$ in Theorems \ref{thm1n}-\ref{thm3n}.

\begin{appendix}

\addtocontents{toc}{\protect\setcounter{tocdepth}{3}}

\section{Most genera have one spinor genera}\label{appendixb}

The discussion below is for classes in genera of indefinite ternary quadratic forms, but the same discussion goes through for spinor genera in genera of definite ternary quadratic forms. For an integer $d \ge 1$ write $d = rs$, where $s$, the \textit{tame part}, is the product of the odd primes dividing $d$ exactly once, and $r = d/s$ is the \textit{deep part}. Let $\omega_3(r)$ be the number of odd primes $p$ with $p^3 \mid r$, and for a genus $G$ of primitive integral ternary quadratic forms let $h(G)$ be its number of classes. Set
\[
h^{(2+)}(d) = \sum_{\substack{G: \, D(G)=d \\ h(G) \ge 2}} h(G),
\]
the number of classes of determinant $d$, each lying in genera with more than one class.

\begin{proposition} For all $X \ge 2$,

\[\sum_{1 \le d \le X} h^{(2+)}(d) \ll X .\]
\end{proposition}

\begin{proof}\ 

\textit{The class group of a genus}:\quad 
By Prop.~4.4 of \cite{GGSW}, with $P = 8 \prod_{\text{odd } p, \, p^3 \mid d} p$ and $\text{Sq}^*( \mathbb{Z}/P\mathbb{Z} ) = ( \mathbb{Z}/P\mathbb{Z} )^{\times} \text{ modulo squares}$,
\[
\Cl(G) \simeq \text{Sq}^*( \mathbb{Z}/P\mathbb{Z} ) \,/\, \langle \text{data of } G_p \ (p \mid 2r); \text{ the tame primes } q \mid s \rangle.
\]
The packet of $G$ (i.e.\ local data at $2$ and at the deep primes) determines $r$, $P$, and the subgroup $NH$ generated by the first list, and we write $A_H = \text{Sq}^*( \mathbb{Z}/P\mathbb{Z} ) / NH$, a group of exponent $2$ and order at most $|\text{Sq}^*( \mathbb{Z}/P\mathbb{Z} )| = 2^{2+\omega_3(r)}$. Two facts follow from \cite[Sec. 4]{GGSW}: all genera with the same packet $H$ and the same determinant $rs$ have the same class number, namely the index
\[
i_H(rs) = [A_H : \langle \text{images of the } q \mid s \rangle],
\]
and the number of such genera is at most $2^{\omega(s)}$. Grouping the sum in the theorem by packets,
\begin{equation}\label{B01}
\sum_{d \le X} h^{(2+)}(d) \le \sum_{H} \sum_{\substack{s \le X/r_H \text{ odd squarefree} \\ (s, \, 2r_H) = 1}} 2^{\omega(s)} i_H(r_H s) \mathbf{1}\{i_H(r_H s) \ge 2\}.
\end{equation}

\noindent \textit{A character condition}:\quad  If $i_H \ge 2$, the images of the tame primes generate a proper subgroup $B < A_H$; the characters of $A_H$ trivial on $B$ are exactly the characters of the nontrivial group $A_H/B$, so some nontrivial character $\chi$ of $A_H$ kills every tame prime. Let $\chi$ be a real quadratic character of modulus dividing $P$, trivial on $NH$, with $\chi(q) = 1$ for all $q \mid s$; there are fewer than $2^{2+\omega_3(r)}$ such candidates. Bounding $i_H$ by $2^{2+\omega_3(r)}$ as well, we have
\begin{equation}\label{B02}
i_H \mathbf{1}\{i_H \ge 2\} \le 2^{2+\omega_3(r)} \sum_{\chi \neq 1} \prod_{q \mid s} \mathbf{1}\{\chi(q) = 1\},
\end{equation}
and the theorem reduces to bounding, for each such $\chi$ and $Y = X/r_H$,
\[
S_\chi(Y) = \sum_{\substack{s \le Y \text{ odd squarefree} \\ \chi(q) = 1 \ \forall \, q \mid s}} 2^{\omega(s)}.
\]

\noindent \textit{The mean-value bound}:\quad $S_\chi(Y)$ is the summatory function of the nonnegative multiplicative function given by $f(q) = 2 \cdot \mathbf{1}\{\chi(q) = 1\}$ on odd primes and $f = 0$ on higher prime powers. We apply the  Halberstam--Richert inequality (\cite{HH}, \cite[Thm.~III.3.5]{tenen}), which states that: 

\noindent for $f \ge 0$ multiplicative, with $f(p) \le A$ and $\sum_{p} \sum_{\nu \ge 2} f(p^\nu) p^{-\nu} \le B$, one has \[\sum_{n \le x} f(n) \le (A+B+1) \frac{x}{\log x} \exp \left( \sum_{p \le x} f(p)/p \right).\] This gives, with $A = 2$, $B = 0$,
\begin{equation}\label{B03}
S_\chi(Y) \le 3 \frac{Y}{\log Y} \exp \left( \sum_{q \le Y,\;\text{prime}} \frac{1 + \chi(q)}{q} \right).
\end{equation}
By Mertens' theorem, $\sum_{q \le Y} 1/q = \log \log Y + O(1)$; it remains to show
\begin{equation}\label{B04}
\exp \left( \sum_{q \le Y} \frac{\chi(q)}{q} \right) \ll \log(2P),
\end{equation}
 which follows by comparing the sum with the Euler product of $L(s_0, \chi)$ with  $s_0 = 1 + 1/\log Y$.

\vspace{0.5em}
\noindent \textit{Summing over characters and packets}:\quad By \eqref{B01}, \eqref{B02}, \eqref{B03}, and $\log(2P) \ll \log(2r)$, we have 
\[
\sum_{d \le X} h^{(2+)}(d) \ll \sum_{H} \frac{4^{2+\omega_3(r_H)} X}{r_H} \log(2r_H) \le X \sum_{r} E(r), ,
\]
where \[\quad E(r) = \frac{m(r) 4^{2+\omega_3(r)} \log(2r)}{r},\] with $m(r)$  the number of packets with deep part $r$, and $r$ running over the admissible deep parts: a power of $2$ times an odd square-full number. Writing $p^{k_p} \parallel r$, \cite[Sec.4]{GGSW} gives $m(r) \ll \prod_{p} (k_p + 1)^2 \le \tau(r)^2$. Since $\tau(r)^2$, $4^{\omega_3(r)}$ and $\log(2r)$ are each $O_\varepsilon(r^\varepsilon)$, we get $E(r) \ll_\varepsilon r^{\varepsilon-1}$.

For $\varepsilon < \frac{1}{2}$, the sum $\sum_{r} r^{\varepsilon-1}$ converges: the powers of $2$ contribute a geometric series, and every odd square-full number is of the form $u^2 v^3$, so that $\sum_{u,v \ge 1} u^{2(\varepsilon-1)} v^{3(\varepsilon-1)} < \infty$. Therefore $\sum_{d \le X} h^{(2+)}(d) \ll X$.
\end{proof}

\section{Background results}\label{prelim}

	We collect here various results that we use in the paper. Let $q$ and $a$ be integers such that 
	$1 \le a \le q$ and $(q,a) = 1$. For notation, see Notation \ref{notation}.
    
	\begin{lemma}\label{prelim-1}  For  $A>0$\leavevmode 
		\begin{enumerate}[label=\upshape\roman*.]
		\item Brun-Titchmarsh Inequality:
		\[\pi(X; q, a) \ll \frac{X}{\varphi(q) \log(3X/q)} \quad \text{ for }  q\leq 2x, \; X\geq 2.\]
		\item Siegel-Walfitz:
		\[ |E(X;q,a)|+ |E^*(X;q,a)| = O_A\left(X (\log{X})^{-A}\right),
		\]
        \item if $X\leq q$ and $a=\pm 1$, then $|E(X;q,a)|\ll \frac{Li(X)}{\phi(q)}$.
		\end{enumerate}
	\end{lemma}

    An arithmetical function $f$ is a Selberg \cite{selberg77} multiplicative function  if $f(n)=\prod_p F_p(\nu_p(n))$, with $F_p(0)=1$ for almost all primes $p$.

\begin{proposition}\label{prelim-4} A Selberg multiplicative function $f$ not identically zero is semi-multiplicative. That is there exists an integer $b\geq 1$ such that (i). $f(n)=0$ if $b\nmid n$, (ii). $f(b)\neq 0 $, and (iii). there exists a multiplicative function $\hat{f}$ such that $f(n)=f(b)\hat{f}(\frac{n}{b})$ for all $n$'s with $b\mid n$.
\end{proposition}
\proof Since $f$ is not identically zero, for each prime $p$ there is an integer $k_p$ such that (i). $k_p=0$ if $F_p(0)\neq 0$, and (ii). $k_p\geq 1$ is the least $j$ satisfying $0\leq j \leq k_p -1$ with $F_p(j)=0$ and $F_p(k_p)\neq 0$. We put $b= \prod_p p^{k_p}$, and $\hat{f}(l) = \frac{f(bl)}{f(b)}$.\qed

\begin{remark} Subject to the conditions in the proposition, 
\[\sum_{d\leq X} f(d) = f(b)\sum_{l\leq \frac{X}{b}} \hat{f}(l),\]
the analysis of which can proceed using Prop.~\ref{prelim-2} below when applicable.
\end{remark}
\medskip

    We now state some results from \cite{GranK} and \cite{BEPS}, and refer the reader to those papers for details and further references.
    
	\begin{proposition}[\cite{GranK}, Thm.~1 and Sec.~7]\label{prelim-2} \leavevmode \\
		  \indent Let $f$ be a multiplicative function satisfying 
		\begin{equation}\label{GK1} \sum_{p\le X} f(p)\log p = \alpha X +  O\left( \frac{X}{(\log X)^A} \right),\end{equation} for some $\alpha\in\C$ and some $A>0$, and either one of the two conditions stated below for some $k>0$. Let $J$ be the largest integer $<A$ and
		\[
		c_0=\tilde{c}_0 = \prod_p \left(1+\frac{f(p)}{p}+\frac{f(p^2)}{p^2}+\cdots\right)\left(1-\frac{1}{p}\right)^{\alpha} .
		\]
		Then
		\[\begin{split}
			\sum_{n\le X} f(n) 
			&= \int_2^X  \sum_{j=0}^J c_j  \frac{ (\log y)^{\alpha-j-1} } {\Gamma(\alpha-j)} dy
			+  O(X (\log X)^{k-1-A}(\log\log X)^{{\bf1}_{A=J+1}}),  \\
			&= X \sum_{j=0}^J \tilde{c}_j  \frac{ (\log X)^{\alpha-j-1} } {\Gamma(\alpha-j)}
			+  O(X (\log X)^{k-1-A}(\log\log X)^{{\bf1}_{A=J+1}})    \ ,
	\end{split}	\]
	where the coefficients $c_j$ and $\tilde{c}_j$ for $j\geq 1$ are given explicitly in \cite{GranK}.
		The implied constants depend at most on   $k$, $A$, and the implicit constant in \eqref{GK1}. The dependence on $A$ comes from both its size, and its distance from the nearest integer. The conditions to be satisfies are one of the follow two:
		\begin{enumerate}[label=\upshape{\arabic*}.]
			\item $|f|\le\tau_k$, the coefficient of $\zeta(s)^k$, for some positive real number $k$, or
			\item $f$ satisfies the following three conditions:
            \vspace{.2cm}
			\begin{enumerate}[label=\upshape\alph*.]
				\item $\displaystyle \quad \sum_{\substack{p\le X \\ \nu\ge1}} \frac{|f(p^\nu)|}{p^\nu} \le k\log\log X+O(1)
				\quad \text{ for all }  X\geq 2,$\vspace{.2cm}
				\item  $\displaystyle \quad \sum_{p\le X} \frac{|f(p)|\log p}{p} \ll \log X,$\vspace{.2cm}
				\item $\displaystyle \quad \sum_{p\le X,\ \nu\ge1} \frac{|f(p^\nu)|^2}{p^\nu} = o_{X\to\infty}(\log X).$
			\end{enumerate}
		\end{enumerate}
		\end{proposition}
\vspace{2pt}
    \begin{corollary}\label{abstract}
    Let $f:\Gen\to\C$ be a Selberg multiplicative mass function. Using the notation in Secs. 2 and 3, let $F(d)$ and $H(d)$ be the associated functions over $\mathbb{N}$, with parameter $b$ so that $F(d)=F(b)H(\frac{d}{b})$. Suppose there exist $\alpha_f\in\C$ and $A>0$ such that
        \[\sum_{p\leq X} H_p(1)\log p=\alpha_f X+O_A\left(\frac{X}{(\log X)^A}\right),\]
    and there is $k>0$ such that $|H_p(r)|\leq \tau_k(p^r)$ for all primes $p$ and all $r\geq1$. Then
    \[\lim_{X\to\infty}\frac{1}{X(\log X)^{\alpha_f-1}}\sum_{D(G)\leq X}f(G)=\frac{F(b)}{b\Gamma(\alpha_f)}\prod_p\left(\sum_{k=0}^\infty\frac{H_p(k)}{p^k}\right)\left(1-\frac{1}{p}\right)^{\alpha_f}\in\C.\]
\end{corollary}

The following is an extension of Prop.~\ref{prelim-2}.
\begin{proposition}\label{simple}
    Let $f:\N\to\C$ be a multiplicative function supported only on squarefree numbers. Assume that there exist some $\alpha\in\C$, some $k\in\N$, some $A>0$, and some small $\delta>0$ such that
    \begin{enumerate}[label=\upshape\arabic*.\ ]
        \item 
    $\displaystyle \quad \sum_{p\leq X}f(p)\log p=\alpha X+O\left(\frac{X}{(\log X)^A}\right),$ \vspace{0.2cm}
        \item 
        $\displaystyle \quad \sum_{p\leq X}\frac{|f(p)|}{p}\leq k\log\log X+O(1),$\vspace{0.2cm}
        \item
        $\displaystyle \quad \sum_{p\leq X}\frac{|f(p)|\log p}{p}\ll\log X, \ \text{and}$ \vspace{0.2cm}
        \item 
        $\displaystyle \quad \sum_{p}\sum_{k=2}^\infty\frac{|f(p)|^k}{p^{k(1-\delta)}}<\infty.$
    \end{enumerate}
    Then we have
    \[\lim_{X\to\infty}\frac{1}{X(\log X)^{\alpha-1}}\sum_{n\leq x}f(n)=\frac{1}{\Gamma(\alpha)}\prod_p\left(1+\frac{f(p)}{p}\right)\left(1-\frac{1}{p}\right)^\alpha.\]
\begin{proof}
    To deduce the result following the argument of \cite[Section 7]{GranK}, under our assumptions (1), (2), and (3),  it suffices to show that
    \[\sum_{n\leq x}|R_f(n)|\ll x^{1-\delta/2},\]
    where we define $R_f$ by $f=\tau_f*R_f$ and where $\tau_f$ is defined by the Euler product
    \[\sum_{n=1}^\infty\frac{\tau_f(n)}{n^s}=\prod_{p}\left(1-\frac{1}{p^s}\right)^{-f(p)}.\]
    In the case at hand, we have $R_f(p)=0$ and for $k\geq2$ we have
    \[R_f(p^k)=\frac{(f(p))_k}{k!}-f(p)\frac{(f(p))_{k-1}}{(k-1)!}.\]
    where $(a)_k=a(a+1)\cdots (a+k-1)$. Since we have $|(a)_k/k!|\leq \max(1,|a|)^k$ for every $a\in\C$ and $k\in\N$, we see that $|R_f(p^k)|\leq 2+2|f(p)|^k$. Our last hypothesis (4) implies that the Euler product for the Dirichlet series associated to $R_f$ converges absolutely for $\Re(s)\geq 1-\delta$. The   result then follows by Kronecker's lemma.
\end{proof}
\end{proposition}
\begin{proposition}[\cite{BEPS}, Lemma~5.1]\label{prelim-7}\leavevmode \\
Suppose $1\leq h_1< h_2 , \ldots $ and $h_n \to \infty$. Suppose the Dirichlet series $\sum_n a_nh_n^{-w}$ has non-negative coefficients and the series converges for $\Re(w)>1$ to $f(w)$. Suppose 
\[ \lim_{\substack{w \to 1 +iy\\\Re(w)>1}}f(w),\]
exists for each $y\neq 0$. And suppose there exist real numbers $A>0$, $r\geq 0$ and $0\leq \theta <1$ such that
\[ f(w) = A(w-1)^{-r-1} +O(|w-1|^{-r-\theta}),\]
for $\Re(w)>1$, and $|w-1|<1$. Then, as $X \to \infty$, we have
\[ \sum_{h_n\leq X}a_n \sim \frac{A}{\Gamma(r+1)}X(\log{X})^r.\]
\end{proposition}
	\medskip
    
	Next we state some results from \cite{bar22}.
	
	\begin{proposition}[\cite{bar22}, Prop.~3.1, Thm.~3.4, Thm.~5.2]\label{prelim-3}\leavevmode \\
		 \indent If  $X \ge 2$, $q \ge 1$, and $a \ge 1$  with 
		$(q,a) = 1$,
		\begin{enumerate}[label=\upshape{\arabic*}.]
		\item One has  \[ 
		\sum_{\substack{p \le X \\ p \equiv a \bmod q}} p^{-1}
		 = \phi(q)^{-1}\log\log X + g(q,a) +  F(X;q,a),\]
         where
         \[
		F(X;q,a)= x^{-1}  E(X;q,a) - \int_{X}^{\infty} u^{-2}E(u;q,a)\; du,
		 \]
		with $g(q,a) = g^*(g,a) + O(\frac{\log{q}}{\phi(q)})$, with $g^*(q,a) =\frac{1}{a}$ if $q\geq 2$ and $a$ a prime, and is zero otherwise.
		\item By Siegel-Walfitz, we have
		\[
		\sum_{\substack{p \le X \\ p \equiv a \bmod q}} p^{-1}
		=
		\phi(q)^{-1}\log\log X 
		+ g(q,a)
		+ O_A((\log X)^{-A}),
		\] the implicit constants depending only on $A$.
		\end{enumerate}
	\end{proposition}

\begin{corollary}\label{prelim-5}
\[ 
\sum_{\substack{p \le X \\ p \equiv a \bmod q}} \frac{\log{p}}{p}
= \phi(q)^{-1}\log X + h(q,a) +  F^*(X;q,a), \]
where
\[
F^*(X;q,a)= \frac{\log{X}}{X}  E(X;q,a) - \int_{X}^{\infty}\frac{\log{(u/e)}}{u^2}  E(u;q,a)\; du =  O_A((\log X)^{-A}),\]
with
\[ h(q,a) =   \int_{2}^{\infty}\frac{\log{(u/e)}}{u^2}  E(u;q,a)\; du - \frac{\log{2}}{\phi(q)}. \]
Here, $h(q,a) = h^*(g,a) + O_\ve(\frac{q^\ve}{\phi(q)})$ for any $\ve >0$, and $h^*(q,a) =\frac{\log{a}}{a}$ if $q\geq 2$ and $a$ a prime, and is zero otherwise. 
\end{corollary}
\begin{proof} The proof is along the lines given in \cite{bar22}, using Brun-Titchmarsh and Siegel-Walfitz Theorems. 
	\end{proof}

\end{appendix}



\hypersetup{urlcolor=black, citecolor=black}
\printbibliography
\end{document}